\documentclass{gAPA2e}
\usepackage{psfrag}          
\usepackage{color}           

\def\del  {\partial}

\def\eps{\epsilon}

\def\R{\mathbb{R}}

\begin{document}
\doi{10.1080/0003681YYxxxxxxxx}
 \issn{1563-504X}
\issnp{0003-6811}
\jvol{00} \jnum{00} \jyear{2008} \jmonth{January}

\markboth{Taylor \& Francis and I.T. Consultant}{Applicable Analysis}


\title{{\itshape Perturbative numeric approach in microwave imaging}}

\author{A. Rozanova-Pierrat$^{\rm a}$ $^{\ast}$\thanks{$^\ast$Corresponding author. Email: anna.rozanova-pierrat@polytechnique.edu
}  \\\vspace{6pt} $^{\rm a}${\em{Laboratoire de la Physique de la Matière Condensée,
École Polytechnique,
Route de Saclay,
91120 Palaiseau, France}} \\\vspace{6pt}\received{v3.3 released May 2008} }

\maketitle

\begin{abstract}
In this paper, we show that \textcolor{black}{using measurements for different frequencies,
and using} ultrasound localized perturbations it is possible to extend the
method of \textcolor{black}{the imaging by elastic deformation} developed by Ammari and al. \textcolor{black}{[Electrical Impedance Tomography by Elastic Deformation
SIAM J. Appl. Math. , 68(6), (2008), 1557--1573.]} to problems of the form
\begin{eqnarray*}
\operatorname{div}(\gamma\nabla u)+k^{2}qu & =0 & \mbox{ in }\Omega,\\
\gamma\frac{\partial u}{\partial n} & =\varphi & \mbox{ on }\partial\Omega,
\end{eqnarray*}
 and to reconstruct by a perturbation method both $\gamma$ and $q$, provided \textcolor{black}{that} $\gamma$ is
coercive and $k$ \textcolor{black}{is not} a resonant frequency.\bigskip

\begin{classcode}31B20; 31A25 \end{classcode}\bigskip

\end{abstract}

\section{Introduction \textcolor{black}{and Notations}}
In the recent years, a lot of attention has been devoted to the reconstruction
of physical parameters of partial differential equations from electromagnetic
measurements. In the case of electrical impedance tomography (EIT) it
is well known that the detection of the conductivity from boundary
measurements is a very ill-conditioned problem. This drawback has
limited its use so far to anomaly detection. In a recent work, Ammari
\textit{et al.~} have shown that combining these measurements with
simultaneous localized ultrasonic perturbations allows to recover
the conductivity with great precision. The purpose of this work is
to show that such an approach can be generalized successfully to the
study of Helmholtz type problems. 

\textcolor{black}{In what follows we use the following notations:}
\begin{itemize}
 \item \textcolor{black}{$\Omega$ is a smooth domain in $\R^n$ with a regular boundary denoted by $\del \Omega$,}
\item \textcolor{black}{$x$ is a point in $\Omega$,}
\item \textcolor{black}{$\Omega'=\{x\in \Omega| \operatorname{dist}(x,\del \Omega)\ge d_0>0\}$ represents the interior points of $\Omega$,}
\item \textcolor{black}{$w\subset \Omega'$ is the region of the localization of the ultrasound perturbations, which is supposed to be small  compared to the size of $\Omega'$,}
\item \textcolor{black}{$|w|$ is the volume of $w$, }
\item \textcolor{black}{$1_w$ denotes the characteristic function corresponding to
the set $w$, i.e., the function which takes the value $1$ on
the set and the value $0$ outside,}
\item \textcolor{black}{$z\in w$ is the centre point of the region of the ultrasound perturbation,}
 \item \textcolor{black}{$k_i>0$ is a frequency,}
\item  \textcolor{black}{ $\gamma(x)$ is the conductivity and is a scalar real-valued function such that~$0<~c_0<~\gamma(x)<C_0$ for all $x\in \overline{\Omega}$, }
\item  \textcolor{black}{$q(x)$ is the permittivity and is a scalar real-valued function such that~$0<~c_0<~q(x)<C_0$ for all $x\in \overline{\Omega}$,} 
\item  \textcolor{black}{$u(x)$ is the potential induced on the boundary by the electromagnetic field $\varphi$ in the absence of ultrasonic perturbations ($u(x)$  and $\varphi(x)$ are  complex-valued functions),} 
\item  \textcolor{black}{$u_w$ is the perturbed potential field induced on the boundary by the electromagnetic field $\varphi$ in the presence of ultrasonic perturbations localized in the domain $w$ ($u_w$ is a complex-valued function),}
\item  \textcolor{black}{$\lambda$ is the amplitude of the ultrasonic perturbation,}
\item  \textcolor{black}{$\gamma_w(x)$ is the perturbed conductivity (real-valued  positive bounded function),}
\item  \textcolor{black}{$\tilde{\gamma}$ is the value of the perturbed conductivity $\gamma_w$ in the area $w$ of the perturbation (real-valued positive bounded function),}
\item  \textcolor{black}{$q_w(x)$ is the perturbed permittivity (real-valued  positive bounded function),}
\item \textcolor{black}{$\tilde{q}$ is the value of the perturbed permittivity $q_w$ in the area $w$ of the perturbation (real-valued positive bounded function),}
\item  \textcolor{black}{$\mbox{{\sffamily{\textbf M}}}_w$ and $\mbox{{\sffamily{\textbf m}}}_w$ are the polarization tensors,}
\item  \textcolor{black}{$N_{\gamma,q}(x,z)$ is the Neumann function for the operator $\operatorname{div}(\gamma(x) \nabla_x)+q(x)$ in $\Omega$ corresponding to a Dirac mass at $z$,}
\item  \textcolor{black}{$W^1_\infty(\Omega)$ is the Sobolev space of the functions $u(x)$ such that $u\in L_\infty(\Omega)$ and $\nabla u \in L_\infty(\Omega)$,}
\item \textcolor{black}{for the complex-valued function $u$, the function $\overline{u}$ denotes its complex-conjugated.}
\end{itemize}

The problem we consider
is the following.
 Let $\gamma\in C^{1}(\Omega)$ and $q\in C^{0}(\Omega)$ \textcolor{black}{be bounded scalar real-valued functions (see the list of notations).} 
For $i=1,2$,
let $u_{i}\in H^{1}(\Omega)$ be such that 
\begin{eqnarray}
\operatorname{div}(\gamma\nabla u_{i})+k_{i}^{2}qu_{i} & =0 & \mbox{ in }\Omega,\label{eq:hhnp}\\
\gamma\frac{\partial u_{i}}{\partial\nu} & =\varphi_i & \textnormal{on }\partial\Omega.\label{nneiman} \end{eqnarray}
 The well-posedness of this problem requires that $k_{i}^{2}$ is
not an eigenvalue of the generalized eigenvalue problem 
\begin{eqnarray}
-\operatorname{div}(\gamma\nabla u) & =\lambda qu & \mbox{ in }\Omega,\label{eq:vp}\\
\gamma\frac{\partial u}{\partial\nu} & =0 & \textnormal{on }\partial\Omega.\nonumber 
\end{eqnarray}
It is well known that this problem admits a countable number of eigenmodes,
with no accumulation point, and that each eigenvalue as a finite multiplicity.
We will assume that $k_{1}$ and $k_{2}$ do not correspond to eigenvalues
of problem~(\ref{eq:vp}).  The generalization of \textcolor{black}{the} method introduced in~\cite{AmFink} is the following.
\textcolor{black}{For  frequency $k_{i}$ being fixed,} we measure the potential $u_{i}$,
solution of problem~(\ref{eq:hhnp})-(\ref{nneiman}), on $\partial\Omega$.

Assume now that ultrasonic waves are localized around a point $z\in\Omega$,
creating a local change in the physical parameters of the medium. \textcolor{black}{Further, we suppose that $q$ and $\gamma$ are known close
to the boundary of the domain, so that ultrasonic probing is limited
to interior points $x$ in $\Omega'$ (see the list of notations), where $d_0$ is very large compared to the radius of the spot of the ultrasonic perturbation. }

We suppose that this deformation affects $\gamma$ and $q$ \emph{linearly}
with respect to the amplitude of the ultrasonic signal. Such an assumption
is reasonable if the amplitude is not too large. Thus, when the electric
potential is measured while the ultrasonic perturbation is enforced,
the equation for the potential is
 \begin{eqnarray}
\operatorname{div}(\gamma_{w}\nabla u_{i,w})+k_{i}^{2}q_{w}u_{i,w} & = & 0\quad \mbox{ in }\quad \Omega,\label{eq:hhp}\\
\gamma\frac{\partial u_{i,w}}{\partial\nu} & = & \varphi_i\quad \textnormal{ on }\quad \partial\Omega.\label{neiman}
\end{eqnarray}
 with
\begin{eqnarray}
\gamma_{w} & = & \gamma+1_{w}(\tilde{\gamma}\lambda-\gamma),\label{gw}\\
q_{w} & = & q+1_{w}(\tilde{q}\lambda-q),\label{qw}
\end{eqnarray}
 where $\lambda$ is the amplitude of the ultrasonic perturbation \textcolor{black}{given by the ratio of the perturbed volume $V^p_w$ of $w$ over the unperturbed one $V_w$ (see~\cite{AmFink}). 
In
other words}


\begin{equation*}
\begin{array}{ll}
 \textcolor{black}{\gamma_w(x)= \left\{
	\begin{array}{ll}
		\gamma(x),
	&
		\quad x\in \Omega\setminus w,
	\\
                \lambda(x)\tilde{\gamma}(x),&
 		\quad x\in w
\end{array}\right.}
&
\textcolor{black}{\qquad q_w(x)=\left\{
	\begin{array}{ll}
		q(x),
	&
		\quad x\in \Omega\setminus w,
	\\
               \lambda(x)\tilde{q}(x),
	&
		\quad x\in w  
       \end{array}\right.}
\end{array}
 \end{equation*}

\textcolor{black}{where $\lambda(x)=V^p_w/V_w$ is a known function. }

The analysis of the change of the Neumann-to-Dirichlet map as a result
of electromagnetic perturbation of small volume follows~\cite{AmFink}. \textcolor{black}{The main differences between the case of the conductivity equation considered in~\cite{AmFink} and our  case of the Helmholtz equation are the following:  this time the boundary data $\varphi$ and the solutions $u_i$ are complex-valued functions in our case while they are real in~\cite{AmFink}) and in our case we need to reconstruct simultaneous two coupled real-valued parameters $\gamma$ and $q$. Therefore we expand  the main ideas of~\cite{AmFink} to our case (see Section~\ref{Dem}). The choice of real $\gamma$ and $q$ implies the existence of eigenfrequencies (see problem~(\ref{eq:vp})) and this gives an additional difficulty in numeric reconstruction. The case of complex $\gamma$ and $q$ which allows to avoid the resonances, will be considered in~\cite{AY}.} 

 The signature of the perturbations on boundary measurements
can be measured by the change of energy on the boundary, namely 
\begin{equation}\label{gee}
 \int_{\partial\Omega}(u_{w}-u)\overline{\varphi}d\sigma=|w|\left[\mbox{{\sffamily{\textbf M}}}_w\left(\frac{\tilde{\gamma}\lambda}{\gamma}\right)(\tilde{\gamma}\lambda-\gamma)\nabla u(z)\cdot\overline{\nabla u(z)}-k^{2}(\tilde{q}\lambda-q)u(z)\cdot\overline{u(z)}\right].
\end{equation}

 Assuming the perturbed region is a ball, the polarization tensor \textcolor{black}{$\mbox{{\sffamily{\textbf M}}}_w\left(\frac{\tilde{\gamma}\lambda}{\gamma}\right)$}
is a scalar, 
$$M_w=(\tilde{\gamma}\lambda-\gamma)/(\tilde{\gamma}\lambda+\gamma).$$ 
Therefore, for a localized
perturbation focused at a point $z$, we read the following data (rescaled
by the volume)
 \begin{equation}
D_{z}(\lambda)=\gamma|\nabla u(z)|^{2}\frac{(\frac{\tilde{\gamma}}{\gamma}\lambda-1)^2}{\frac{\tilde{\gamma}}{\gamma}\lambda+1}-k^{2}q|u(z)|^{2}(\frac{\tilde{q}}{q}\lambda-1).\label{Dzz}
\end{equation}
 \textcolor{black}{We notice  that the data $D_{z}(\lambda)$ from~(\ref{Dzz}) can be measured for given $\lambda$ and $k$ thanks to the identity:
$$D_{z}(\lambda)=\frac{1}{|w|}\int_{\partial\Omega}(u_{w}-u)\overline{\varphi}d\sigma.$$}

The parameters $\frac{\tilde{\gamma}}{\gamma}(z)$ and $\frac{\tilde{q}}{q}(z)$ are
unknown, but the amplitude $\lambda$ is known.
Varying the position of localization, we are  able to recover
this localized internal data everywhere inside the domain. Thanks to the following lemma~\cite{yveP,AY}, 
\begin{lemma}\label{P1}
If the data $D_{z}$ is known for four distinct values of $\lambda$,
chosen independently of $\gamma$ and $q$, then one can recover $
\gamma(z)|\nabla u(z)|^{2}\mbox{ and }q(z)|u(z)|^{2}.$
\end{lemma}
  we can find directly the functions $J(z)=\gamma(z) |\nabla u(z)|^{2}$ and $j(z)=q(z)|u(z)|^{2}$ for the unique solution $u$ of problem~(\ref{eq:hhnp})-(\ref{nneiman}). 

\textcolor{black}{The proof of  Lemma~\ref{P1} 
 is simply a study of functions of one variable, which is detailed
in  Appendix~\ref{App}.}

%

The rest of the paper is organized as follows:
in Section~\ref{Dem} we prove formula~(\ref{gee}),
 in Section~\ref{Algo} we describe a reconstruction method by perturbations and  in Section~\ref{num} we  give and analyse our numeric results, obtained for two different frequencies and one boundary data in the form of a plane wave. 

\section{Proof of  asymptotic expansion~(\ref{gee})}\label{Dem}

We suppose that $k^2$ do not correspond to eigenvalues of problem~(\ref{eq:vp}).
\textcolor{black}{To prove the asymptotic expansion~(\ref{gee}), we first need the following Proposition:}
\begin{proposition}\label{Prop}
 We have the following identities
\begin{eqnarray}
  &&\gamma\frac{\del(u_w-u)}{\del n}=0,
  \label{neiman0}\\
  && \operatorname{div}(\gamma_w\nabla( u_w-u))+k^2q_w(u_w-u)\nonumber\\
&&=-\operatorname{div}(1_w(\gamma_w-\gamma)\nabla u)-k^21_w(q_w-q)u,\label{spas}\\
  &&
   \operatorname{div}(\gamma\nabla(u_w-u))+k^2 q(u_w-u)\nonumber\\
&&=-\operatorname{div}(1_w(\gamma_w-\gamma)\nabla u_w)-k^2 1_w(q_w-q)u_w,\label{fe}
\end{eqnarray}
\end{proposition}
\textcolor{black}{Thanks to Proposition~\ref{Prop}, we can estimate the difference between the perturbed and unperturbed solutions $u_w-u$ in $L_2(\Omega)$ by a norm of $u$ in the perturbed region $w$ and by a power of the small volume $|w|$ bigger than $0.5$.}

 \begin{lemma}\label{LEMEST} Suppose that $\Omega\subset R^n$ contains a subset of
$\Omega'\subset \Omega$ of class $C^2$, such that $\operatorname{dist}
(\Omega',\del\Omega)>d_0>0$, and such that $w\subset \Omega'$. 
 Let
$q,$ $\gamma \in L_\infty(\Omega)$ be positive functions,
satisfying
$0<c_0<q(x),\gamma(x)<C_0<+\infty$  a. e.  $x\in \Omega$,
 and $k^2$ is not a Neumann eigenvalue for problem~(\ref{eq:vp}). 
 Then for the
functions $u_w$ and $u\in W^1_\infty$ verifying Eq.~(\ref{neiman0})
and Eq.~(\ref{spas})
we have
\begin{equation}\label{estim3}
\|u_w-u\|_{H^1(\Omega)}\le C |w|^{\frac 12}|u|_{W^1_\infty(w)}.
\end{equation}
Therefore, thanks to relation~(\ref{fe}),
 for $m=\max\{2,n\}$ and all  $\kappa$ satisfying
$0<\kappa<\frac{2}{m}$  there exists a positive constant $C>0$ depending
only on $\Omega'$, $d_0$, $c_0$, and $C_0$, such that
\begin{equation}\label{estim2}
\|u_w-u\|_{L_2(\Omega)}\le C |w|^{\frac 12
+\kappa}|u|_{W^1_\infty(w)}.
\end{equation}
\end{lemma}
 \textbf{Proof.} 
The proof of  estimate~(\ref{estim3}) follows the proofs of Lemma 15.1 and Proposition 15.2 from~\cite{AmKg}. Indeed,
as soon as $q(x)/\gamma(x)k^2$ is not an eigenvalue for the operator $-\triangle$ in $L_2(\Omega)$ with the homogeneous Neumann boundary condition, in our case  problem~(\ref{eq:hhnp})-(\ref{nneiman}) has a unique weak solution $u$ in $H^1(\Omega)$ (for every $\phi\in H^{-1/2}(\del \Omega)$).
$A_\delta$ is uniformly continuous and uniformly coercive on $H^1\times H^1$.
The embedding $H^1(\Omega)\Subset L_2(\Omega)$ is still compact because $\del \Omega \in C^2$ and $\Omega$ is compact.

For passing to the perturbed problem, we change $\delta$ on $w$ (repeat the procedure from~\cite{AmKg}) and obtain with the help of relation~(\ref{spas}) the desired estimate~(\ref{estim3}).

Let us prove estimate~(\ref{estim2}).
Select $v$ as the solution to
\begin{eqnarray*}
&& \operatorname{div}(\gamma \nabla v)+k^2 q v=u_w-u\\
&& \gamma \frac{\del v}{\del n}|_{\del \Omega}=0.
\end{eqnarray*}
For this $v$ we have $\|v\|_{H^2(\Omega)}\le
C\|u_w-u\|_{L_2(\Omega)}$, and 
\begin{eqnarray}
&&\int_\Omega |u_w-u|^2dx=-\int_\Omega \gamma(x) \nabla(u_w-u)\nabla \overline{v} dx + \int_\Omega\omega^2 q(x)(u_w-u)\overline{v}dx\nonumber\\
&&=\int_\Omega \overline{v}\left[ \operatorname{div} (\gamma(x)\nabla (u_w-u))+k^2 q(x)(u_w-u)\right]dx\nonumber\\
&&=-\int_\Omega \overline{v}\left[\operatorname{div} (1_w(\gamma_w-\gamma)\nabla u_w)+k^2 1_w(q_w-q)u_w\right]dx\nonumber\\
&&\le \left|\int_\Omega\overline{v}\left[\operatorname{div} (1_w(\gamma_w-\gamma)\nabla u_w)+k^2 1_w(q_w-q)u_w\right]dx\right| \nonumber\\
&&\le \left|\int_\Omega 1_w(\gamma_w-\gamma)\nabla\overline{v}\nabla u_wdx \right|+\left|\int_\Omega k^2 1_w(q_w-q)\overline{v} u_wdx\right|\nonumber\\
&&\le C\left(\left(\int_w|\nabla u_w|^qdx\right)^{\frac 1q}\left(\int_\Omega |\nabla v|^pdx\right)^{\frac 1p}+\left(\int_w| u_w|^{\tilde{q}}dx\right)^{\frac {1}{\tilde{q}}}\left(\int_\Omega | v|^{\tilde{p}}dx\right)^{\frac{ 1}{\tilde{p}}}\right)\nonumber\\
 && \le C_p\left(\left(\int_w |\nabla u_w|^{q}dx\right)^{\frac 1{q}} \|v\|_{H^2(\Omega)}+\left(\int_w |u_w|^{\tilde{q}}dx\right)^{\frac 1{\tilde{q}}} \|v\|_{H^1(\Omega)}\right)\nonumber\\
 &&\le C_p  \left(\left(\int_w |\nabla u_w|^{q}dx\right)^{\frac 1{q}}+\left(\int_w |u_w|^{\tilde{q}}dx\right)^{\frac 1{\tilde{q}}}\right)
 \|u_w-u\|_{L_2(\Omega)}\label{esti1}
\end{eqnarray}
 provided $p$, $\tilde{p}$ and $q$, $\tilde{q}$ are related by $\frac 1q +\frac 1p=1$ and $\frac{ 1}{\tilde{q}} +\frac{1}{\tilde{p}}=1$. We use Sobolev's Embedding Theorem to provide the inclusions
$H^2\subset W^1_p$ and $H^1\subset L_{\tilde{p}}$. We require that $q,$ $\tilde{q}>\frac{2m}{m+2},$  so that $1<p,\; \tilde{p}<\frac{2m}{m-2}$.
For any $1<\tilde{q}<2$ (see~\cite[p.164]{yveV}) we have
\begin{eqnarray}&&\|u_w\|_{L_{\tilde{q}}(w)}\le \|u_w-u\|_{L_{\tilde{q}}(w)}+\|u\|_{L_{\tilde{q}}(w)}\le\nonumber\\
&&\le \left(\int_w 1dx\right)^{\frac {1}{\tilde{q}} -\frac {1}{2}}\|u_w-u\|_{L_2(w)}+|w|^{\frac 1{\tilde{q}}}\|u\|_{L_\infty(w)}\nonumber\\
&& \le C |w|^{\frac 1{\tilde{q}}}\|u\|_{L_\infty(w)},\label{estin}
\end{eqnarray}
and for any $1<q<2$ we obtain
\begin{equation}\label{yv}
 \|\nabla u_w\|_{L_q(w)}\le C_q|w|^{\frac{1}{q}}\|\nabla u\|_{L_\infty}.
\end{equation}
A combination of estimations~(\ref{esti1}),~(\ref{yv}) and~(\ref{estin}) yields
\[\|u_w-u\|_{L_2(\Omega)}\le C|w|^{\frac 1q} (|u|_{L_\infty(w)}+|\nabla u|_{L_\infty(w)})\]
for any $\frac {2m}{m+2}<q<2$. In other words\textcolor{black}{,} $\frac 12<\frac 1q<
\frac{m+2}{2m}=\frac 12+\frac 2m$ with $m=\max\{2,n\}$, from where
we can take $\tilde{q}=\frac 1q=\frac 12+\kappa$ for $0<\kappa<\frac
2m$. $\Box$

\textcolor{black}{In addition of estimates~(\ref{estim3}) and~(\ref{estim2}), let us show that the difference  $u-u_w$ can be totally described  by an integral expression over $w$.
}

\begin{proposition}
Suppose that $k^2$ is not the Neumann eigenvalue for
$\operatorname{div}(\gamma(x)\nabla_x )+q(x)$ on $w$. Let $N_{\gamma q}(x,z)$ be the Neumann function for
$\operatorname{div}(\gamma(x)\nabla_x )+q(x)$ in $\Omega$ corresponding to a Dirac
mass at $z$. That is $N_{\gamma q}$ is the solution to
\begin{equation}\left\{%
\begin{array}{ll}
   \operatorname{div}(\gamma(x)\nabla_x N_{\gamma q}(x,z))+k^2 q(x) N_{\gamma q}(x,z)=-\delta_z, & \hbox{in } \Omega, \\
    \gamma\frac{\del N_{\gamma q}}{\del \nu}=0 & \hbox{ on } \del \Omega. \\
\end{array}%
\right.\label{NN}
\end{equation}
Then, by definition of $N_{\gamma q}$ (which is a real function!), the function $U$ defined by
\[U(x)=\int_{\del \Omega} N_{\gamma q}(x,z)\varphi(z)d\sigma(z)\]
is the solution of  system~(\ref{eq:hhnp})-(\ref{nneiman}). 
 Therefore,
the solutions $u$ and $u_w$ of  systems~(\ref{eq:hhnp})-(\ref{nneiman}) and~(\ref{eq:hhp})-(\ref{neiman}) satisfy
\begin{equation}
(u-u_w)(z)=\int_{w}(\gamma_w-\gamma)(z)\nabla u_w(z) \nabla_z N_{\gamma q}(z,x)dz+\int_w k^2(q-q_w)(z)u_w(z) N_{\gamma q}(z,x)dz. \label{asim}
\end{equation}
\end{proposition}
 \textbf{Proof.}
Note that the Neumann function $N_{\gamma q}(x,z)$  is
defined as a function of $x\in \overline{\Omega}$ for each fixed
$z\in \Omega$. Since $k^2$ is not the Neumann
eigenvalue for $\operatorname{div}(\gamma(x)\nabla_x )+q(x)$ on $w$, the direct problem~(\ref{eq:hhnp})
admits a unique solution $u$ \textcolor{black}{(see~\cite{AmKg})}. Thus, the solution $u$
is represented by the formula
\[u(x)=\int_{\del \Omega} N_{\gamma q}(x,z)\varphi(z)d\sigma(z).\]
We notice  that
\begin{equation}\label{uwt}
\operatorname{div}(\gamma(x)\nabla u_w)+k^2q(x)u_w=-\operatorname{div}(1_w(\gamma_w-\gamma)(x)\nabla u_w)-k^21_w(q_w-q)(x)u_w,
\end{equation}
We multiply relation~(\ref{uwt}) by $N_{\gamma q}$  and integrate \textcolor{black}{over $\Omega$:}
\begin{eqnarray*}
&&\int_{\del \Omega} \varphi(z)N_{\gamma q}(x,z)d\sigma(z)-\int_\Omega
\gamma(z)\nabla u_w(z) \nabla N_{\gamma q}(x,z)dz+\int_\Omega k^2q(z)
u_w(z)N_{\gamma q}(z,x)dz=\\
&&=\int_\Omega 1_w(z)(\gamma_w-\gamma)(z)\nabla u_w(z) \nabla N_{\gamma q}(x,z)dz - \int_\Omega k^2 1_w(z)(q_w-q)(z)u_w(z)N_{\gamma q}(x,z)dz.
\end{eqnarray*}
\textcolor{black}{Therefore,} using $\int_{\del \Omega}  u_w(z) \gamma(z) \nabla
N_{q}(x,z)d\sigma(z)=0$, \textcolor{black}{from the following equality}
\begin{eqnarray*}
&&u(z)+\int_\Omega
u_w(z)\left(\operatorname{div}(\gamma(z)\nabla N_{\gamma q}(x,z))  + k^2q(z)N_{\gamma q}(z,x)\right)dz=\\
&&=\int_\Omega 1_w(z)(\gamma_w-\gamma)(z)\nabla u_w(z) \nabla N_{\gamma q}(x,z)dz +\int_\Omega k^2 1_w(z)(q-q_w)(z)u_w(z)N_{\gamma q}(x,z)dz,
\end{eqnarray*}
we obtain Eq.~(\ref{asim}).$\Box$

Multiplying Eq.~(\ref{asim})  by $\overline{\varphi}(x)$ and integrating over $\del \Omega$, we find
\begin{eqnarray*}&&\int_{\del \Omega}(u-u_w)\overline{\varphi} d\sigma(x)=\int_{w}(\gamma_w-\gamma)(z)\nabla u_w(z) \nabla_z\left(\int_{\del \Omega}N_{\gamma q}(z,x)\overline{\varphi}(x)d\sigma(x)\right)dz\\
&&+\int_{w}k^2(q-q_w)(z) u_w(z) \left(\int_{\del \Omega}N_{\gamma q}(z,x)\overline{\varphi}(x)d\sigma(x)\right)dz,
\end{eqnarray*}
which gives
\begin{equation}
 \int_{\del \Omega}(u-u_w)\overline{\varphi} d\sigma(x)=\int_{w}(\gamma_w-\gamma)(z)\nabla u_w(z) \nabla\overline{u}dz+\int_{w}k^2(q-q_w)(z)\ u_w(z) \overline{u}dz.\label{key}
\end{equation}

\begin{remark}\cite{yveV}
Consider a sequence of sets $w_\eps\subset\subset \Omega$.
 Since the family of functions
$\frac{1}{|w_\eps|}1_{w_\eps}$ is bounded in $L_1(\Omega)$, it
follows from a combination of the Banach-Alaoglu Theorem and Riesz
Representation Theorem that we may find a regular, positive Borel
measure $\mu$, and a subsequence $w_{\eps_k}$, with
$|w_{\eps_k}|\rightarrow 0$, such that
$\frac{1}{|w_\eps|}1_{w_\eps}\rightarrow d\mu$.
\end{remark}
\textcolor{black}{Finally, thanks to the a priori estimations~(\ref{estim3}),~(\ref{estim2}) and the representation formula~(\ref{key}),  we establish the main result:}
\begin{lemma}
Assume that $u\in W^1_\infty(w)$. Consider a sequence of sets
$w\subset\subset \Omega$ such that $\frac{1}{|w|}1_{w}$ converges in
the sense of measures to a probability measure $d\mu$ as $|w|$ tends
to zero. Then,
\begin{equation}\label{repres}
\int_{\del \Omega} (u_w-u)\overline{\varphi}
d\sigma=\int_w \mbox{{\sffamily{\textbf
M}}}_w \textcolor{black}{(\tilde{\gamma}\lambda-\gamma)}|\nabla u|^2 dx-k^2|w|\int_{w} \mbox{{\sffamily{\textbf m}}}_w \textcolor{black}{(\tilde{q}\lambda-q)} |u|^2dx+O(|w|^{1+\kappa}).
\end{equation}
The exponent $\kappa$ only depends on $\Omega_1$, $\sup_\Omega
|q_w|$, $\sup_\Omega
|\gamma_w|$, $\inf_\Omega |q_w|$ and $\inf_\Omega |\gamma_w|$. The remainder term has the form
\[|O(|w|^{1+\kappa})|\le C|w|^{1+\kappa}\|u\|_{W^1_\infty(w)}\|\nabla \psi\|_{L_\infty(\Omega)},\]
where $C$  depends only on $\Omega_1$, $\sup_\Omega
|q_w|$, $\sup_\Omega
|\gamma_w|$, $\inf_\Omega |q_w|$ and $\inf_\Omega |\gamma_w|$. Finally, with a hypothesis that $w$ is a ball, \textcolor{black}{the polarization tensors $ \mbox{{\sffamily{\textbf
M}}}_w$ and $\mbox{{\sffamily{\textbf
m}}}_w$ become} the scalar functions $M_w$ and $m_w$\textcolor{black}{, which are
given by}
\[ M_w=\frac{1}{|w|}1_w(x)\textcolor{black}{\left(\frac{\frac{\tilde{\gamma}}{\gamma}\lambda(x)-1}{\frac{\tilde{\gamma}}{\gamma}\lambda(x)+1}\right)}\quad \hbox{and}\quad  m_w=\frac{1}{|w|}1_w(x).\]
\end{lemma}
\textbf{Proof.} Suppose that $k^2$ is not a Neumann eigenvalue for problem~(\ref{eq:vp}).
We have relation~(\ref{key}). We are looking for an approximation of the terms of Eq.~(\ref{key}) depending on $u_w$ by a function depending on $u$.
In the same way as in~\cite{AmFink}, we introduce the solution $\zeta_w$ of the following problem
\begin{equation*}
\left\{\begin{array}{l}
  \operatorname{div} (\gamma_w(x)\nabla\zeta_w)+k^2q_w(x)\zeta_w=\operatorname{div} (\gamma(x)\nabla_x \zeta)+k^2q(x) \zeta \quad \hbox{in }\Omega,\\
\gamma\frac{\del \zeta_w}{\del n}=\gamma \frac{\del \zeta}{\del n} \quad \hbox{on } \del \Omega. 
\end{array}
\right.
\end{equation*}
Corresponding to $\zeta_w$, we define in the unperturbed case $\zeta=x+C+i(x+\tilde{C})$, where $C$ and $\tilde{C}$ are  constants in $\R^d$ for $x\in \R^d$.
This time all functions $\zeta$, $\zeta_w$, $u$ and $u_w$ are complex.
The choice of $C$ ($\tilde{C}$) will be  discussed later. 
Thanks to Lemma~\ref{LEMEST}, for $\zeta_w-\zeta$ we still have an analogue version of  Proposition 3.1 of~\cite[p.6]{AmFink}:
\begin{proposition}
 Consider a sequence of sets $w\subset \subset \Omega,$ such that $\frac{1}{|w|}1_w$ converges in the sense of measures to a probability measure $d\mu$ as $|w|$ tends to zero. Then, the corrector $\frac{1}{|w|}1_w \frac{\del \zeta_w}{\del x_j}$ converges in the sense of measures to $M_jd\mu$ ($M_j$ is a scalar function).
Furthermore, it satisfies
\[\|\nabla (\zeta_w-\zeta)\|_{L_2(\Omega)}\le C|w|^{\frac{1}{2}} \quad \hbox{and } \|\zeta_w-\zeta\|_{L_2(\Omega)}\le C |w|^{\frac{1}{2}+\kappa},\]
where the constants $\kappa>0$ and $C>0$ depend only on $\Omega_1$, $\sup_\Omega
|q_w|$, $\sup_\Omega
|\gamma_w|$, $\inf_\Omega |q_w|$ and $\inf_\Omega |\gamma_w|$. 
\end{proposition}

The rest of the proof follows the  analogous one given in details in~\cite{AmFink}.
This time the remaining term is bounded by
 \[|O(|w|^{1+\kappa})|\le C|w|^{1+\kappa}\|u\|_{W^1_\infty(w)}\|\psi\|_{W^1_\infty(\Omega)}.\]
We also remark (see~\cite{AmFink} for the notations) that the choice of $\psi_i=\frac{\del}{\del x_i}u\star \eta$ (where $\eta$ is the standard mollifier) determine the constants $C=(C_1,\dots,C_d)$ and $\tilde{C}=(\tilde{C}_1,\dots,\tilde{C}_d)$ in the definition of the function $\zeta(z)$:
\[C_j+i\tilde{C}_j=\frac{\overline{u}(z_0)-z_0\frac{\del} {\del x_i}\overline{u}(z_0)}{\frac{\del} {\del x_i}\overline{u}(z_0)} \quad \hbox{for a } z_0\in w,\]
which ensures that $\zeta(z)\overline{\psi} \approx \overline{u}(z)$.
 
Finally, we deduce 
\begin{eqnarray*}
 &&\int_{\del \Omega}(u-u_w)\overline{\varphi} d\sigma(x)\\
&&=|w|\int_{w}\mbox{{\sffamily{\textbf M}}}_w(\gamma_w-\gamma)(z)|\nabla u(z)|^2 dz+|w|\int_{w}k^2\frac{1_w}{|w|}(q-q_w)(z)| u(z)|^2dz +O(|w|^{1+\kappa}),
\end{eqnarray*}
 with $M_w=\frac{1}{|w|}1_w(x)\textcolor{black}{\left(\frac{\frac{\tilde{\gamma}}{\gamma}\lambda(x)-1}{\frac{\tilde{\gamma}}{\gamma}\lambda(x)+1}\right)}$  if $w$ is a sphere.
This proves relation~(\ref{gee}) and 
 provide the existence of a known function $D_z(\lambda)$ \textcolor{black}{from~(\ref{Dzz}).} $\Box$


\section{Reconstruction $\gamma$ and $q$ by a perturbative method. Numeric algorithm}\label{Algo}

We consider the system of Helmholtz equations with different frequencies $k_1\ne k_2$:
\begin{eqnarray}
 &&\operatorname{div}(\gamma(x)\nabla u_{k_1})+k_1^2 q(x)u_{k_1}=0 \qquad \hbox{in } \Omega,\\
&&\operatorname{div}(\gamma(x)\nabla u_{k_2})+k_2^2 q(x)u_{k_2}=0 \qquad \hbox{in } \Omega,\\
  && u_{k_1}= u_{k_2}=\psi \qquad \hbox{on } \partial\Omega.\label{neiman}
\end{eqnarray}
The data $\psi$ is the Dirichlet data measured as a response to the
current $\varphi$ in  absence of elastic deformation. We take
$\psi=e^{i\arctan{y/x}}$, which represents \textcolor{black}{a} plane wave. 

We use the following formulas
\begin{tabular}{ll}
$\gamma(x) |\nabla u_{k_1} |^2=J_{k_1}(x)$ and & $q(x)|u_{k_2}|^2=j_{k_2}(x).$
  \end{tabular}
Thus, we can approximate our problem by  system~(\ref{AQJ0}) and~(\ref{AQJ})
\begin{eqnarray}
 &&\operatorname{div}\left(\frac{J_{k_1}(x)}{|\nabla u_{k_1}|^2}\nabla u_{k_1}\right)+k_1^2 q(x)u_{k_1}=0 \qquad \hbox{in } \Omega,\label{AQJ0}
\\
&&\operatorname{div}(\gamma(x)\nabla u_{k_2})+k_2^2 \frac{j_{k_2}(x)}{|u_{k_2}|^2}u_{k_2}=0 \qquad \hbox{in } \Omega, \label{AQJ}
\end{eqnarray}
\textcolor{black}{where it is supposed that 
$$|\nabla u_{k_1}|^2>0 \quad \hbox{and} \quad  |u_{k_2}|^2>0 \quad \hbox{for all} \quad x\in \Omega.$$
}
Let us explain the steps of the numeric algorithm. The method uses two sub-algorithms to reconstruct $\gamma$ for a fixed $q$ (constant for the ultrasound perturbation)  and to reconstruct $q$ for a fixed $\gamma$ (constant).

First we notice that we have two frequencies {$k_1$} and {$k_2$}.

{\color{black}\textbf{Step 0.} We construct the functions {$J_{k_1}$ and $j_{k_2}$}.

\textbf{Step 1.} We take an initial guess $q_0$ and $\gamma_0$.

\textbf{Step 2.}  In the aim of updating first $\gamma_0$ we solve the linear system for chosen $q_0$ and $\gamma_0$ and the frequency $k_1$:
\[\begin{array}{lr}
\left\{\begin{array}{l}
\operatorname{div}(\gamma_0\nabla u_{k_1})+k_1^2q_0 u_{k_1}=0\\
 {u_{k_1}}|_{\del \Omega}=\psi
\end{array}\right.
\end{array}\]
We obtain the solution of this system which we denote by ${u_0}_{k_1}$.
Knowing the approximate solution ${u_0}_{k_1}$, we calculate the error on $\gamma$:
$${E_0}_{k_1}=\frac{J_{k_1}}{|\nabla {u_0}_{k_1}|^2}-\gamma_0.$$

\textbf{Step 3.} We verify the condition $|{E_0}_{k_1}|<\eps_{\textrm{precision}}$ for a given positive constant $\eps_{\textrm{precision}}$, which gives the desired order of the precision of the final result. 
If  $|{E_0}_{k_1}|$ is smaller than $\eps_{\textrm{precision}}$,  we  take  $\gamma\equiv \gamma_0$ and go to  Step~5 for the reconstruction of $q$, otherwise we go to  Step~4.

\textbf{Step 4.} We apply the algorithm described in details in Subsection~\ref{AlgoG} to determine the correctors  $ \delta\gamma_1$ and $\delta {u_1}_{k_1}$ for a fixed $q_0$ and  to update $\gamma_0$ using formula~(\ref{updateG}).

\textbf{Step 5.} In the aim of updating $q_0$, we solve the following linear system with  the frequency $k_2$ for a chosen $q_0$ and $\gamma_0$ updated on Step~4:
\[\begin{array}{lr}
\left\{\begin{array}{l}
\operatorname{div}({\gamma_0}\nabla {u_{k_2}})+{k_2^2q_0} {u_{k_2}}=0\\
 {u_{k_2}}|_{\del \Omega}=\psi
\end{array}\right.
\end{array}\]
We obtain the solution of this system which we denote by ${u_0}_{k_2}$.
Knowing the approximate solution ${u_0}_{k_2}$, we calculate the error on $q$:
$${e_0}_{k_2}=\frac{j_{k_2}}{|{u_0}_{k_2}|^2}-q_0.$$

\textbf{Step 6.} We verify the condition $|{e_0}_{k_2}|<\eps_{\textrm{precision}}$. 
If  $|{e_0}_{k_2}|$ is smaller than $\eps_{\textrm{precision}}$, we  take  $q\equiv q_0$ and finish the algorithm, otherwise we do  Step~7.

\textbf{Step 7.}
 We apply the algorithm described in details in Subsection~\ref{AlgoQ} to determinate the correctors  $\delta q_1$ and $\delta {u_1}_{k_2}$ for a fixed $\gamma_0$ and to update $q_0$ using formula~(\ref{updateQ}). Next we go to Step~2.
}

\subsection{Algorithm of reconstruction \textcolor{black}{of} $\gamma$ for a constant $q$}~\label{AlgoG}
\textbf{Step 1.} We start from  an initial guess $\gamma_0,$ and solve the corresponding Dirichlet problem for \textcolor{black}{the} Helmholtz equation
\begin{eqnarray*}
&& \operatorname{div}(\gamma_0(x) \nabla u_0) +k^2 q u_0=0,\\
&&u_0|_{\del \Omega}=\psi.
\end{eqnarray*}
Solving the direct problem for $\psi=e^{i\arctan{y/x}}$, we
obtain $u_0$.

\textbf{Step 2.} We have seen that our inverse problem is asymptotically approached by the  direct problem
\begin{equation}
    \left\{%
\begin{array}{ll}
    \operatorname{div}(\frac{J(x)}{|\nabla u|^2}\nabla u)+\omega^2qu=0, & \hbox{in } \Omega, \\
    u=\psi, & \hbox{on } \del \Omega. \\
\end{array}%
\right.
\end{equation}
We compute the difference
\begin{equation}\label{errG}
E_0:=\frac{J(x)}{|\nabla u_0|^2}-\gamma_0
\end{equation}
and verify
\begin{equation}\label{precisionG}
|E_0|<C_{\textrm{prec}},
\end{equation}
where $C_{\textrm{prec}}$ is our wished order of \textcolor{black}{the} precision. If 
condition~(\ref{precisionG}) holds, we finish our algorithm and \textcolor{black}{set}
$\gamma\equiv \gamma_0.$ \textcolor{black}{Otherwise} we \textcolor{black}{go to} the next step.

\textbf{Step 3.} We use now the expression
\[(\gamma_0+\delta \gamma_1)|\nabla(u_0+\delta u_1)|^2=J(x),\]
having \textcolor{black}{the} goal to approximate the known function $J(x)$ with the help
of the small correctors $\delta u_1$ and $\delta \gamma_1$. \textcolor{black}{We suppose that $\delta \ll 1$ and  that $\delta\max\limits_x|\gamma_1|$ and $\delta\max\limits_x |u_1|$ are of the order of $\delta$.}

\textcolor{black}{By expanding} the expression, we obtain
\begin{eqnarray*}
&&\delta
\gamma_1\left(1+2\delta\left(\frac{\nabla(\operatorname{Re}u_0)\nabla(\operatorname{Re}u_1)+\nabla(\operatorname{Im}u_0)\nabla(\operatorname{Im}u_1)}{|\nabla
u_0|^2}\right)+\delta^2\frac{|\nabla u_1|^2}{|\nabla u_0|^2}\right)
=\frac{J(x)}{|\nabla u_0|^2}-\gamma_0\\
&&
-2\delta\frac{\gamma_0\left(\nabla(\operatorname{Re}u_0)\nabla(\operatorname{Re}u_1)+\nabla(\operatorname{Im}u_0)\nabla(\operatorname{Im}u_1)\right)}{|\nabla
u_0|^2}-\delta^2\gamma_0\frac{|\nabla u_1|^2}{|\nabla u_0|^2}.
\end{eqnarray*}
We consider only terms of order \textcolor{black}{not} smaller than $\delta$:
\[\delta \gamma_1=\frac{J(x)}{|\nabla u_0|^2}-\gamma_0
-2\delta\frac{\gamma_0\left(\nabla(\operatorname{Re}u_0)\nabla(\operatorname{Re}u_1)+\nabla(\operatorname{Im}u_0)\nabla(\operatorname{Im}u_1)\right)}{|\nabla
u_0|^2}.\] \textcolor{black}{To find} the corrector $\tilde{u}_1=\delta u_1$, we
\textcolor{black}{expand} the following equation
\[ \operatorname{div}((\gamma_0+\delta\gamma_1 )\nabla(u_0+\delta u_1))+k^2q(u_0+\delta u_1)=0.\]

\textcolor{black}{By considering} the terms of order \textcolor{black}{not} smaller than $\delta$ and
\textcolor{black}{by} replacing   $\delta \gamma_1$ by the approximated formula, we can find $\tilde{u}_1$  as the solution of the following problem
\begin{eqnarray}
&&\operatorname{div}\left[\gamma_0\left(\nabla \tilde{u}_1-2
\frac{\nabla u_0}{|\nabla u_0|^2}\left(\nabla\operatorname{Re}u_0\nabla\operatorname{Re}\tilde{u}_1+\nabla\operatorname{Im}u_0\nabla\operatorname{Im}\tilde{u}_1\right)\right)\right]\label{Et}\\
&&+\operatorname{div}(E_0\nabla \tilde{u}_1)+\operatorname{div}(E_0\nabla u_0)+k^2q\tilde{u}_1=0,\nonumber\\
&& \tilde{u}_1|_{\del \Omega}=0.\label{Ett}
\end{eqnarray}
Let us \textcolor{black}{define} \[\bm{GU}_0=\left(\begin{array}{c}\nabla\operatorname{Re}u_0 \\
                            \nabla\operatorname{Im}u_0
                           \end{array}\right) \quad \hbox{and} \quad \bm{GU}_1=\left(\begin{array}{c}\nabla\operatorname{Re}\tilde{u}_1 \\
                            \nabla\operatorname{Im}\tilde{u}_1
                           \end{array}\right),\]
and suppose that \[\bm U_0=\left(\begin{array}{c}\operatorname{Re}u_0 \\
                            \operatorname{Im}u_0
                           \end{array}\right) \quad \hbox{and} \quad \bm U_1=\left(\begin{array}{c}\operatorname{Re}\tilde{u}_1 \\
                            \operatorname{Im}\tilde{u}_1
                           \end{array}\right),\]
thus we have
\[\nabla\operatorname{Re}u_0\nabla\operatorname{Re}\tilde{u}_1+\nabla\operatorname{Im}u_0\nabla\operatorname{Im}\tilde{u}_1=\bm{GU}_0\cdot \bm{GU}_1^T.\]
We also use the relation \[|\nabla u_0|^2=|\bm{GU}_0|^2.\]

We solve \textcolor{black}{problem~(\ref{Et})-(\ref{Ett}) for the real and imaginary parts of $\title{u}_1$} and using our notations we obtain the system
\begin{eqnarray*}
&&\operatorname{div}\left[\gamma_0\left(\bm{GU}_1-2\frac{\bm{GU}_0}{|\bm{GU}_0|}\left(\frac{\bm{GU}_0}{|\bm{GU}_0|}\cdot \bm{GU}_1\right)\right)\right]\\
&&+\operatorname{div}(E_0 \bm{GU}_1)+\operatorname{div}(E_0\bm{GU}_0)+k^2q\bm{U}_1=0,\\
&& \bm{U}_1|_{\del \Omega}=0.
\end{eqnarray*}
The vector $\bm{\theta}_0=\frac{\bm{GU}_0}{|\bm{GU}_0|}$ is \textcolor{black}{a} unit vector. We can rewrite our system in the form
\begin{eqnarray*}
&&\operatorname{div}\left[\gamma_0\left(\bf{Id}-2\bm{\theta_0}\otimes\bm{\theta_0}\right)\bm{GU}_1\right]+\operatorname{div}(E_0 \bm{GU_1})+\operatorname{div}(E_0\bm{GU}_0)+k^2q\bm{U}_1=0,\\
&& \hbox{ or using eigenvectors}\\
&&\operatorname{div}\left[\gamma_0\left(\bm{\theta}_0^{\perp}\otimes\bm{\theta}_0^{\perp}-\bm{\theta}_0\otimes\bm{\theta}_0\right)\bm{GU}_1\right]+\operatorname{div}(E_0 \bm{GU}_1)+\operatorname{div}(E_0\bm{GU}_0)+k^2q\bm{U}_1=0.
\end{eqnarray*}
We suppose that $\bm{GU}_1\parallel \bm{\theta}_0$ and obtain
 \begin{eqnarray}
&&-\operatorname{div}\left(\gamma_0\nabla \operatorname{Re}\tilde{u}_1\right)+\operatorname{div}(E_0\nabla \operatorname{Re} \tilde{u}_1)+
\operatorname{div}(E_0\nabla \operatorname{Re}u_0)+k^2q\operatorname{Re}\tilde{u}_1=0,\label{GU1cora}\\
&&\operatorname{Re}\tilde{u}_1|_{\del \Omega}=0;\nonumber\\
&&-\operatorname{div}\left(\gamma_0\nabla \operatorname{Im}\tilde{u}_1\right)+\operatorname{div}(E_0\nabla \operatorname{Im}\tilde{u}_1)+
\operatorname{div}(E_0\nabla \operatorname{Im}u_0)+k^2q\operatorname{Im}\tilde{u}_1=0,\label{GU1corb}\\
&& \operatorname{Im}\tilde{u}_1|_{\del \Omega}=0.\nonumber
\end{eqnarray}
This gives $\tilde{u}_1$. 

\textbf{Step 4.} We calculate
\begin{equation}
 \tilde{\gamma}=\gamma_0+\delta \gamma_1=\frac{1}{|\nabla u_0|^2}\left
(J(x)-2\gamma_0(\nabla\operatorname{Re}u_0\nabla\operatorname{Re}\nabla\tilde{u}_1+\nabla\operatorname{Im}u_0\nabla\operatorname{Im}\nabla\tilde{u}_1)\right).\label{updateG}
\end{equation}

We \textcolor{black}{set} now $\gamma_0\equiv \tilde{\gamma}$, \textcolor{black}{and}  return to the
first step to find the corresponding $u_0$ and repeat the procedure.

\subsection{Algorithm of reconstruction \textcolor{black}{of} $q$ for a constant $\gamma$}\label{AlgoQ}
\textbf{Step 1.}  We start from  an initial guess $q_0,$ and solve the corresponding Dirichlet problem for \textcolor{black}{the} Helmholtz equation
\begin{eqnarray*}
&& \gamma \triangle u_0 +k^2 q_0(x) u_0=0,\\
&&u_0|_{\del \Omega}=\psi.
\end{eqnarray*}
Solving the direct problem for $\psi=e^{i\arctan{y/x}}$, we obtain $u_0$.

\textbf{Step 2.} We have seen that our inverse problem is asymptotically approached by the  direct problem\begin{equation}\label{OH}
    \left\{%
\begin{array}{ll}
    \gamma\triangle u+\omega^2\frac{j(x)}{|u|^2}u=0, & \hbox{in } \Omega, \\
     u=\psi, & \hbox{on } \del \Omega. \\
\end{array}%
\right.
\end{equation}
We compute the difference
\begin{equation}\label{err}
\eps_0:=\frac{j(x)}{|u_0|^2}-q_0
\end{equation}
and verify
\begin{equation}\label{precision}
|\eps_0|<C_{\textrm{prec}},
\end{equation}
where $C_{\textrm{prec}}$ is our wished order of precision. If  condition~(\ref{precision}) holds, we finish our algorithm and \textcolor{black}{set} $q\equiv q_0.$ Otherwise we \textcolor{black}{go} to the next step.

\textbf{Step 3.} We use now the expression
\[(q_0+\delta q_1)|u_0+\delta u_1|^2=j(x),\]
having \textcolor{black}{the} goal to approximate the known function $j(x)$ with the help
of the small correctors $\delta u_1$ and $\delta q_1$.

\textcolor{black}{By expanding} the expression, we obtain
\begin{eqnarray*}
&&\delta
q_1\left(1+2\delta\left(\frac{\operatorname{Re}u_0\operatorname{Re}u_1+\operatorname{Im}u_0\operatorname{Im}u_1}{|u_0|^2}\right)+\delta^2\frac{|u_1|^2}{|u_0|^2}\right)
=\frac{j(x)}{|u_0|^2}-q_0\\
&&
-2\delta\frac{q_0\left(\operatorname{Re}u_0\operatorname{Re}u_1+\operatorname{Im}u_0\operatorname{Im}u_1\right)}{|u_0|^2}-\delta^2q_0\frac{|u_1|^2}{|u_0|^2}.
\end{eqnarray*}
\textcolor{black}{As in Section~\ref{AlgoG}, we suppose that $\delta \ll 1$ and  that $\delta\max\limits_x|q_1|$ and $\delta\max\limits_x |u_1|$ are of the order of $\delta$. Consequently,} we consider only terms of order not smaller than $\delta$:
\[\delta q_1=\frac{j(x)}{|u_0|^2}-q_0
-2\delta\frac{q_0\left(\operatorname{Re}u_0\operatorname{Re}u_1+\operatorname{Im}u_0\operatorname{Im}u_1\right)}{|u_0|^2}.\]
\textcolor{black}{To find} the corrector $\tilde{u}_1=\delta u_1$, we \textcolor{black}{expand} the
following equation
\[ \gamma \triangle(u_0+\delta u_1)+k^2(q_0+\delta
q_1)(u_0+\delta u_1)=0.\]

Considering the terms of order no smaller than $\delta$ and
replacing   $\delta q_1$ by the approximated formula, we find
$\tilde{u}_1$ as a solution of the following problem
\begin{eqnarray}
&&\gamma \triangle \tilde{u}_1+k^2\frac{j(x)}{|u_0|^2}
\tilde{u}_1-2k^2 u_0
\frac{q_0\left(\operatorname{Re}u_0\operatorname{Re}\tilde{u}_1+\operatorname{Im}u_0\operatorname{Im}\tilde{u}_1\right)}{|u_0|^2}=-\eps_0 k^2u_0,\label{QU1cor}\\
&& \tilde{u}_1|_{\del \Omega}=0.\nonumber
\end{eqnarray}
We solve the problem and obtain \textcolor{black}{$\tilde{u}_1$}.

\textbf{Step 4.} We calculate
\begin{equation}
 \tilde{q}=q_0+\delta q_1=\frac{1}{|u_0|^2}\left
(j(x)-2q_0(\operatorname{Re}u_0\operatorname{Re}\tilde{u}_1+\operatorname{Im}u_0\operatorname{Im}\tilde{u}_1)\right).\label{updateQ}
\end{equation}

We \textcolor{black}{set} now $q_0\equiv \tilde{q}$, \textcolor{black}{and} return to the
first step to find the corresponding $u_0$ and repeat the procedure.

\section{Numerical results}\label{num}
\begin{figure}[h!]
 \begin{center}
\subfigure[]{\psfrag{X}[c]{{\tiny Number of iterations}}
	\psfrag{Y}[c]{{\tiny Error}}
\psfrag{Gainf}{{\tiny $\|\gamma-\tilde{\gamma}\|_{L_\infty}$}}
\psfrag{Ga1}{{\tiny $\|\gamma-\tilde{\gamma}\|_{L_1}$}}
\psfrag{Ga2}{{\tiny $\|\gamma-\tilde{\gamma}\|_{L_2}$}}
\psfrag{Gerr2}{{\tiny $\|J/|\nabla u|^2-\gamma\|_{L_2}$}}
\psfrag{Gerrinf}{{\tiny $\|J/|\nabla u|^2-\gamma\|_{L_\infty}$}}
\psfrag{Gj0}{{\tiny $\|J/|\nabla u|^2-\gamma\|^2_{L_2}$}}
\psfrag{Gj1}{{\tiny $\|J-\gamma|\nabla u|^2\|^2_{L_2}$}}
  \includegraphics[width=0.45\linewidth]{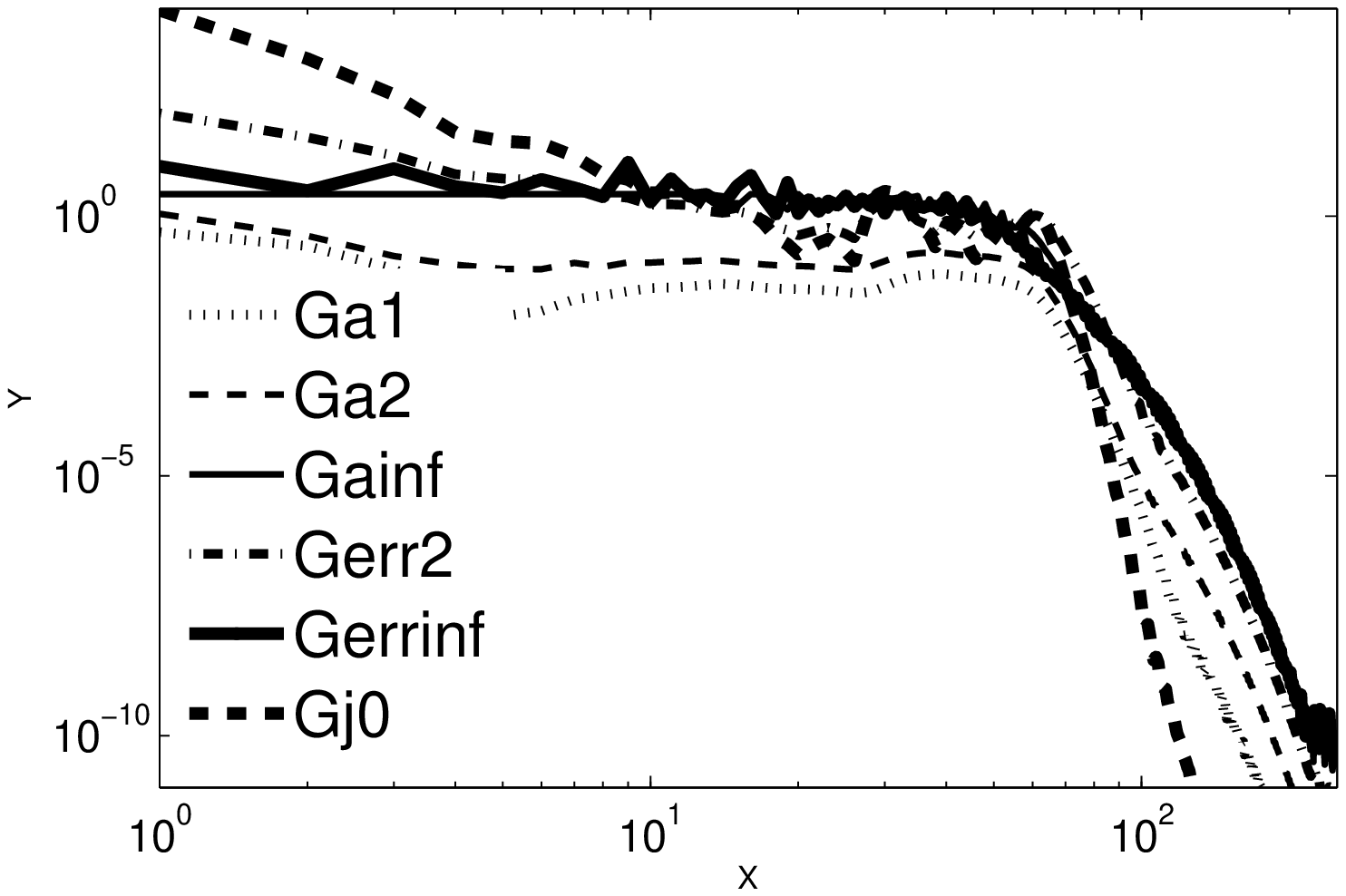}}
\subfigure[]{\psfrag{X}[c]{\tiny{Number of iterations}}
	\psfrag{Y}[c]{{\tiny Error}}
\psfrag{Qainf}{{\tiny $\|q-\tilde{q}\|_{L_\infty}$}}
\psfrag{Qa1}{{\tiny $\|q-\tilde{q}\|_{L_1}$}}
\psfrag{Qa2}{{\tiny $\|q-\tilde{q}\|_{L_2}$}}
\psfrag{Qerr2}{{\tiny $\|j/| u|^2-q\|_{L_2}$}}
\psfrag{Qerrinf}{{\tiny $\|j/|u|^2-q\|_{L_\infty}$}}
\psfrag{Qj0}{{\tiny $\|j/| u|^2-q\|^2_{L_2}$}}
\psfrag{Qj1}{{\tiny $\|j-q| u|^2\|^2_{L_2}$}}
\includegraphics[width=0.45\linewidth]{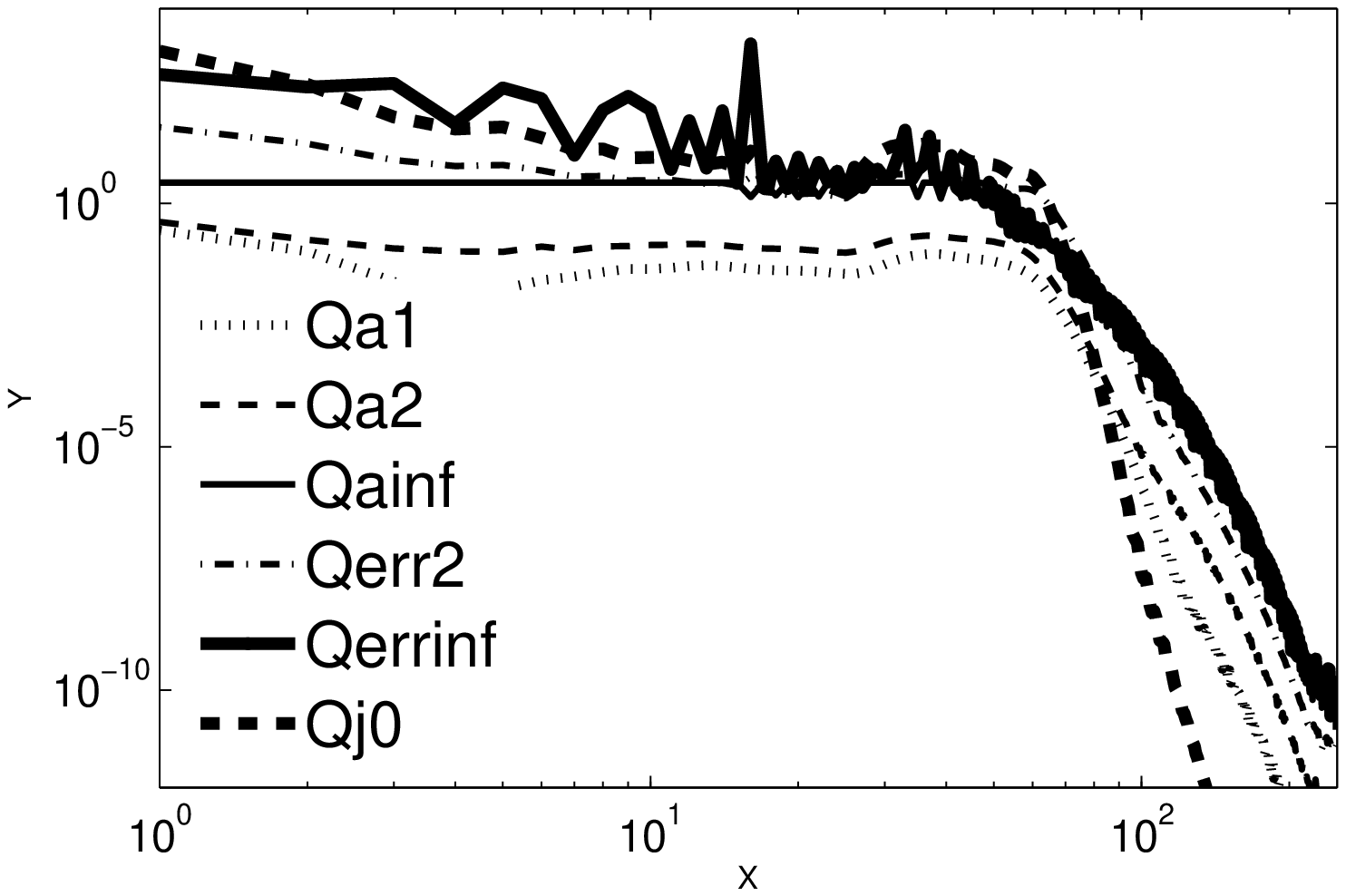}}
 \end{center}
\caption{Convergence results for a ``perfect'' mesh on (a) $\gamma$  and (b) $q$.}\label{ConvResPI}
\end{figure}

To study the efficiency of this approach, we have tested this method on various problems and domains, using the partial differential equation solver FreeFem++~\cite{Free}. We present here one such test. The domain $\Omega$ is a disk of radius $8$ centred at the origin, which contains three inclusions: a triangle, an $L$-shaped domain and an ellipse, which represents a convex object, a non-convex object, and an object with a smooth boundary respectively.
\begin{figure}[h!]
 \begin{center}
\subfigure[]{
\resizebox*{3cm}{!}{
 \includegraphics[width=3cm]{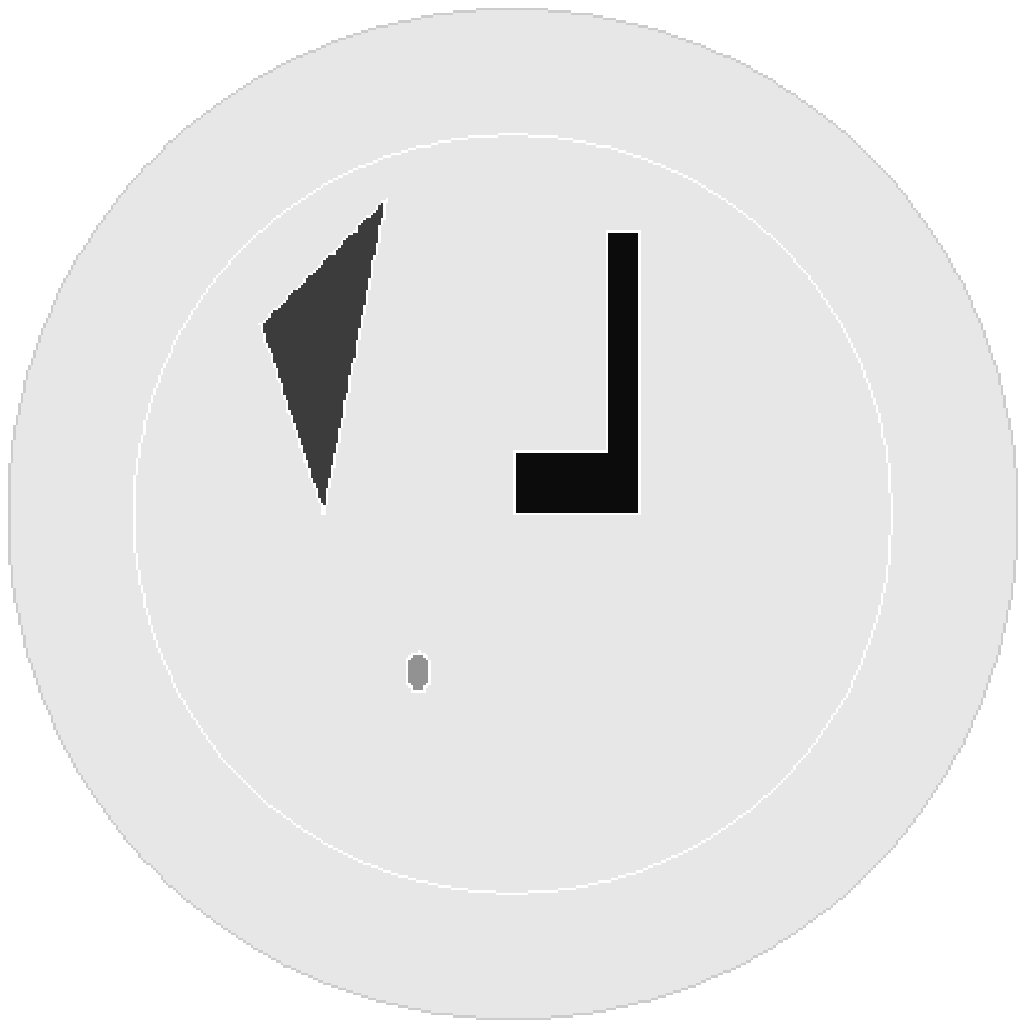}}}
\subfigure[]{
\resizebox*{3cm}{!}{
\includegraphics[width=3cm]{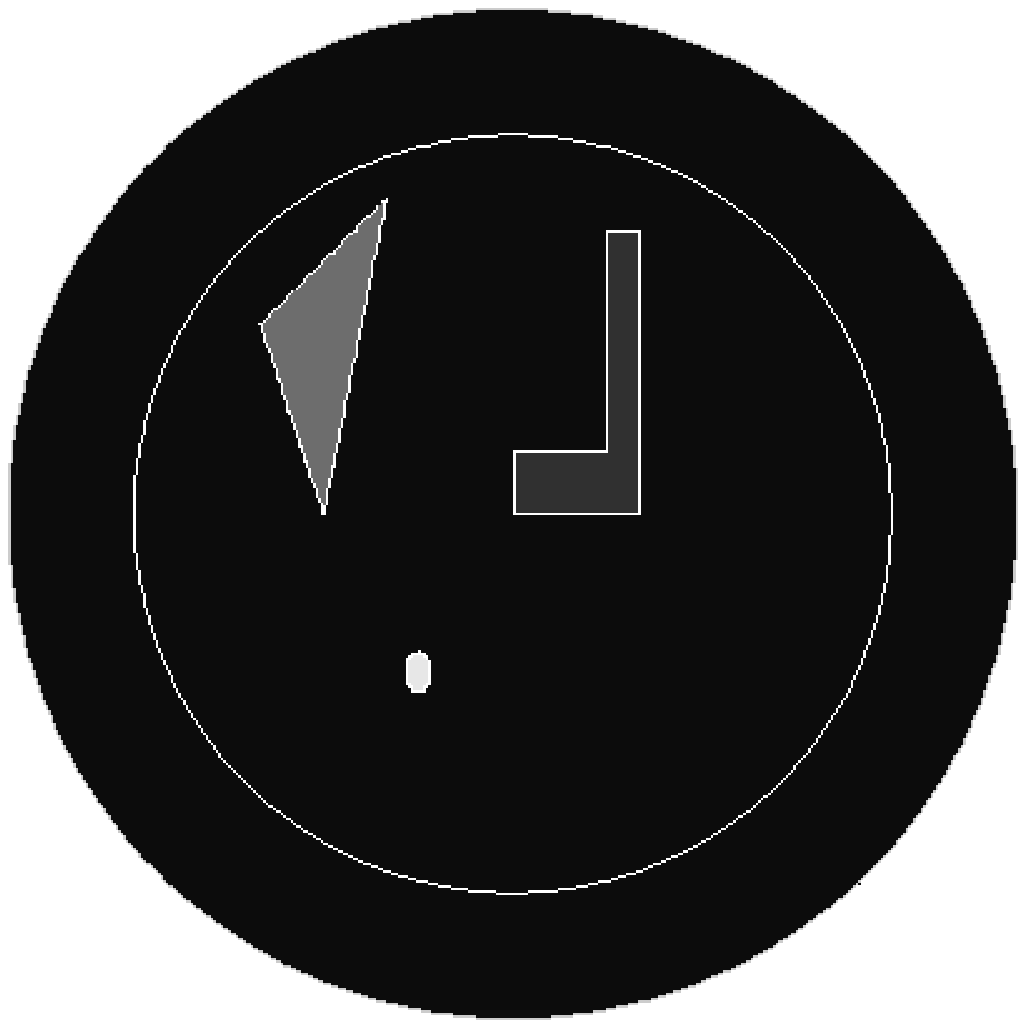}}}
\subfigure[]{
\resizebox*{3cm}{!}{
\includegraphics[width=3cm]{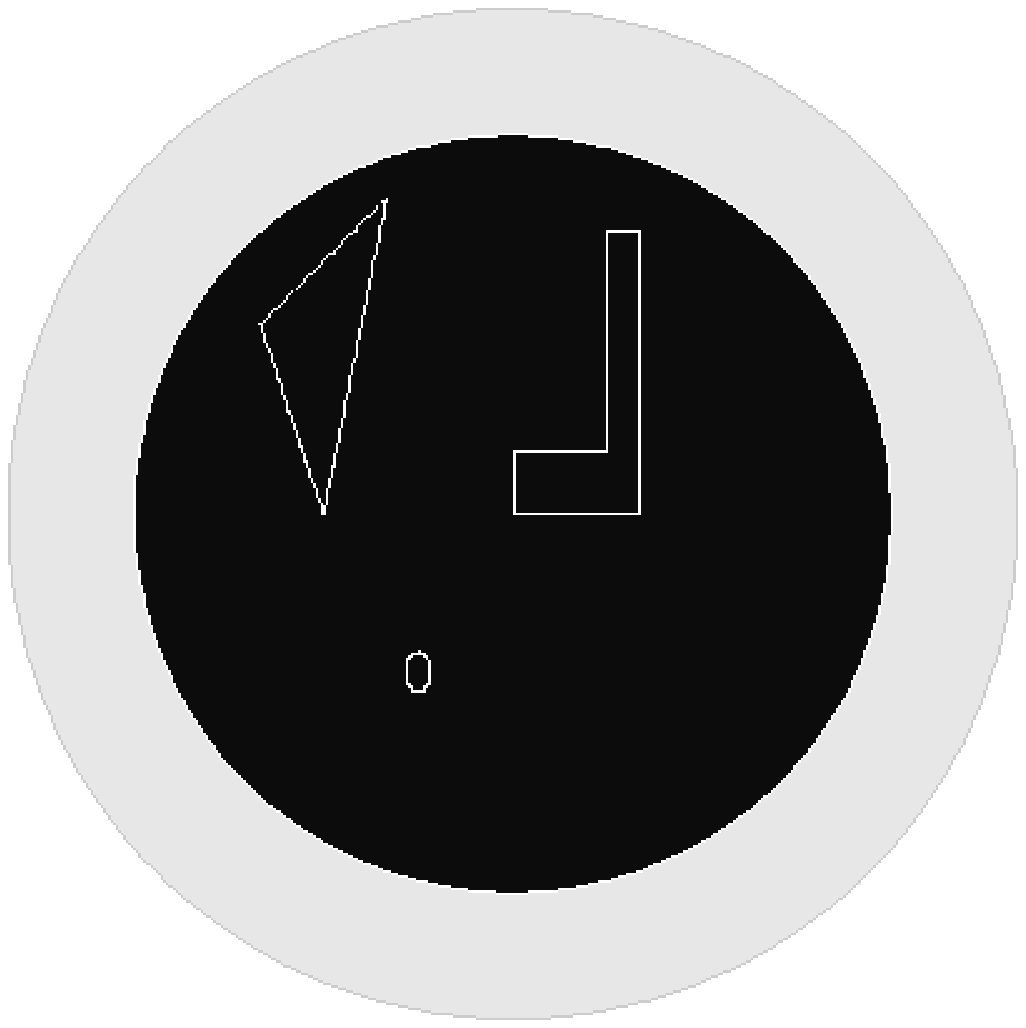}}}
\subfigure[]{
\resizebox*{3cm}{!}{
\includegraphics[width=3cm]{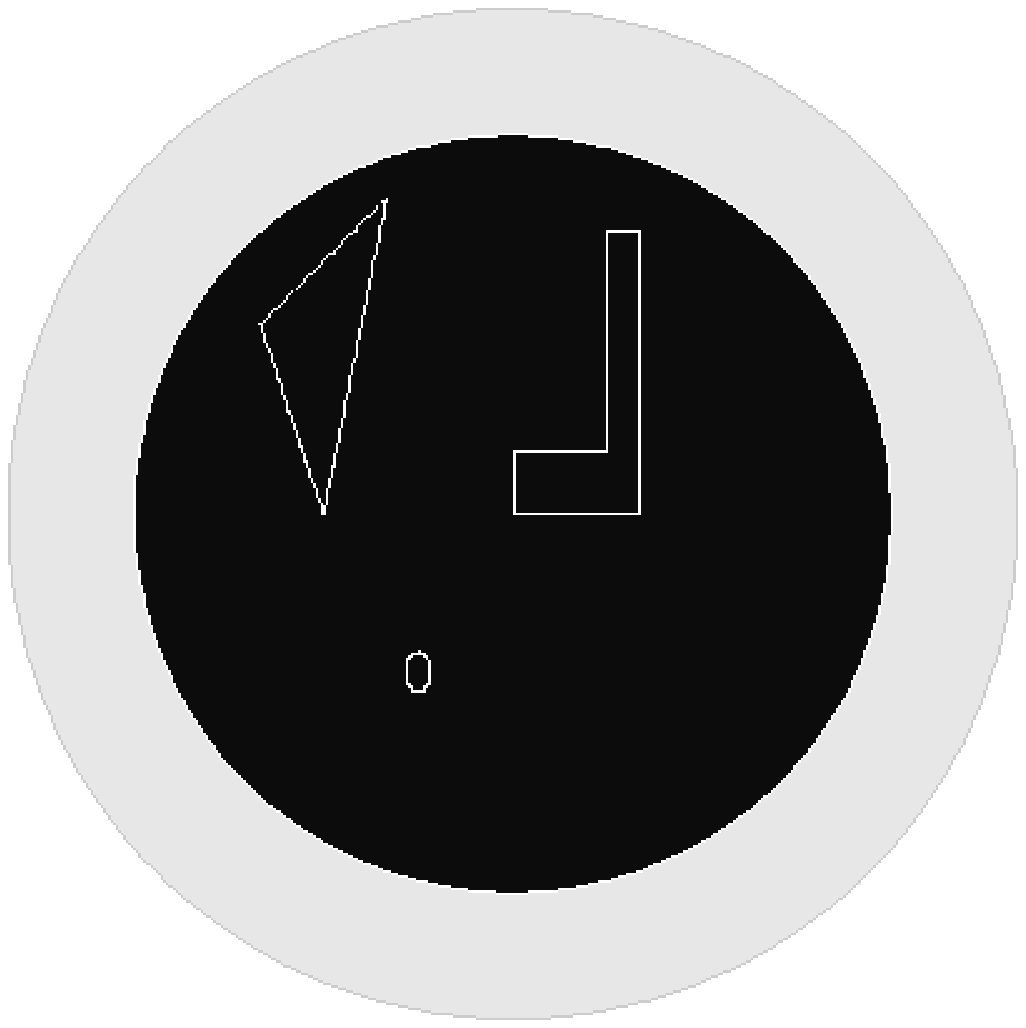}}}
\end{center}
\caption{(a) Distribution of the conductivity $\gamma$. (b) Distribution of the permittivity $q$. (c) Initial guess for $\gamma$. (d) Initial guess for $q$.}\label{PRoj}
\end{figure}
In  Figure~\ref{PRoj} (a) (respectively (b)) the background conductivity (respectively permittivity) is equal to $1$ (\textcolor{black}{respectively} $3$),  the conductivity (respectively permittivity) takes the value $2.5$ (respectively $2$) in the triangle, $1.75$ (respectively $1$) in the ellipse and $3.05$ (respectively $2.55$) in the L-shaped domain for the two frequencies $k_1=\pi\cdot 10^{3}$ and $k_2=\pi\cdot 10^{-3}$.  We purposely choose values corresponding to small and large contrast with the background.
The initial guess in Figure~\ref{PRoj} (c) (respectively (d))  is equal to $3.5$ (respectively $11.5$) inside the disk of  \textcolor{black}{radius $6$ centred} at the origin, and equal to the supposedly known conductivity (permittivity) $1$ (respectively $3$) near the boundary (outside the disk of radius $6$).
\begin{figure}
 \begin{center}
\subfigure[]{
\resizebox*{3cm}{!}{
\includegraphics[width=3cm]{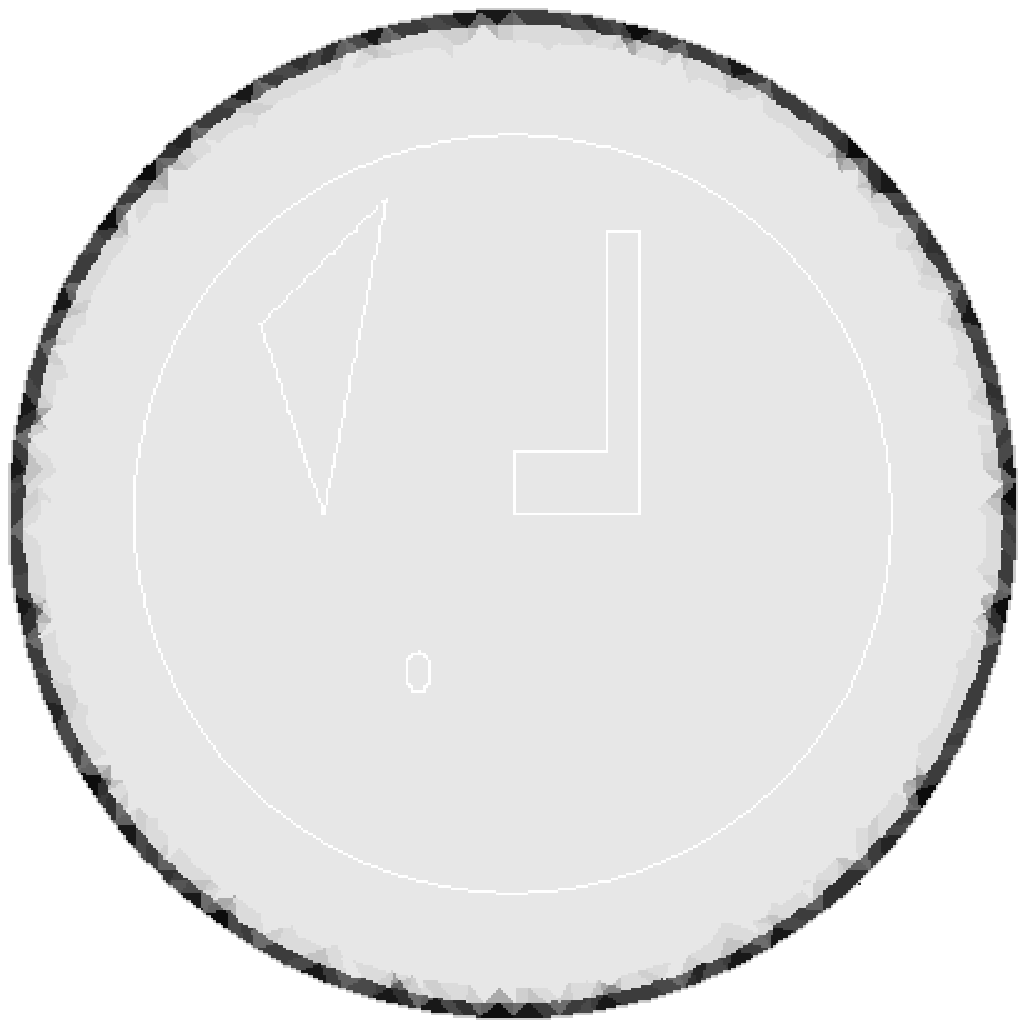}}}
\subfigure[]{
\resizebox*{3cm}{!}{
\includegraphics[width=3cm]{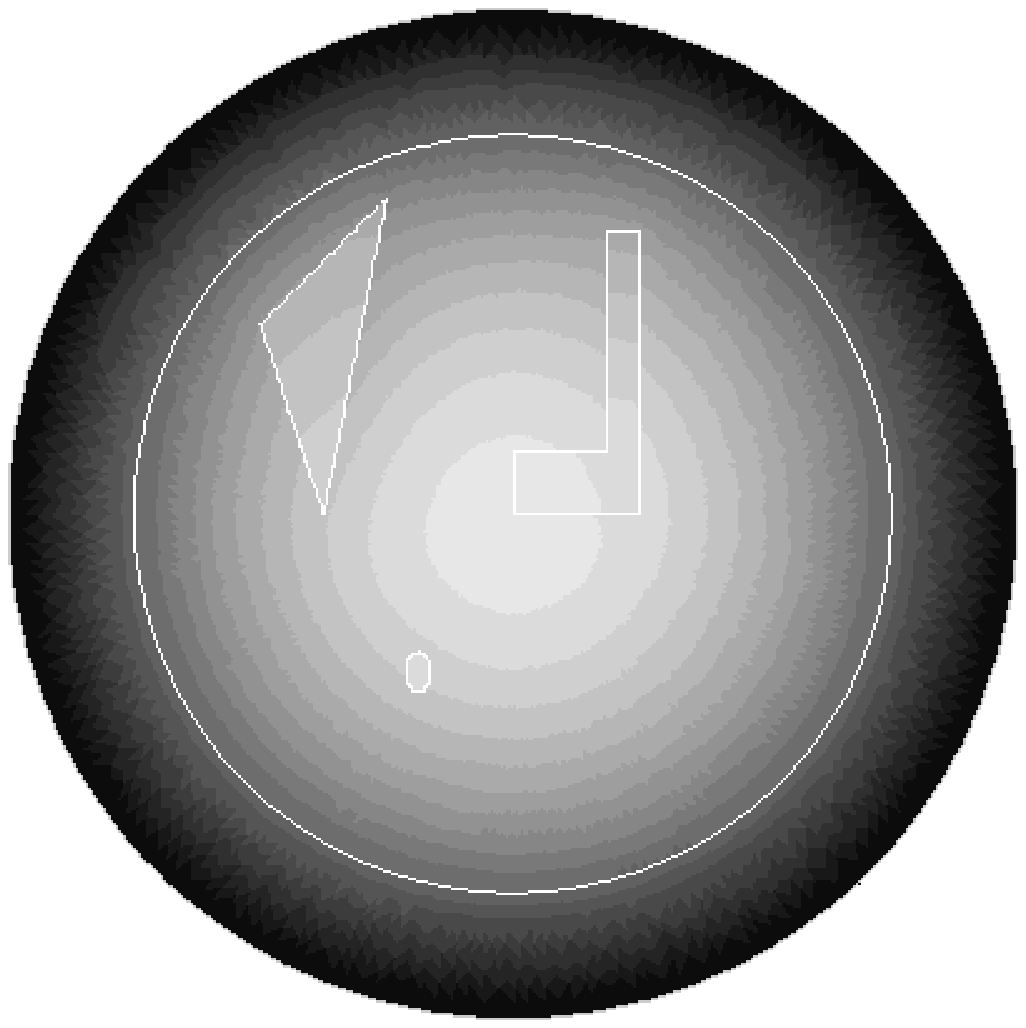}}}
\subfigure[]{
\resizebox*{3cm}{!}{
\includegraphics[width=3cm]{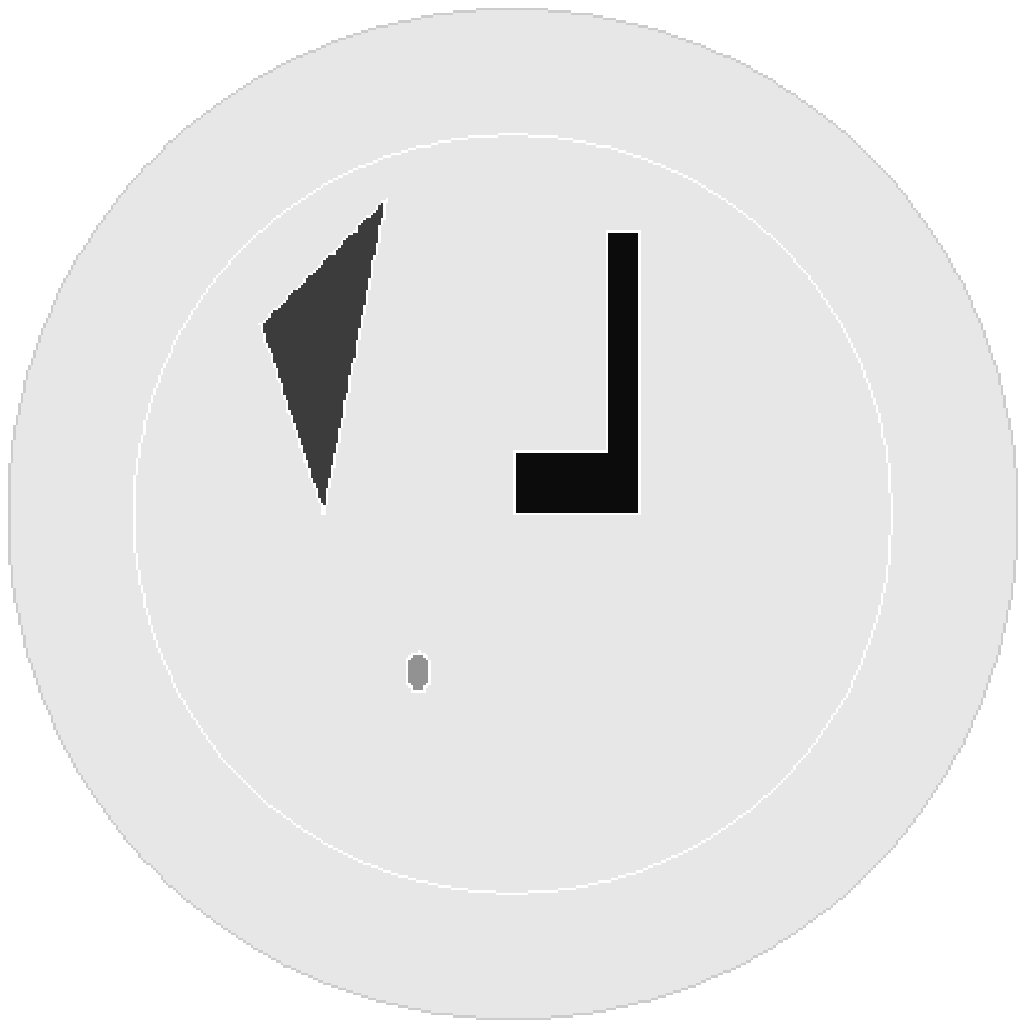}}}
\subfigure[]{
\resizebox*{3cm}{!}{\includegraphics[width=3cm]{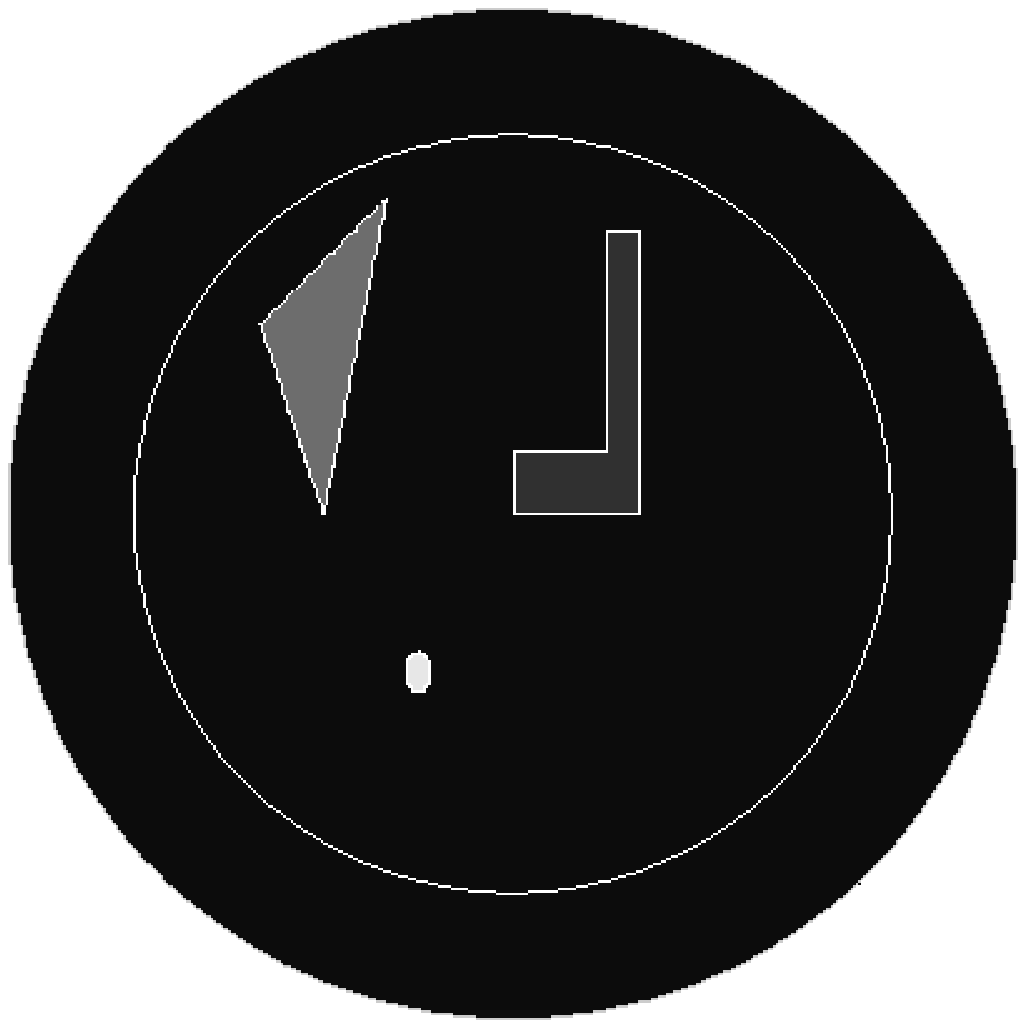}}}
 \end{center}
\caption{Reconstruction test with a ``perfect'' mesh. (a) \textcolor{black}{Collected} data $J$ for the reconstruction of $\gamma$. (b)  \textcolor{black}{Collected} data $j$ for the reconstruction of $q$. (c)  \textcolor{black}{Reconstructed} conductivity $\gamma$. (d) \textcolor{black}{Reconstructed} permittivity $q$.}\label{ExacteRes}
\end{figure}
Figure~\ref{ExacteRes} shows the result of the reconstruction when perfect measures (with ``infinite'' precision) are available. For all presented numerical results we use as boundary potential $\psi=e^{i \arctan (x/y)}$. Figure~\ref{ExacteRes} (a) and (b) represents the collected data, $J(x)$ and $j(x)$. For known values of the contrast, we remark that \textcolor{black}{through} we can 'see' the structure of the permittivity. On Figures~\ref{ExacteRes} (c) and (d), the reconstructed conductivity and permittivity are represented: they perfectly match the target.
 Figure~\ref{ConvResPI} (a) (respectively~\ref{ConvResPI} (b)) presents different errors as functions of the iteration for $\gamma$ (respectively $q$). 
\begin{figure}[h!]
 \begin{center}
  \subfigure[]{
\resizebox*{3cm}{!}{\includegraphics[width=3cm]{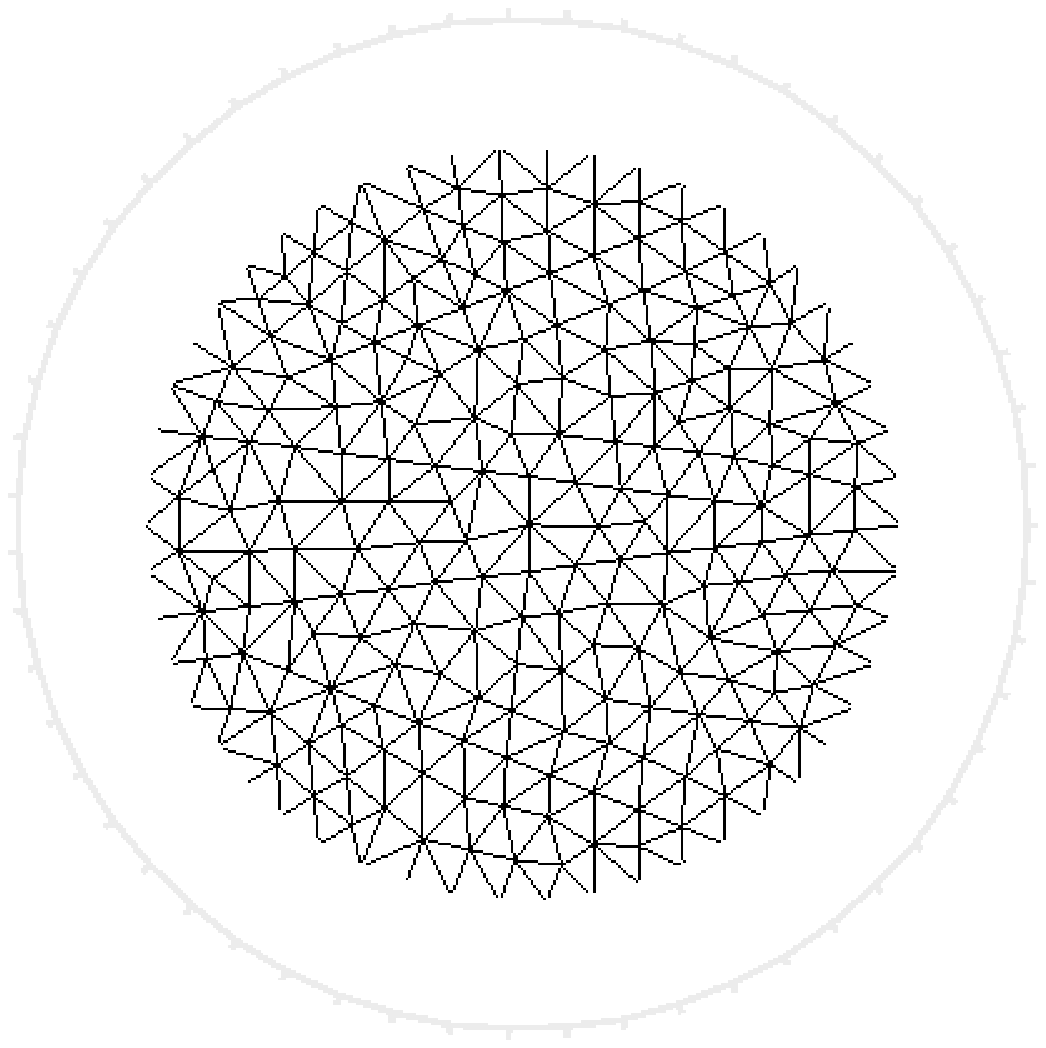}}}
\subfigure[]{
\resizebox*{3cm}{!}{
\includegraphics[width=3cm]{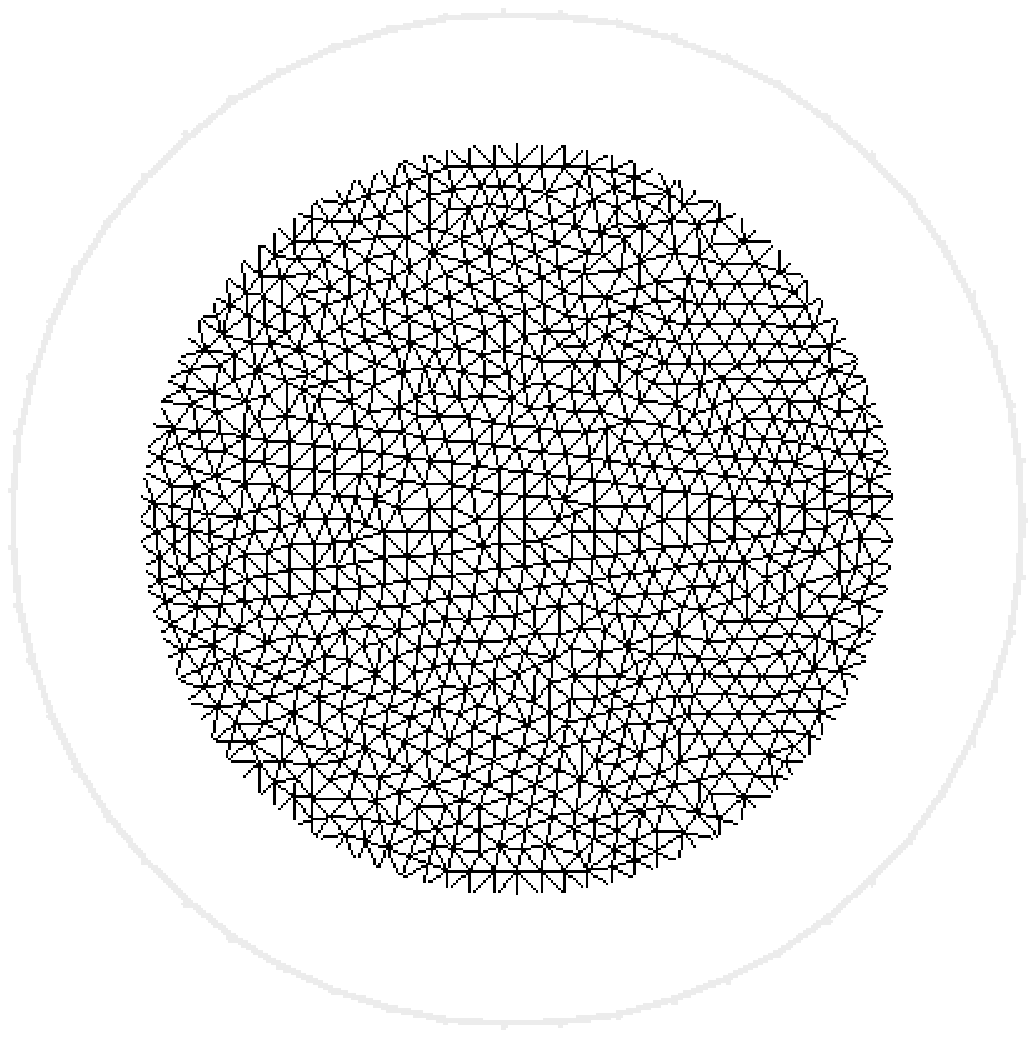}}}
\subfigure[]{
\resizebox*{3cm}{!}{
\includegraphics[width=3cm]{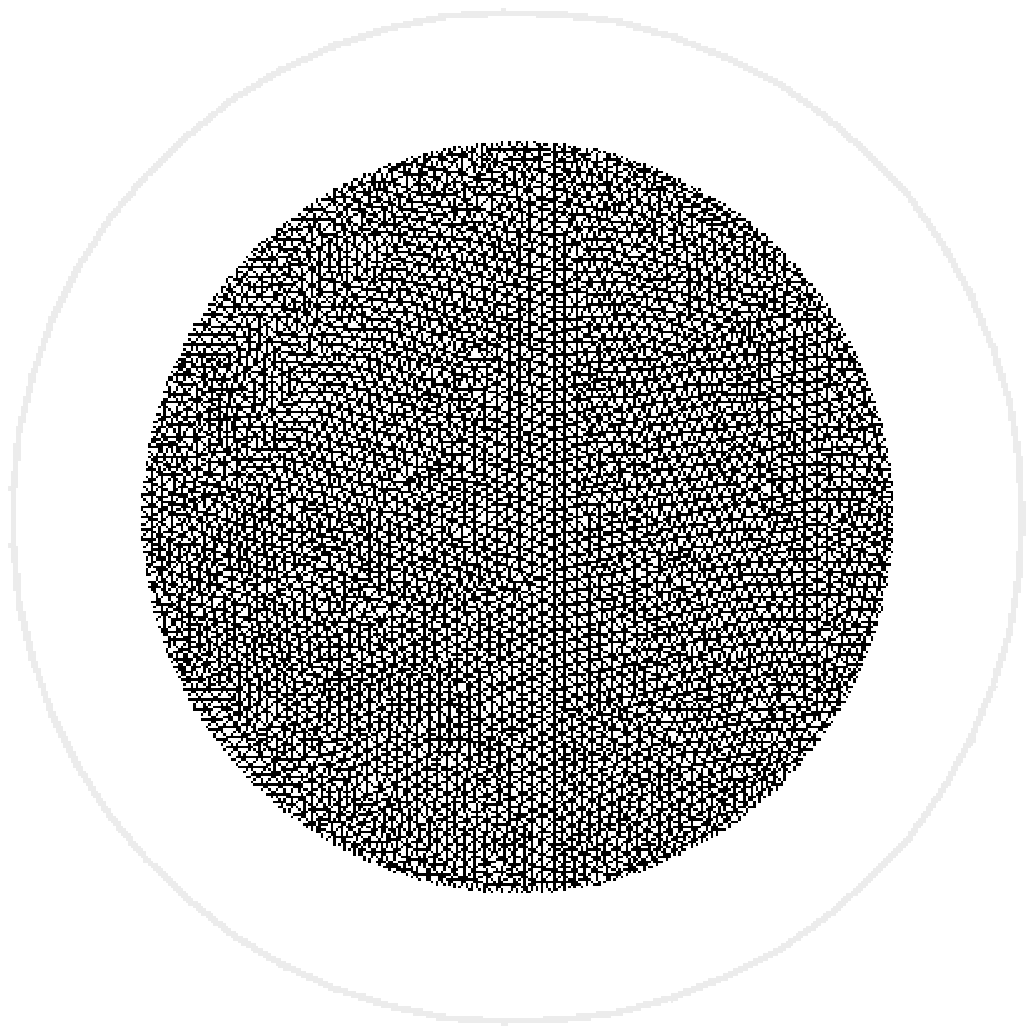}}}
\subfigure[]{
\resizebox*{3cm}{!}{\includegraphics[width=3cm]{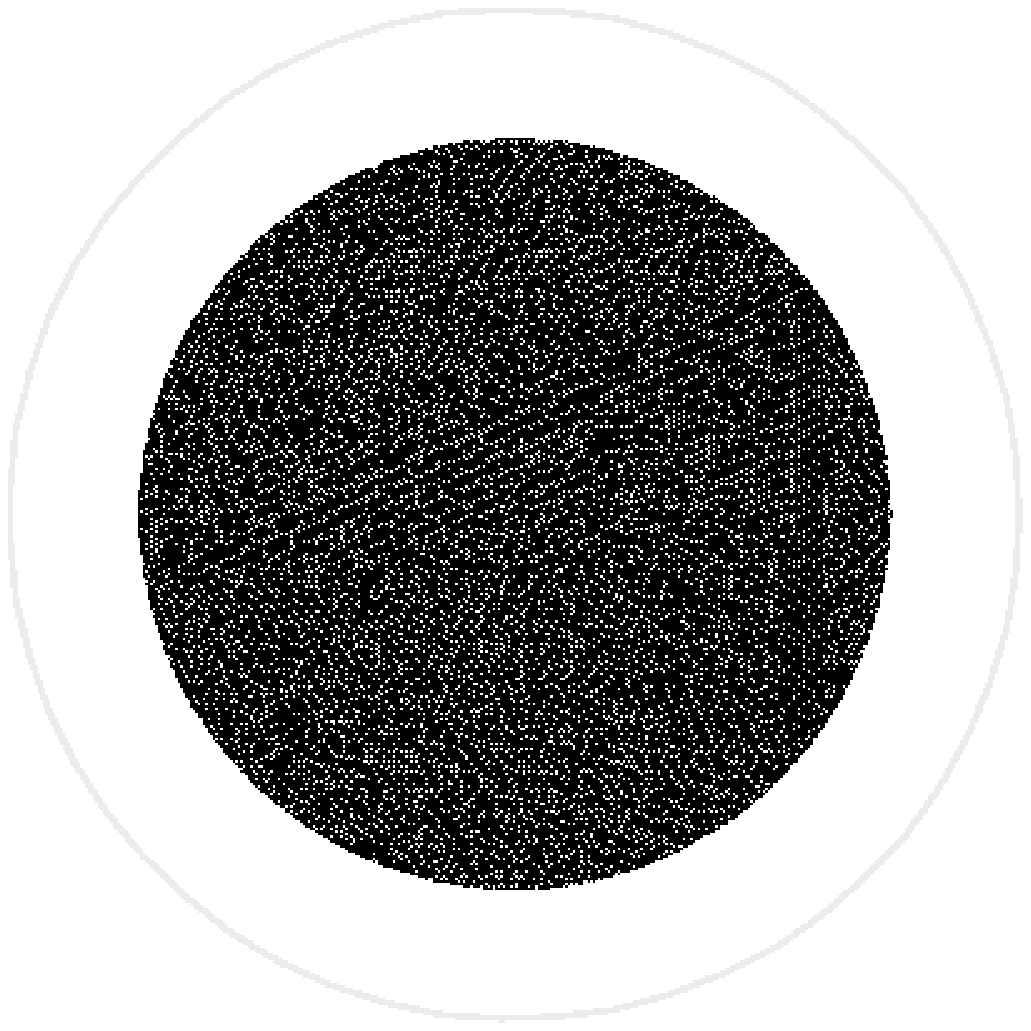}}}
 \end{center}
\caption{``Imperfect'' meshes for: (a) $50$ boundary points; (b) $100$ boundary points; (c) $200$ boundary points; and (d) $400$ boundary points.}\label{Mesh}
\end{figure}

\begin{figure}[h!]
 \begin{center}
\subfigure[]{
\resizebox*{3cm}{!}{\includegraphics[width=3cm]{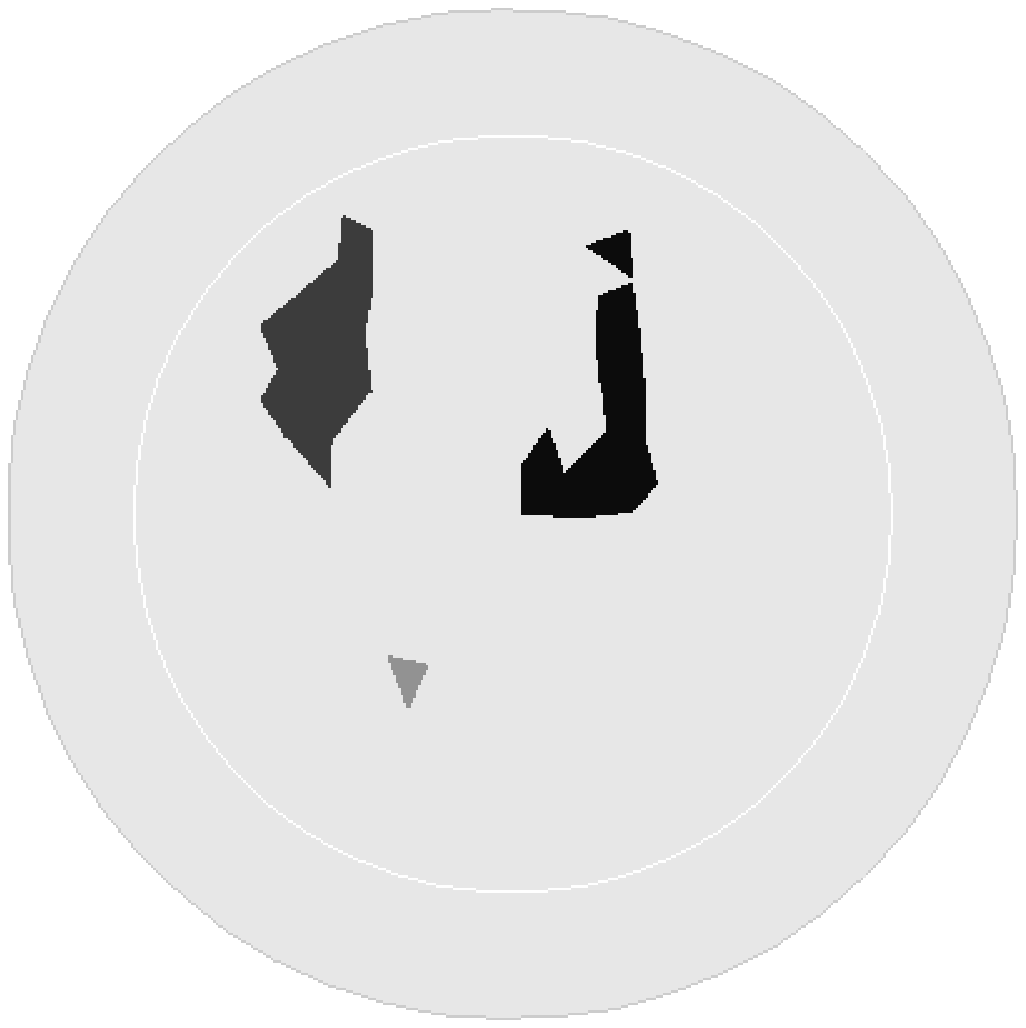}}}
\subfigure[]{
\resizebox*{3cm}{!}{
\includegraphics[width=3cm]{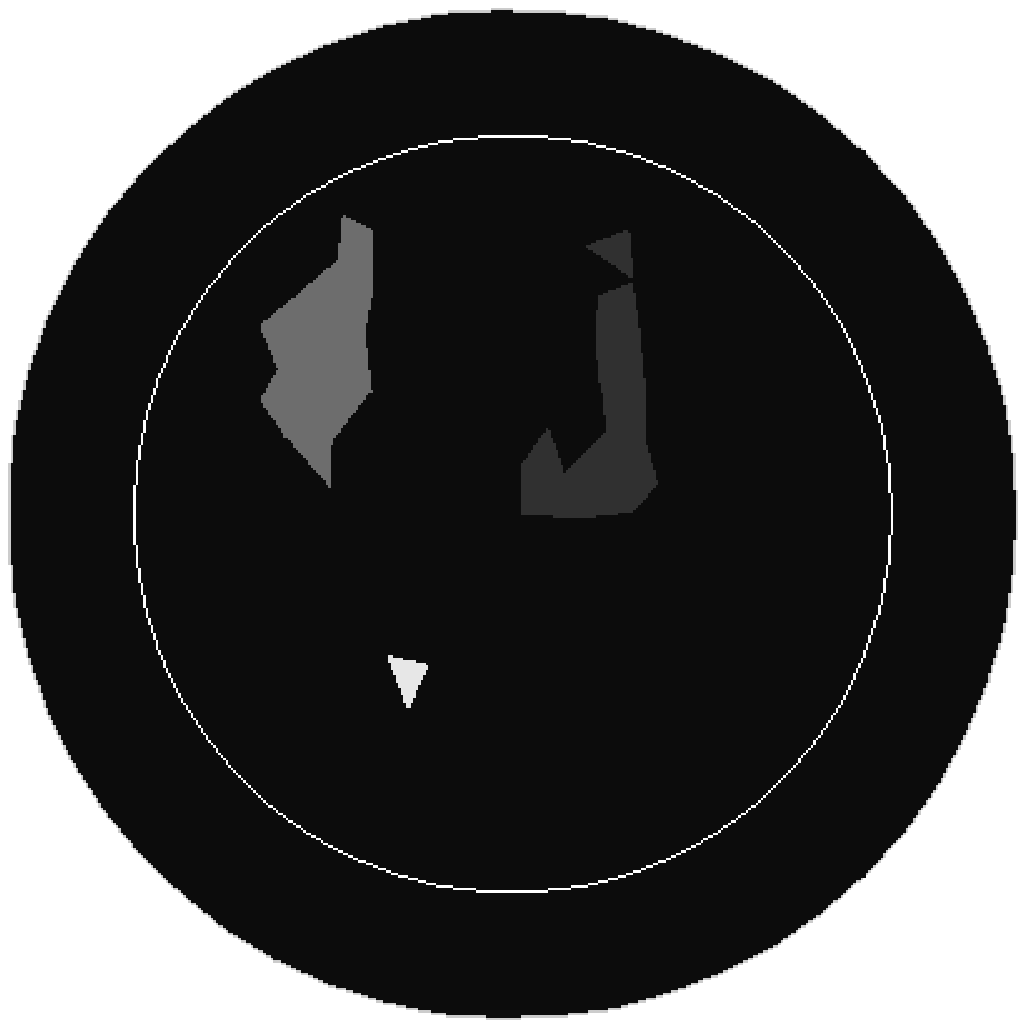}}}
\subfigure[]{
\resizebox*{3cm}{!}{
\includegraphics[width=3cm]{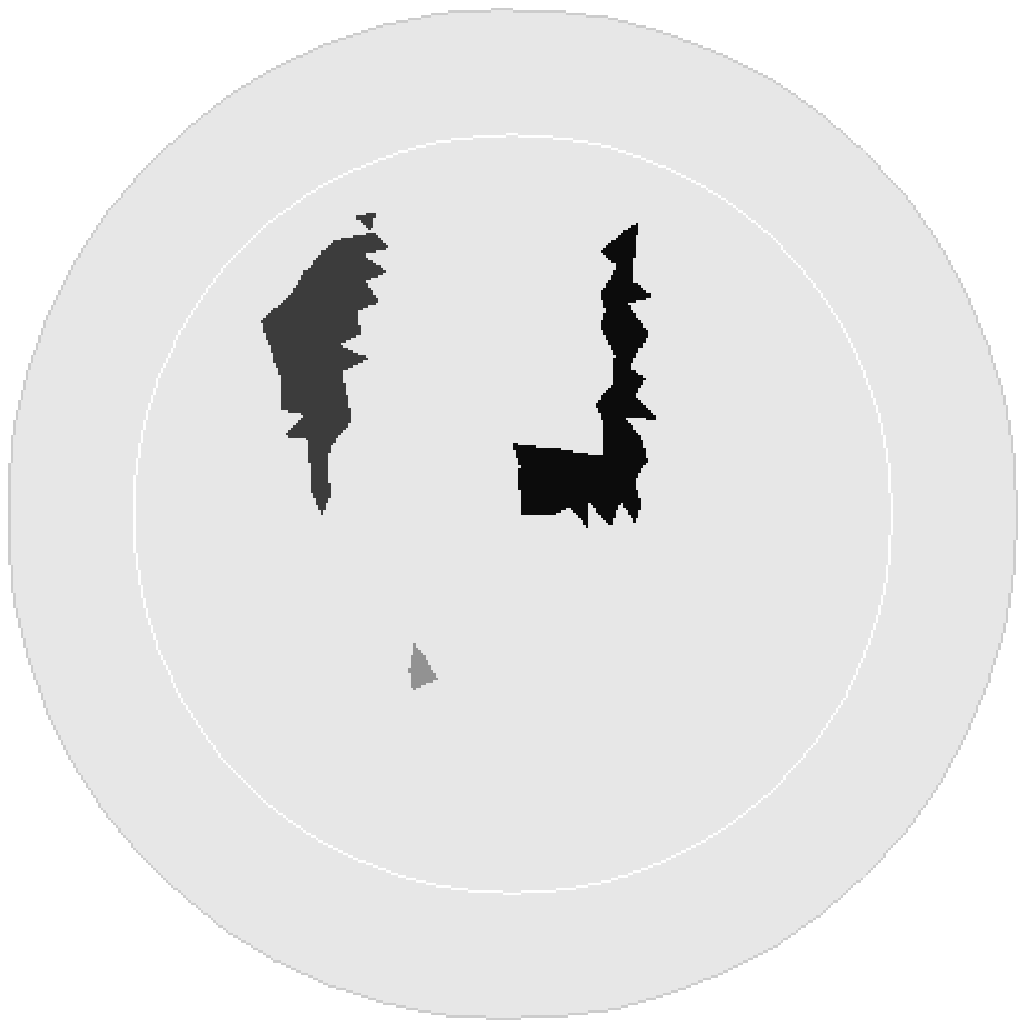}}}
\subfigure[]{
\resizebox*{3cm}{!}{\includegraphics[width=3cm]{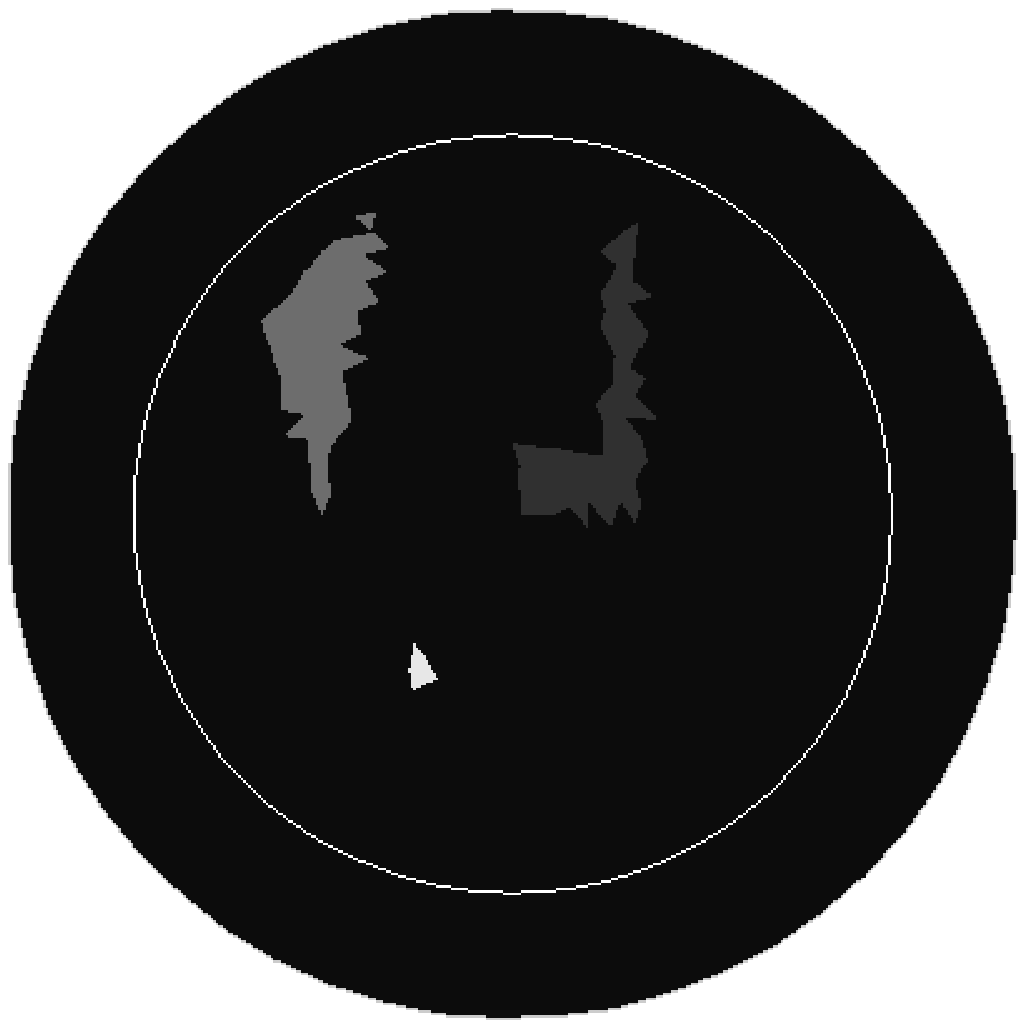}}}
\subfigure[]{
\resizebox*{3cm}{!}{\includegraphics[width=3cm]{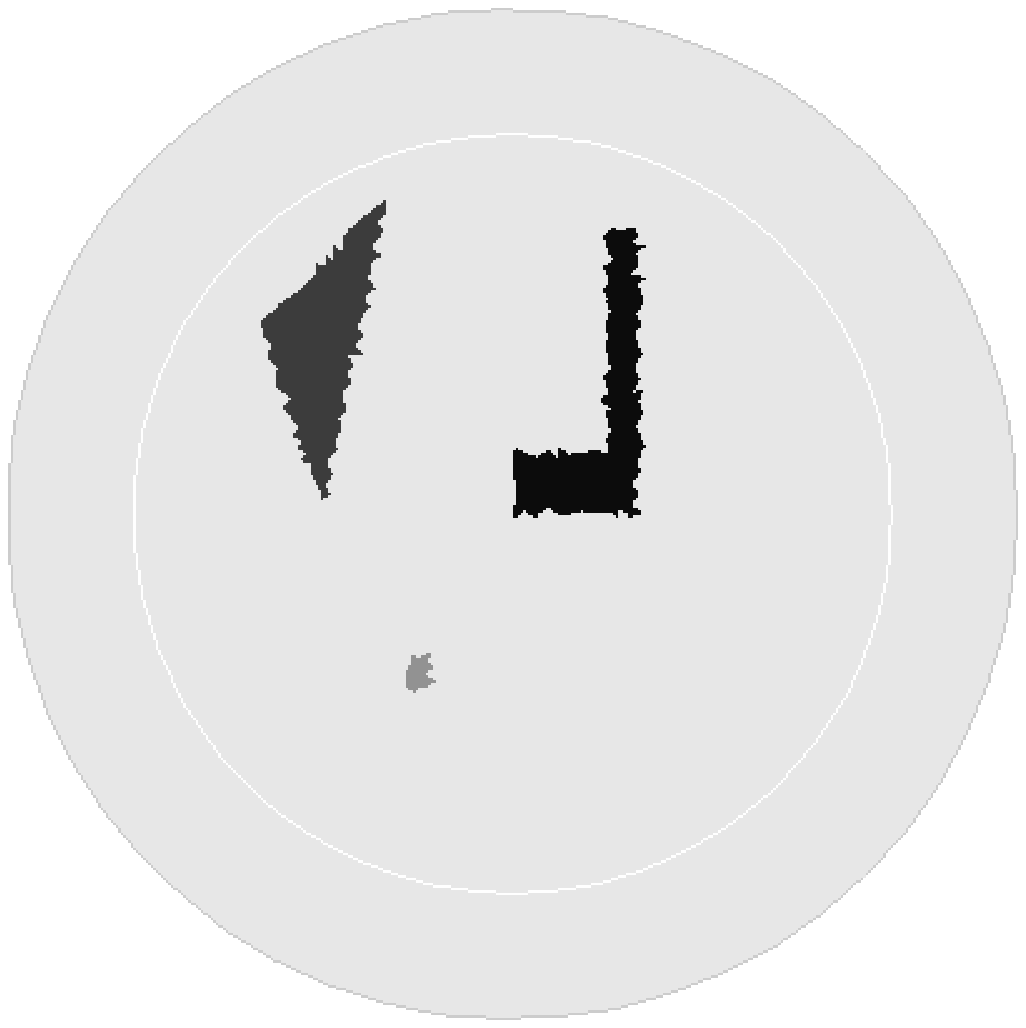}}}
\subfigure[]{
\resizebox*{3cm}{!}{
\includegraphics[width=3cm]{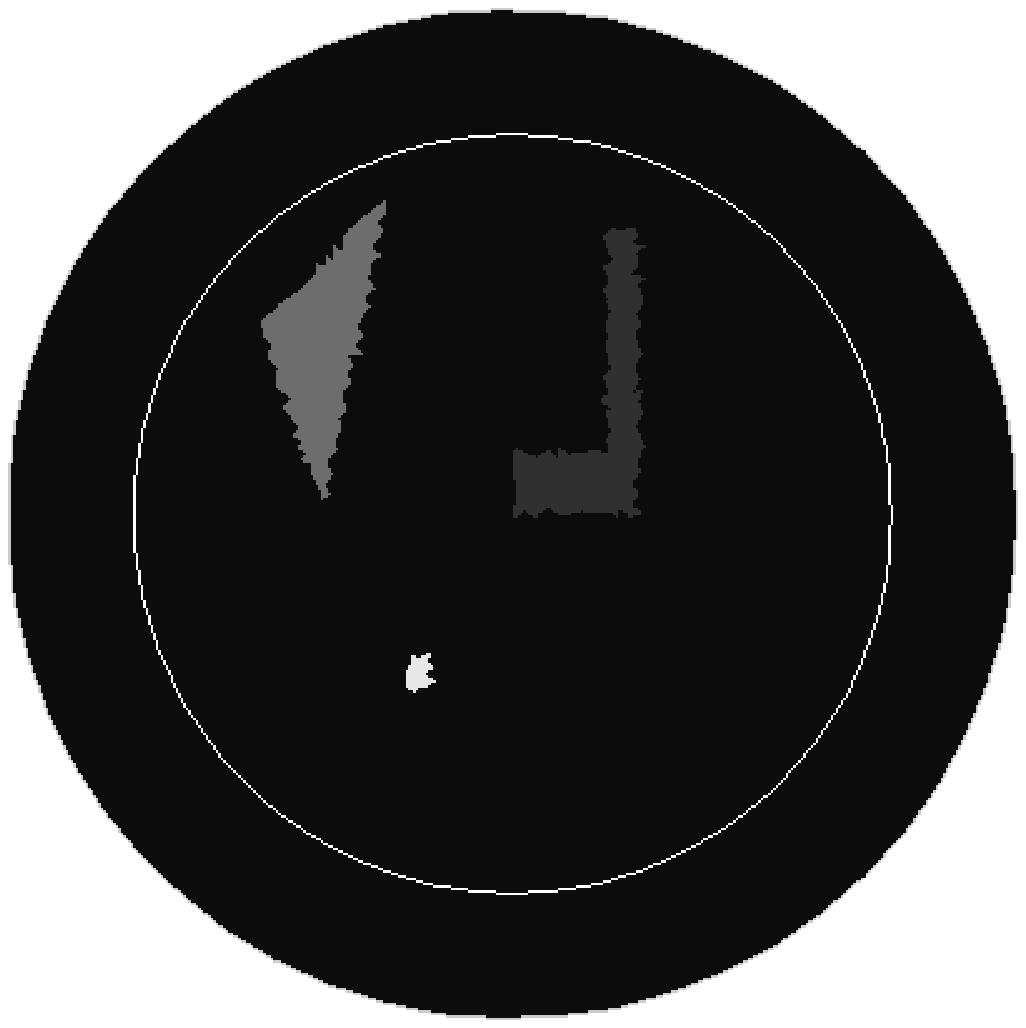}}}
\subfigure[]{
\resizebox*{3cm}{!}{
\includegraphics[width=3cm]{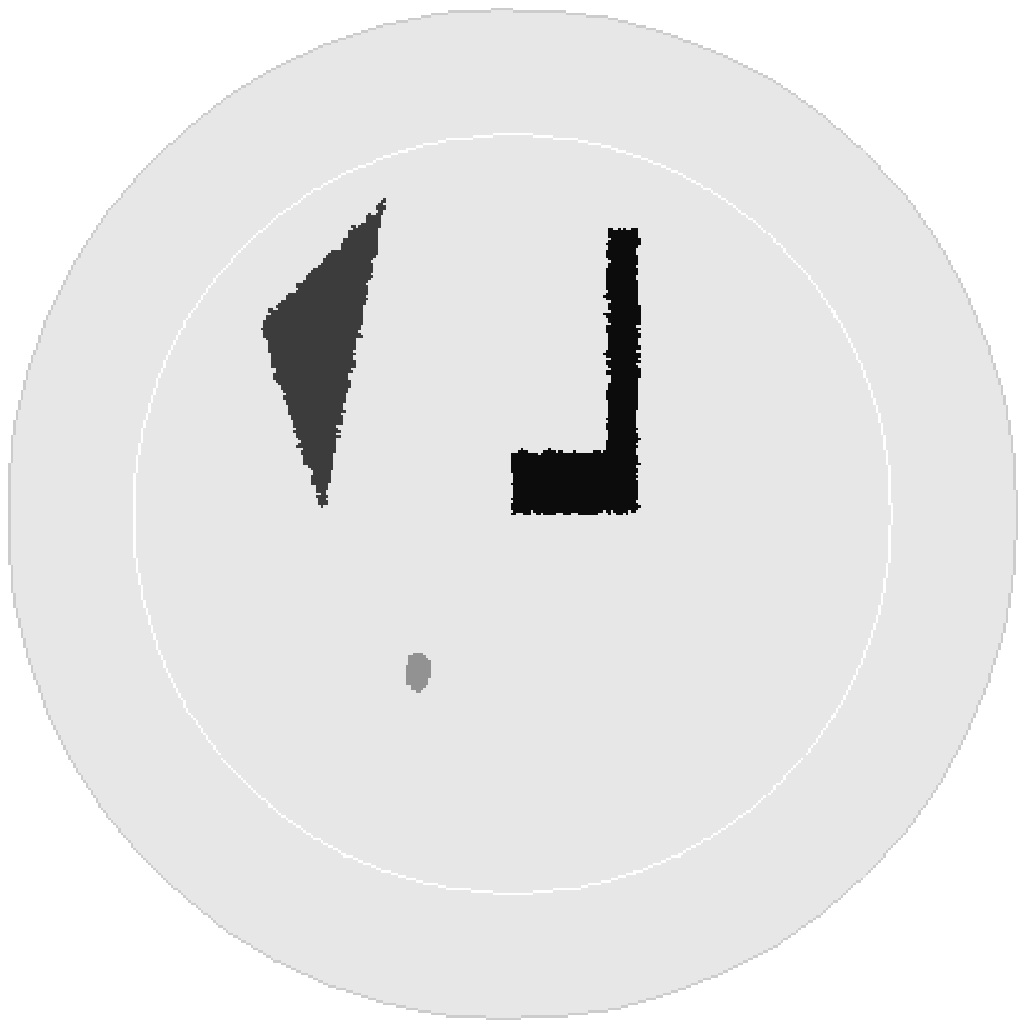}}}
\subfigure[]{
\resizebox*{3cm}{!}{\includegraphics[width=3cm]{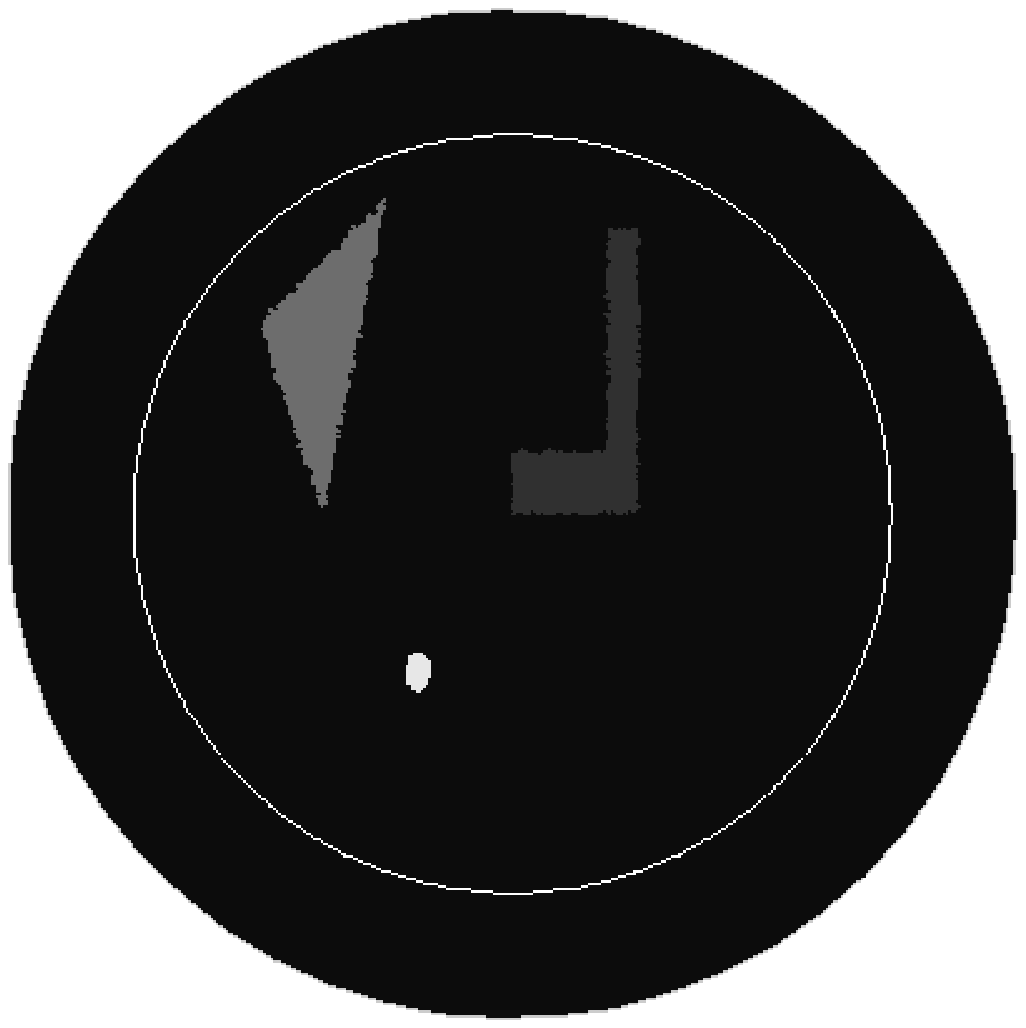}}}
\end{center}
\caption{Reconstruction tests for different ``imperfect'' meshes\textcolor{black}{:} (a) $\gamma$ and (b) $q$ using a regular mesh with $50$ boundary points\textcolor{black}{,} (c) $\gamma$ and (d) $q$ using a regular mesh with $100$ boundary points\textcolor{black}{,} (e) $\gamma$ and (f) $q$ using a regular mesh with $200$ boundary points \textcolor{black}{and} (g) $\gamma$ and (h) $q$ using a regular mesh with $400$ boundary points.}\label{DiferEl}
\end{figure}

\begin{figure}[h!]
 \begin{center}
 \psfrag{X}[c]{{\small Number of iterations}}
	\psfrag{Y}[c]{{\small Error}}
\psfrag{errinfQ400}{{\tiny Error $L_\infty$ on $q$ for 400 electrodes}}
\psfrag{errinfG400}{{\tiny Error $L_\infty$ on $\gamma$ for 400 electrodes}}
\psfrag{errinfQ200}{{\tiny Error $L_\infty$ on $q$ for 200 electrodes}}
\psfrag{errinfG200}{{\tiny Error $L_\infty$ on $\gamma$ for 200 electrodes}}
\psfrag{errinfQ100}{{\tiny Error $L_\infty$ on $q$ for 100 electrodes}}
\psfrag{errinfG100}{{\tiny Error $L_\infty$ on $\gamma$ for 100 electrodes}}
\psfrag{erinfQ50}{{\tiny Error $L_\infty$  on $q$ for 50 electrodes}}
\psfrag{errinfG50}{{\tiny Error $L_\infty$  on $\gamma$ for 50 electrodes}}
\includegraphics[width=\linewidth]{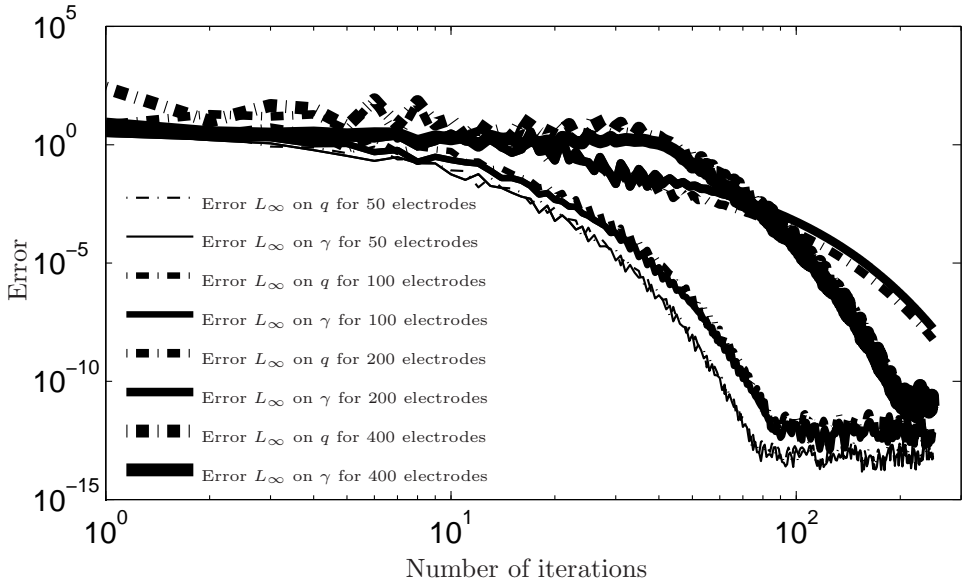}
\end{center}
\caption{Convergence results. Errors $\|j/|u|^2-q\|_{L_\infty}$ and $\|J/|\nabla u|^2-\gamma\|_{L_\infty}$ for meshes with different number of boundary points on $q$ and $\gamma$.}\label{Conv}
\end{figure}

We have also considered imperfect data. We  ran the  reconstruction algorithm with the same conditions, but now assume that the data was measured at the nodes of a regular mesh on the disk, with $50$, $100$, $200$ and $400$ boundary points (see meshes on Figure~\ref{Mesh}). Figure~\ref{DiferEl} shows the obtained reconstructions, which still perfectly match the target. The convergence result for different number of boundary points is given on Figure~\ref{Conv} for  the errors $\|j/|u|^2-q\|_{L_\infty}$ and $\|J/|\nabla u|^2-\gamma\|_{L_\infty}$.  We can observe that the convergence is exponential and that it is even more faster for  meshes of $50$ and $100$ boundary points than for meshes of $200$ or $400$ boundary points. 
%
%
\begin{figure}[h!]
\begin{center}
\subfigure[]{
\psfrag{d1}{{\tiny $k_1=10$, $k_2=10^{-1}$}}
\psfrag{d2}{{\tiny $k_1=10^2$, $k_2=10^{-2}$}}
\psfrag{d3}{{\tiny $k_1=10^3$, $k_2=10^{-3}$}}
\psfrag{d4}{{\tiny $k_1=10^5$, $k_2=10^{-5}$}}
\psfrag{d5}{{\tiny $k_1=10^7$, $k_2=10^{-7}$}}
	\psfrag{X}[c]{{\tiny Number of iterations}}
	\psfrag{Y}[c]{{\tiny Error}}
	\includegraphics[width=0.45\linewidth]{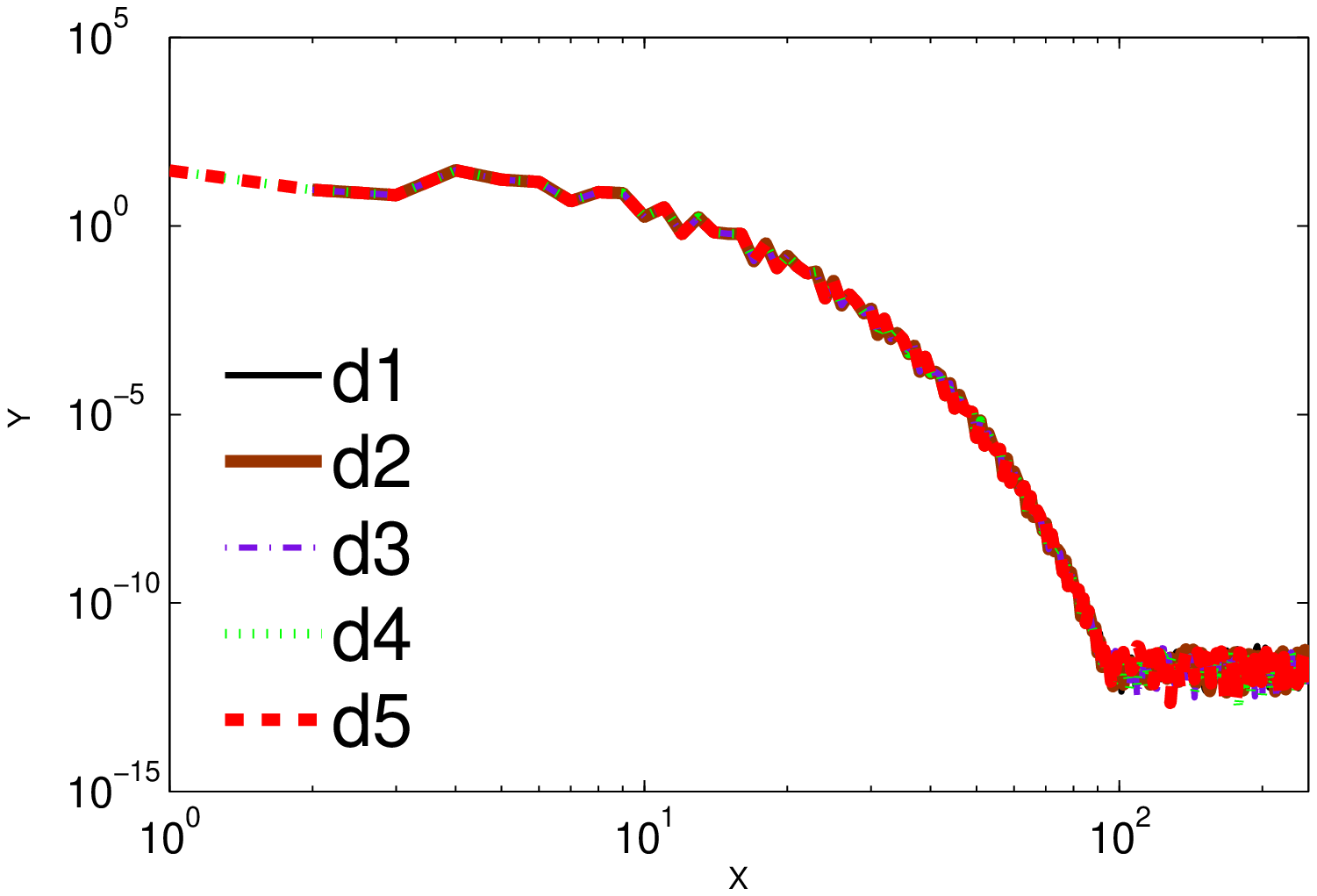}}
\subfigure[]{
\psfrag{X}[c]{{\tiny Number of iterations}}
	\psfrag{Y}[c]{{\tiny Error}}
\psfrag{100qn1}{{\tiny $k_2=10^{-1}$, for $q$}}
\psfrag{100qn2}{{\tiny $k_2=10^{-2}$, for $q$}}
\psfrag{100Gn1}{{\tiny $k_1=10$, for $\gamma$}}
\psfrag{100Gn2}{{\tiny $k_1=10^2$, for $\gamma$}}
	\includegraphics[width=0.45\linewidth]{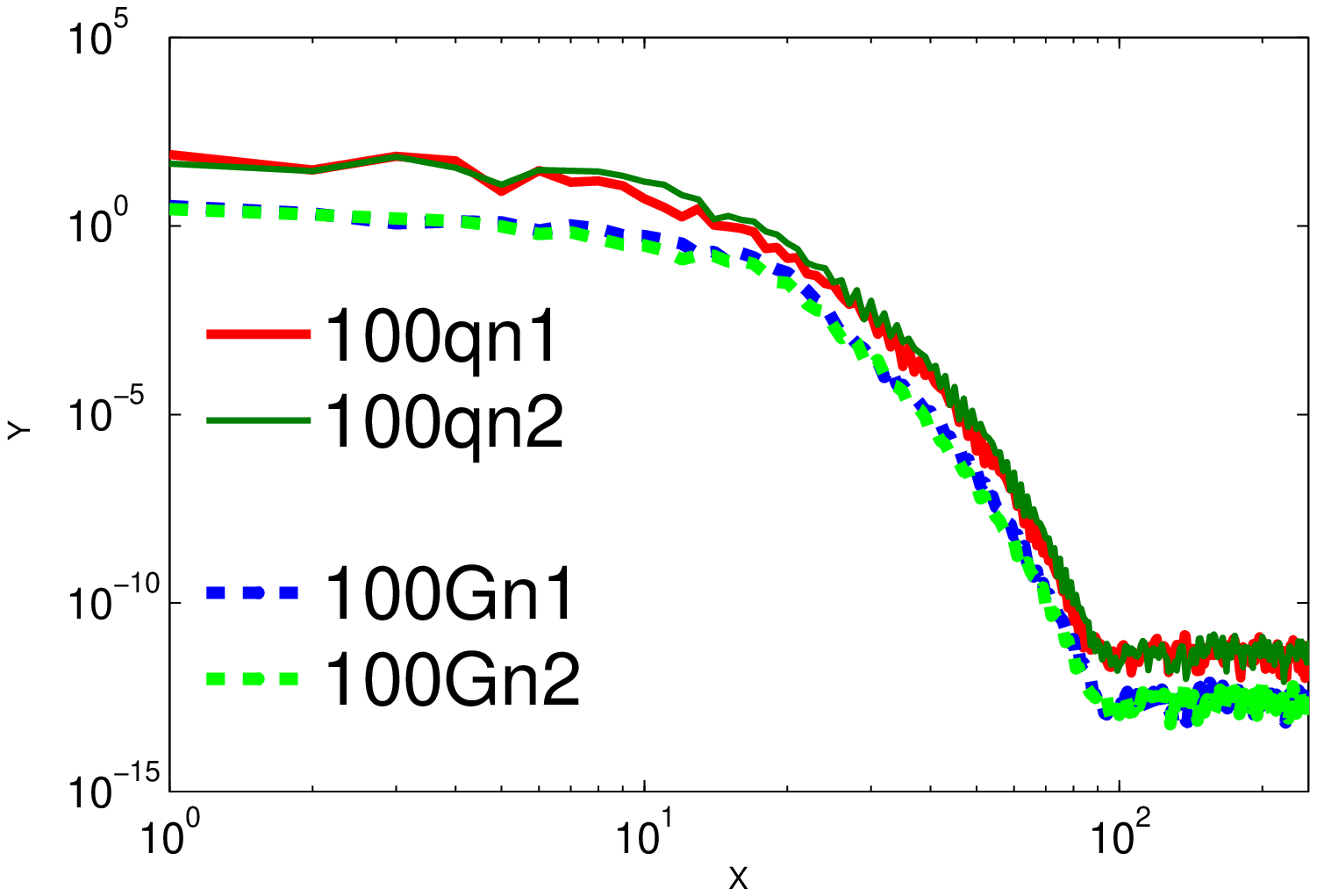}}
\\
\subfigure[]{
\psfrag{X}[c]{{\tiny Number of iterations}}
	\psfrag{Y}[c]{{\tiny Error}}
\psfrag{n1}{{\tiny $k_1=10$, $k_2=10^{-1}$}}
\psfrag{n2}{{\tiny $k_1=10^2$, $k_2=10^{-2}$}}
\psfrag{n3}{{\tiny $k_1=10^3$, $k_2=10^{-3}$}}
\psfrag{n5}{{\tiny $k_1=10^5$, $k_2=10^{-5}$}}
\psfrag{n7}{{\tiny $k_1=10^7$, $k_2=10^{-7}$}}
	\hspace{-0.5cm}\includegraphics[width=0.45\linewidth]{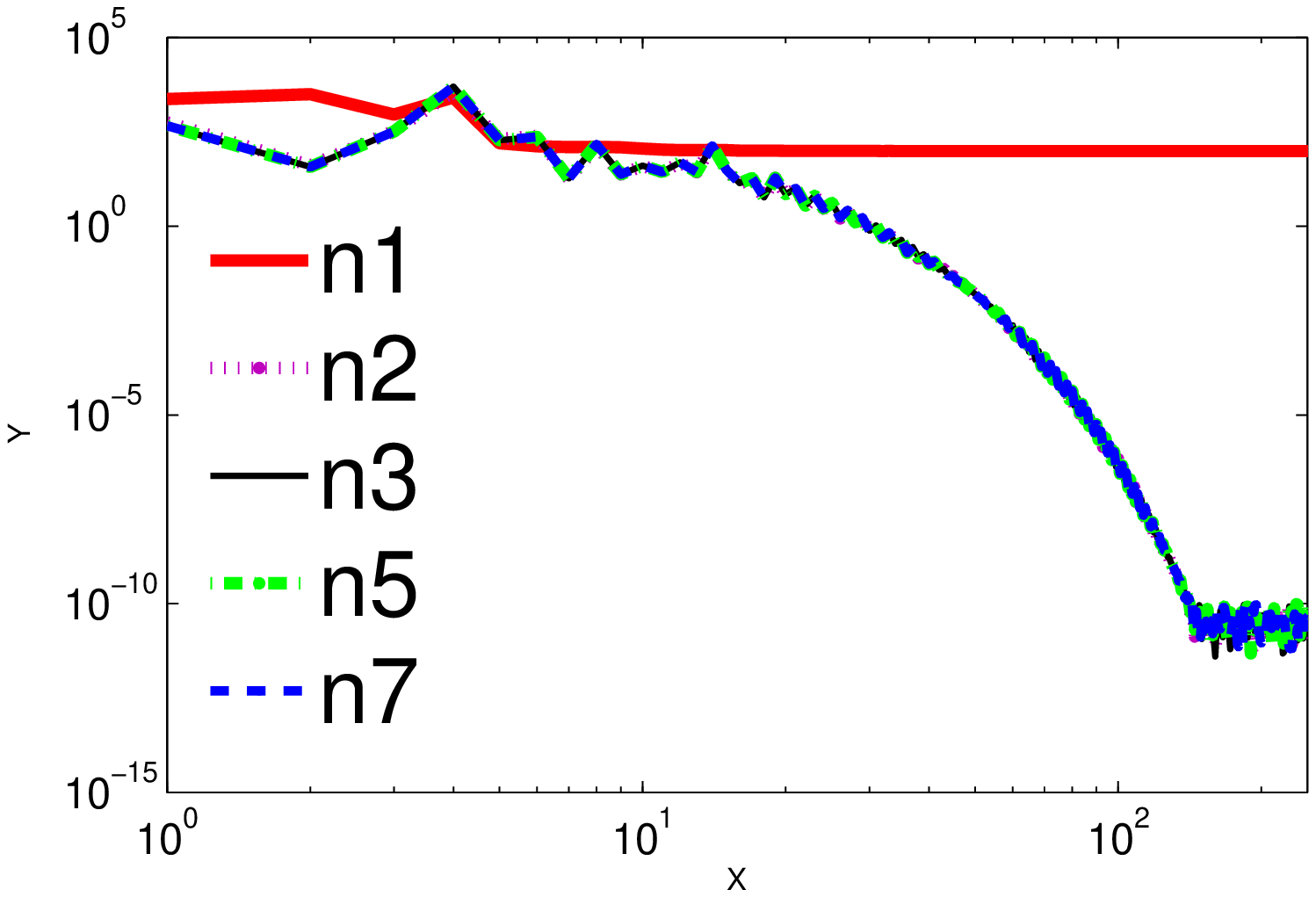}}
\subfigure[]{
\psfrag{X}[c]{{\tiny Number of iterations}}
	\psfrag{Y}[c]{{\tiny Error}}
\psfrag{qn1}{{\tiny  $k_2=10^{-1}$ for $q$}}
\psfrag{200qn2}{{\tiny  $k_2=10^{-2}$ for $q$}}
\psfrag{200Gn1}{{\tiny $k_1=10$ for $\gamma$}}
\psfrag{200Gn2}{{\tiny $k_1=10^2$ for $\gamma$}}
	\includegraphics[width=0.45\linewidth]{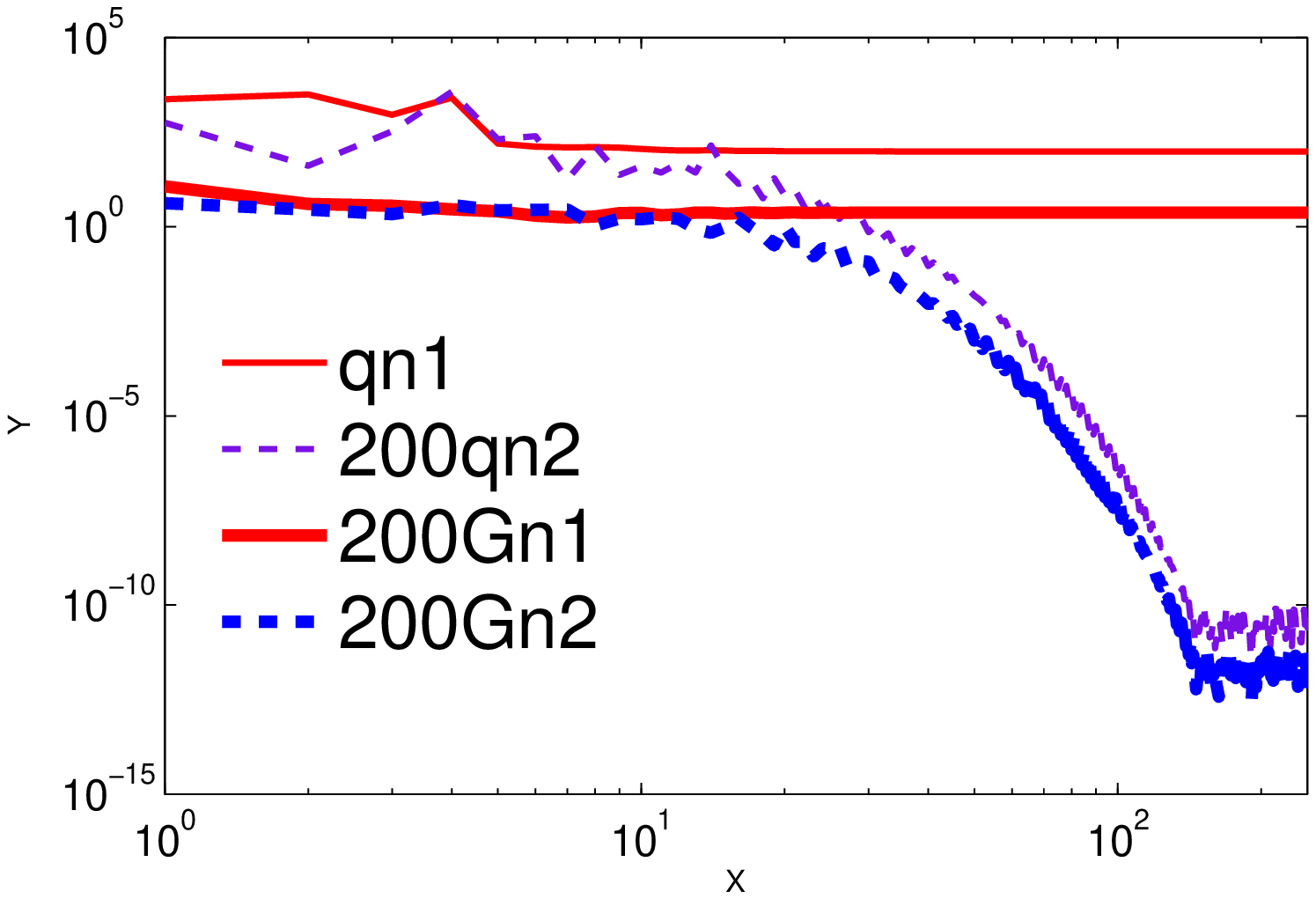}}
\end{center}
\caption{\label{Puel1} Dependence of the errors $\|J_{k_1}/|\nabla u|^2-\gamma\|_{L_\infty}$ and $\|j_{k_2}/|\nabla u|^2-q\|_{L_\infty}$  on the values of $k_1$, $k_2$. (a) Case of 50 points on the boundary.  (b) Case of 100 points on the boundary. (c) Case of 200 points on the boundary. \textcolor{black}{(d) Case of 400 points on the boundary.}}\label{Puel2}
\end{figure}

This  better convergence for meshes with $50$ and $100$ boundary points can be illustrated by the following example. For all types of mesh we can perfectly reconstruct $\gamma$ and $q$ by the perturbative method if one of the chosen frequency (for the reconstruction of  $\gamma$) is big enough and the second frequency (for the reconstruction of  $q$) is small enough. In the previous examples, the frequencies were chosen equal to $k_1=\pi\times10^3$ and $k_2=\pi\times10^{-3}$. During our numeric simulations, we have noticed that the smaller $|k_1-k_2|$ becomes, less efficient the convergence.
More precisely, the algorithm does not converge for   $|k_1-k_2|\le 10$ for the case of meshes with $200$ and $400$ boundary points (see Figure~(\ref{Puel1})).

\begin{figure}[h!]
\begin{center}
\subfigure[]{
\psfrag{m11}{{\tiny $m=1$}}
\psfrag{m21}{{\tiny $m=2$}}
\psfrag{m31}{{\tiny $m=3$}}
\psfrag{m51}{{\tiny $m=5$}}
\psfrag{m71}{{\tiny $m=7$}}
	\psfrag{iter}[c]{{\tiny Number of iterations}}
	\psfrag{E}[c]{{\tiny $\max | \tilde{u}_1|^2$}}
	\includegraphics[width=0.45\linewidth]{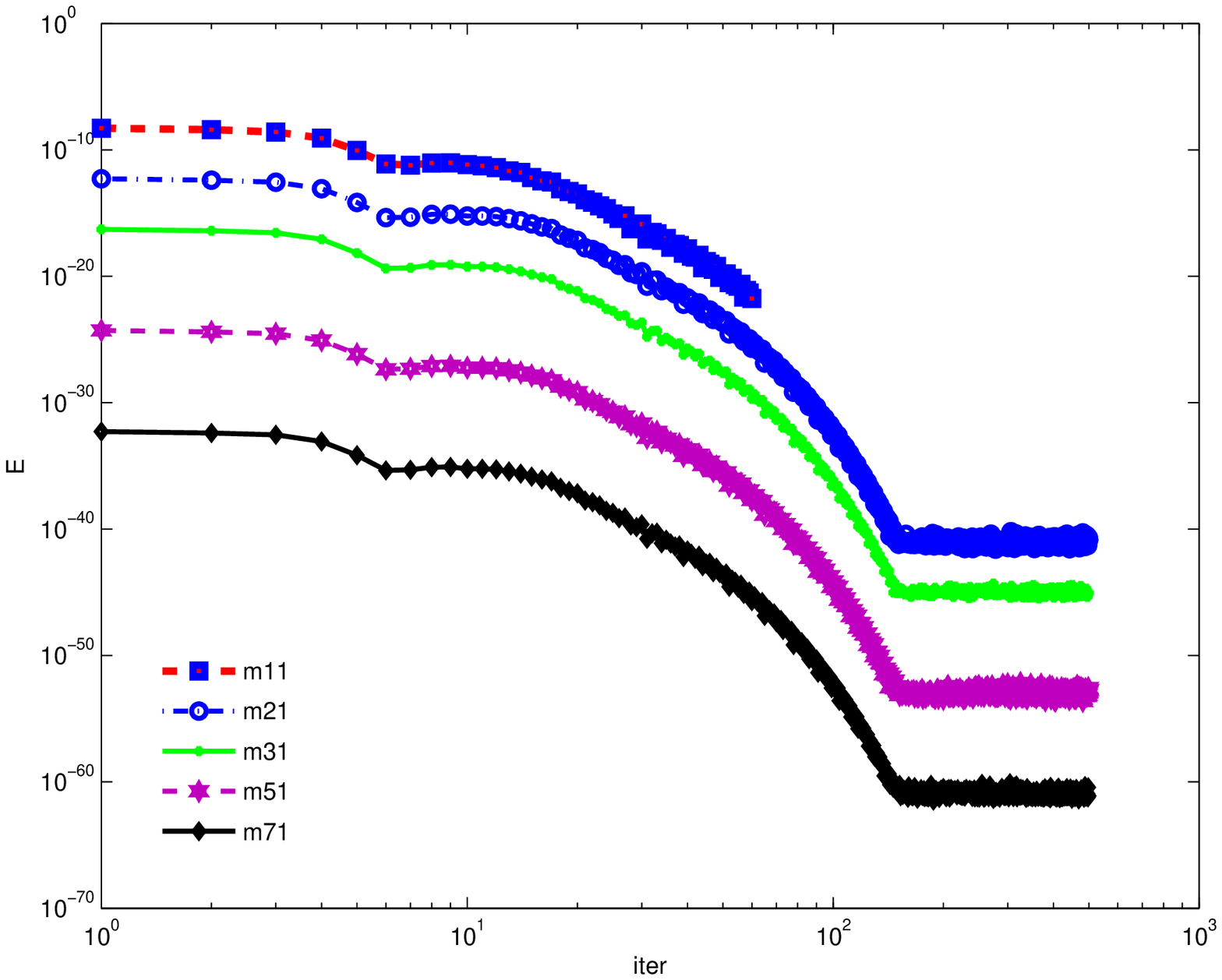}}
\subfigure[]{
\psfrag{iter}[c]{{\tiny Number of iterations}}
	\psfrag{E}[c]{{\tiny $\max | \tilde{u}_1|^2$}}
\psfrag{m11}{{\tiny $m=1$}}
\psfrag{m21}{{\tiny $m=2$}}
\psfrag{m31}{{\tiny $m=3$}}
\psfrag{m51}{{\tiny $m=5$}}
\psfrag{m71}{{\tiny $m=7$}}
	\includegraphics[width=0.45\linewidth]{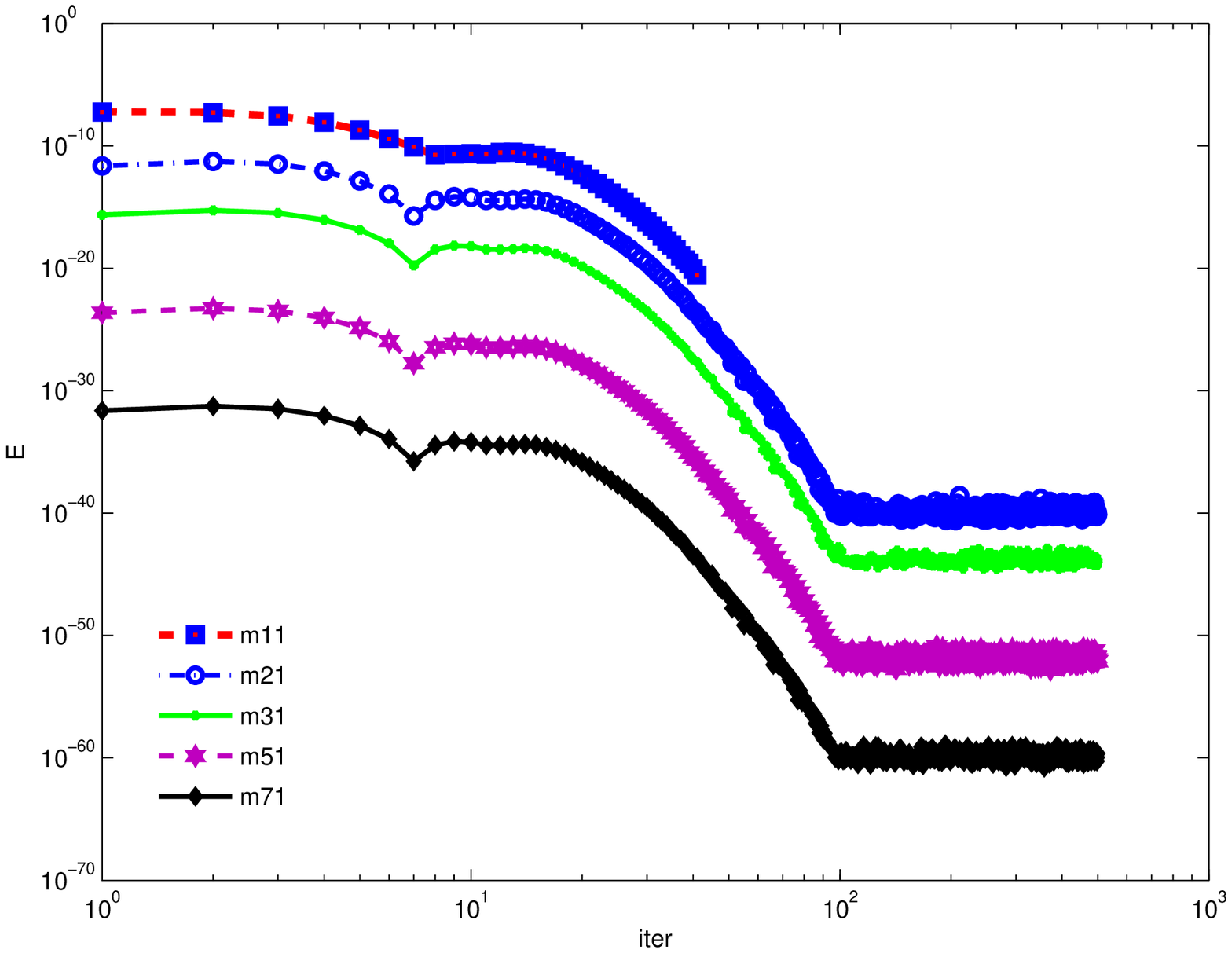}}
\\
\subfigure[]{
\psfrag{iter}[c]{{\tiny Number of iterations}}
	\psfrag{E}[c]{{\tiny $\max |\tilde{u}_1|^2$}}
\psfrag{m11}{{\tiny $m=1$}}
\psfrag{m21}{{\tiny $m=2$}}
\psfrag{m31}{{\tiny $m=3$}}
\psfrag{m51}{{\tiny $m=5$}}
\psfrag{m71}{{\tiny $m=7$}}
	\hspace{-0.5cm}\includegraphics[width=0.45\linewidth]{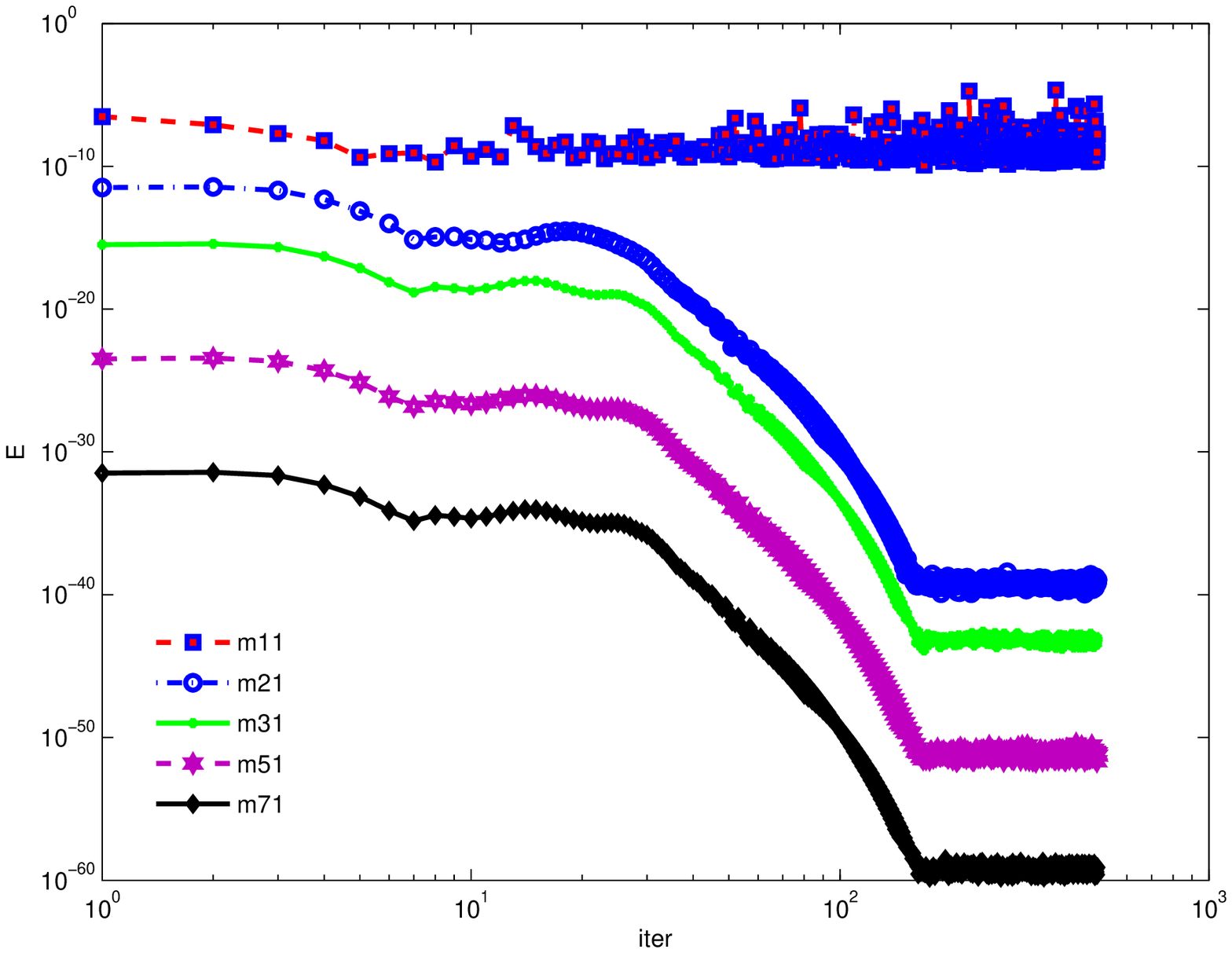}}
\end{center}
\caption{\label{u1O1} \textcolor{black}{Plot of   $\max | \tilde{u}_1|^2$ versus the number of iterations for different values of $k_1=\pi \times 10^m$, where $\tilde{u}_1$ is the corrector in the reconstruction of $\gamma$.   (a) Case of 50 points on the boundary.  (b) Case of 100 points on the boundary. (c) Case of 200 points on the boundary.}}\label{u1O1}
\end{figure}

\begin{figure}[h!]
\begin{center}
\subfigure[]{
\psfrag{m11}{{\tiny $m=1$}}
\psfrag{m21}{{\tiny $m=2$}}
\psfrag{m31}{{\tiny $m=3$}}
\psfrag{m51}{{\tiny $m=5$}}
\psfrag{m71}{{\tiny $m=7$}}
	\psfrag{iter}[c]{{\tiny Number of iterations}}
	\psfrag{E}[c]{{\tiny $\min | \nabla u|^2$}}
	\includegraphics[width=0.45\linewidth]{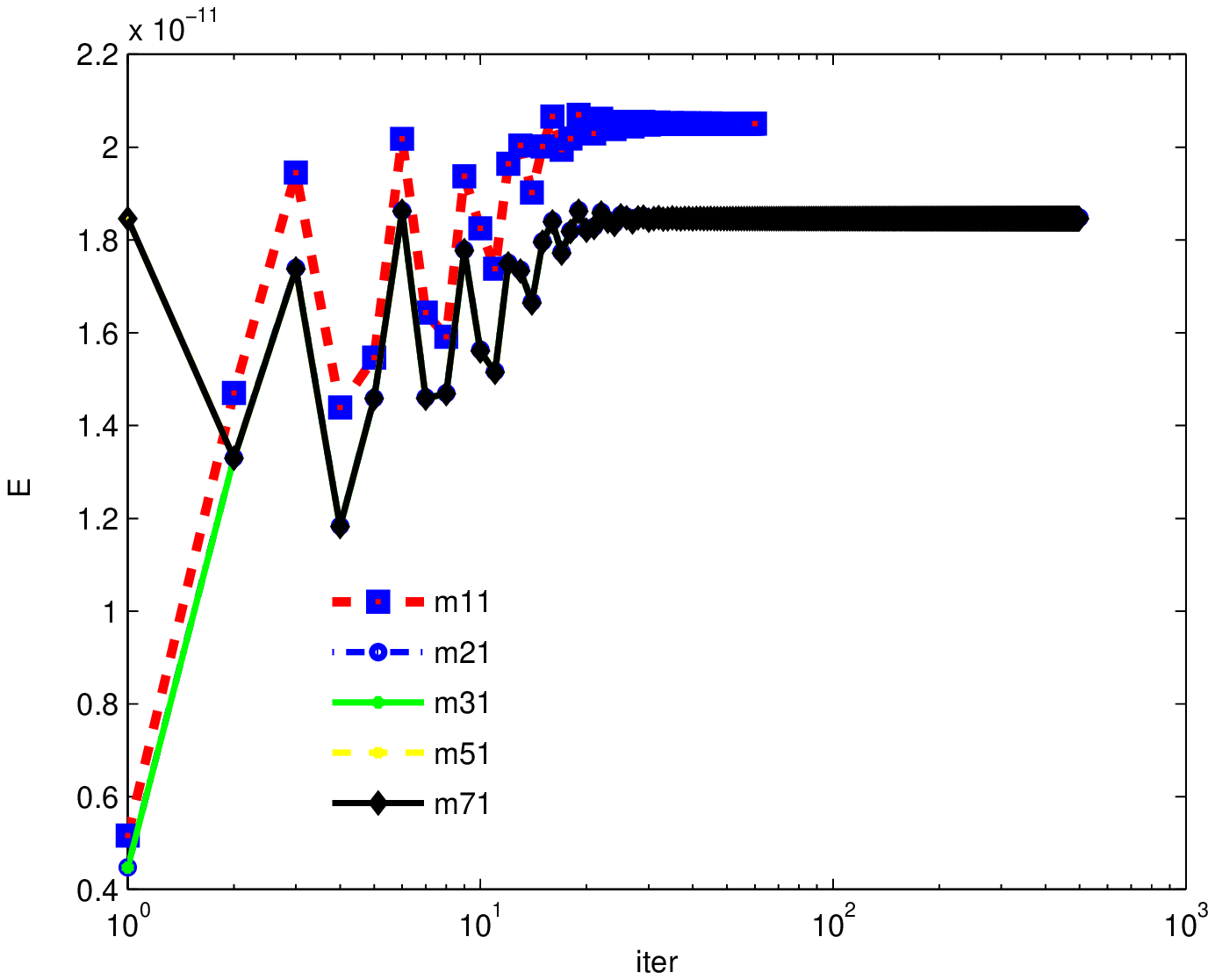}}
\subfigure[]{
\psfrag{m11}{{\tiny $m=1$}}
\psfrag{m21}{{\tiny $m=2$}}
\psfrag{m31}{{\tiny $m=3$}}
\psfrag{m51}{{\tiny $m=5$}}
\psfrag{m71}{{\tiny $m=7$}}
	\psfrag{iter}[c]{{\tiny Number of iterations}}
	\psfrag{E}[c]{{\tiny $\min |u|^2$}}
	\includegraphics[width=0.45\linewidth]{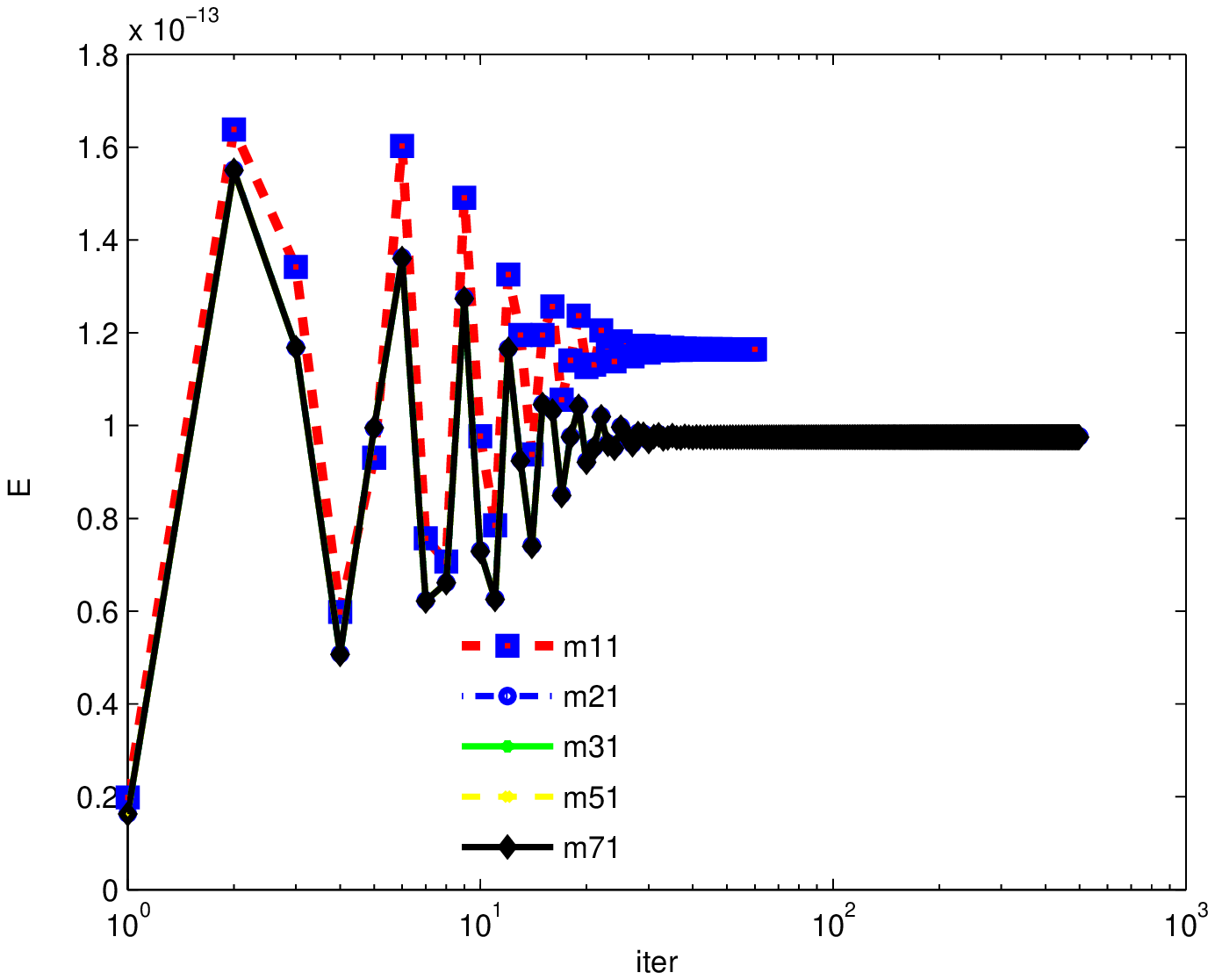}}
\\
\subfigure[]{
\psfrag{iter}[c]{{\tiny Number of iterations}}
	\psfrag{E}[c]{{\tiny $\min | \nabla u|^2$}}
\psfrag{m11}{{\tiny $m=1$}}
\psfrag{m21}{{\tiny $m=2$}}
\psfrag{m31}{{\tiny $m=3$}}
\psfrag{m51}{{\tiny $m=5$}}
\psfrag{m71}{{\tiny $m=7$}}
	\includegraphics[width=0.45\linewidth]{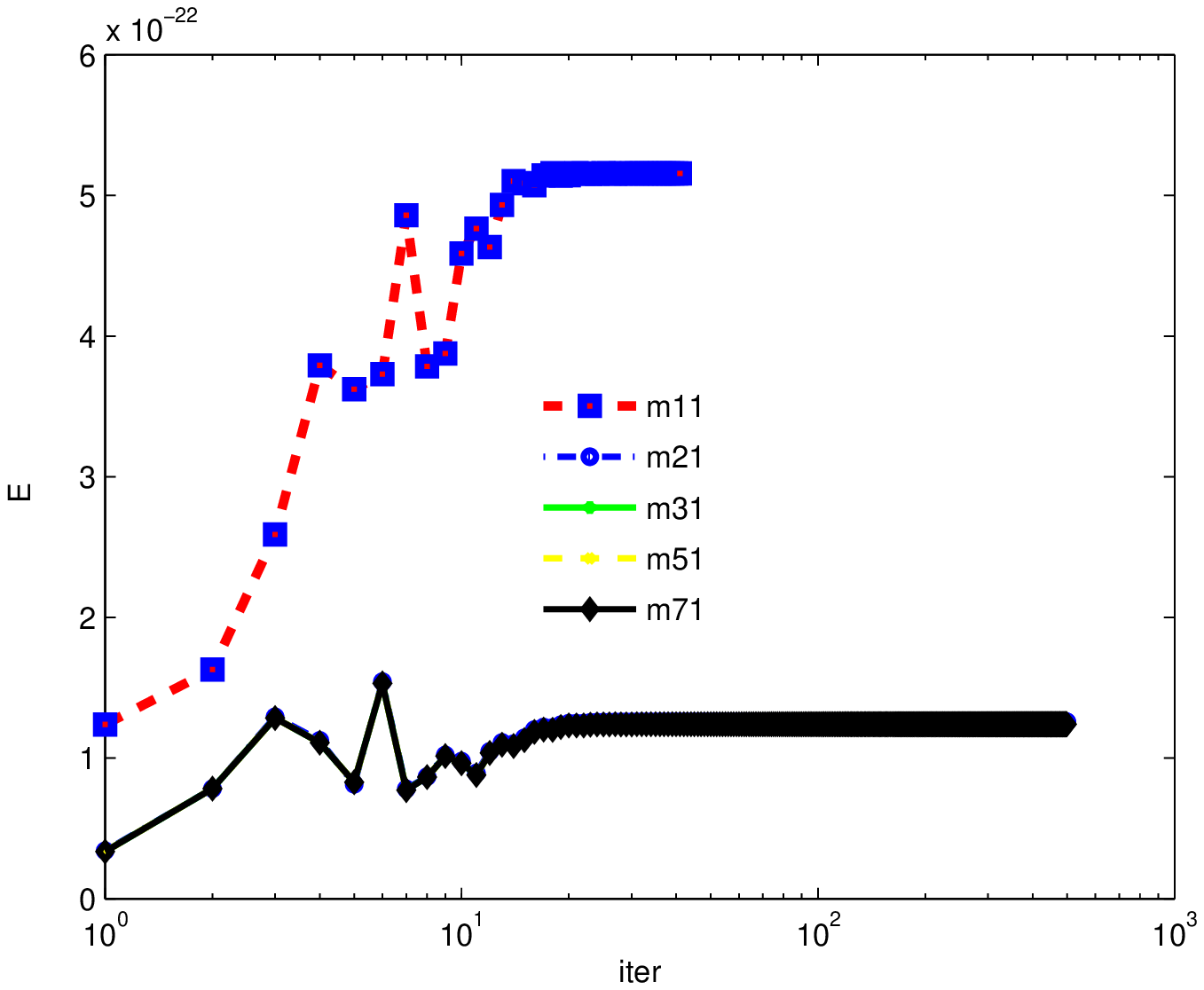}}
\subfigure[]{
\psfrag{iter}[c]{{\tiny Number of iterations}}
	\psfrag{E}[c]{{\tiny $\min |u|^2$}}
\psfrag{m11}{{\tiny $m=1$}}
\psfrag{m21}{{\tiny $m=2$}}
\psfrag{m31}{{\tiny $m=3$}}
\psfrag{m51}{{\tiny $m=5$}}
\psfrag{m71}{{\tiny $m=7$}}
	\includegraphics[width=0.45\linewidth]{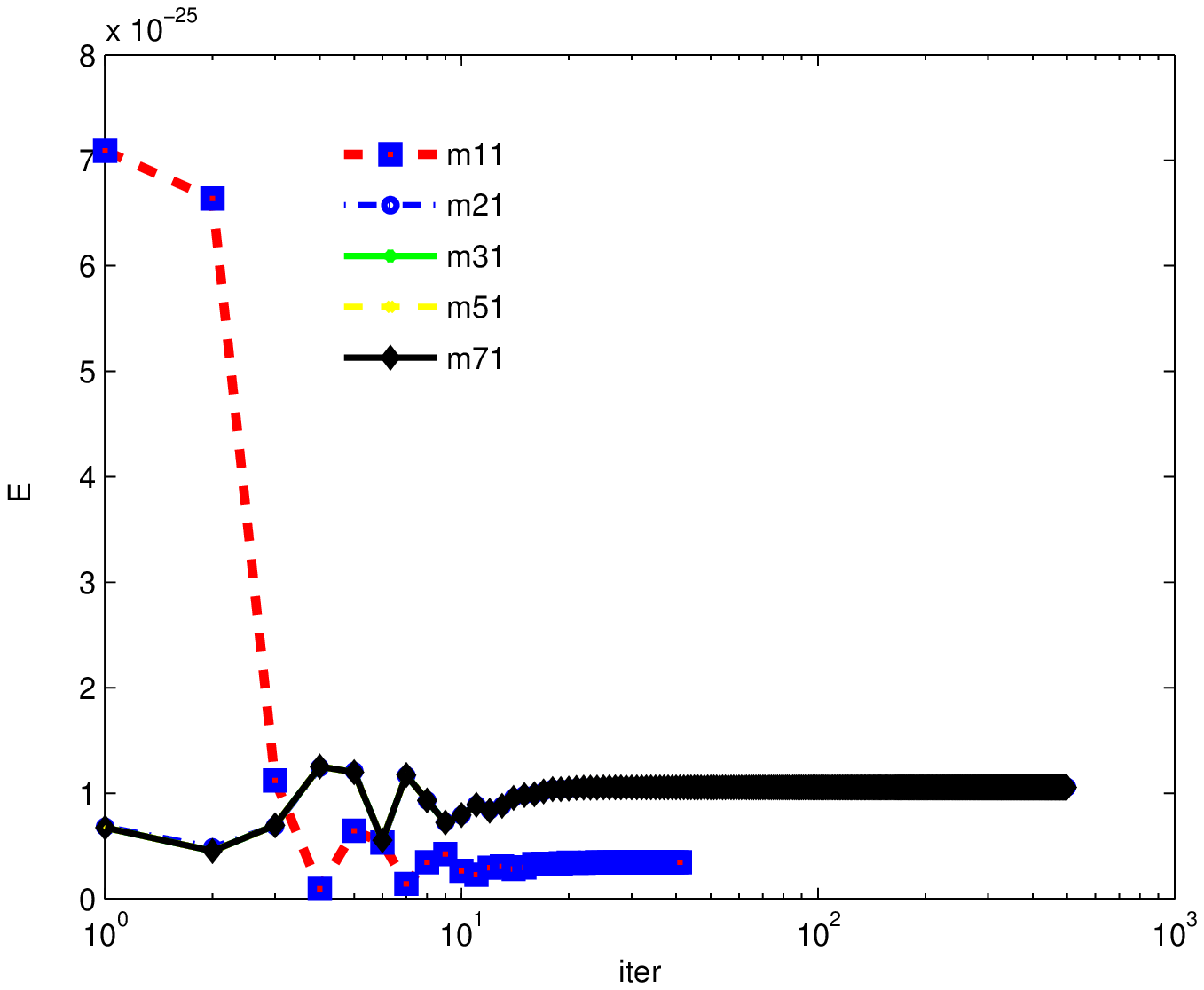}}
\\
\subfigure[]{
\psfrag{iter}[c]{{\tiny Number of iterations}}
	\psfrag{E}[c]{{\tiny $\min | \nabla u|^2$}}
\psfrag{m11}{{\tiny $m=1$}}
\psfrag{m21}{{\tiny $m=2$}}
\psfrag{m31}{{\tiny $m=3$}}
\psfrag{m51}{{\tiny $m=5$}}
\psfrag{m71}{{\tiny $m=7$}}
	\hspace{-0.5cm}\includegraphics[width=0.45\linewidth]{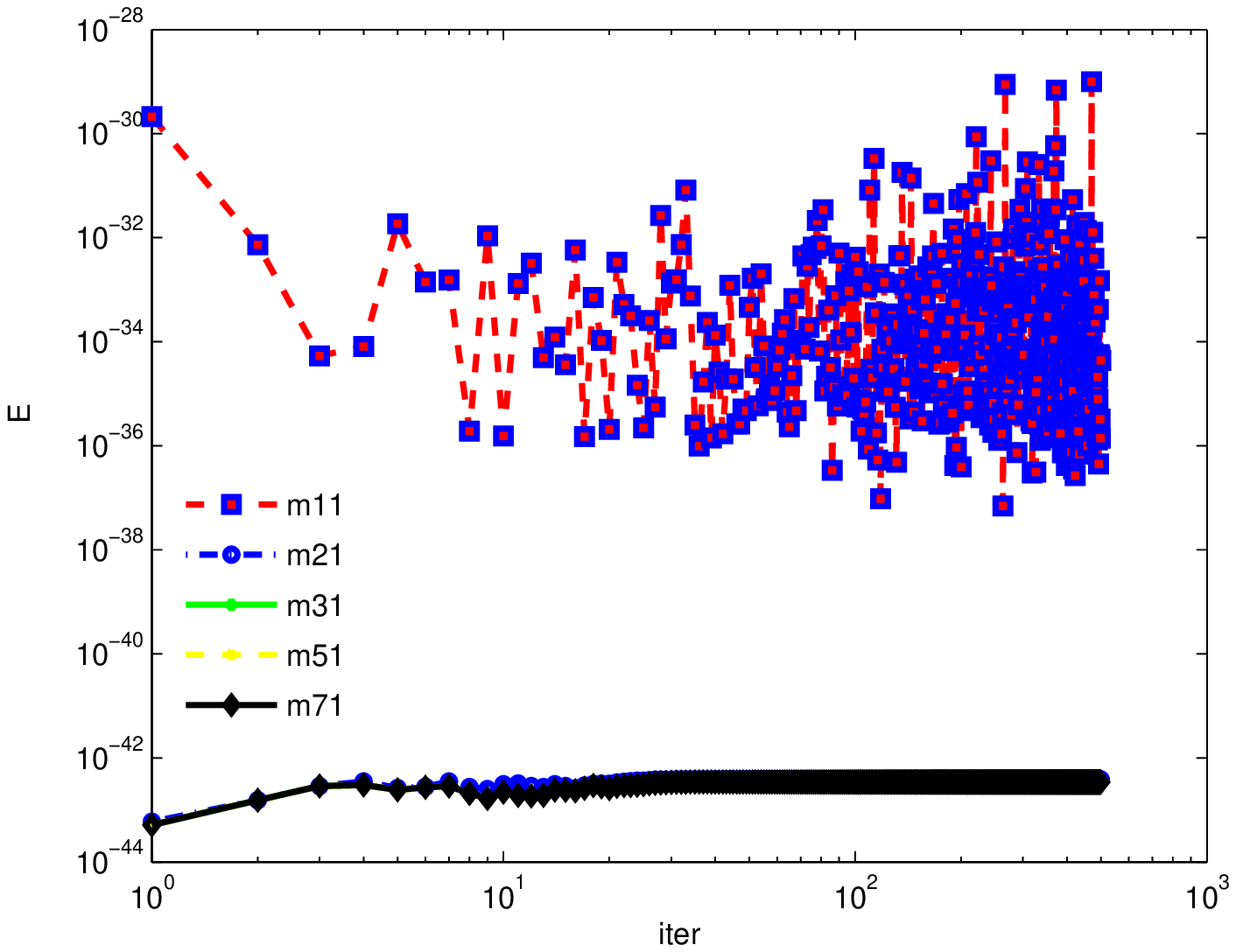}}
\subfigure[]{
\psfrag{iter}[c]{{\tiny Number of iterations}}
	\psfrag{E}[c]{{\tiny $\min | u|^2$}}
\psfrag{m11}{{\tiny $m=1$}}
\psfrag{m21}{{\tiny $m=2$}}
\psfrag{m31}{{\tiny $m=3$}}
\psfrag{m51}{{\tiny $m=5$}}
\psfrag{m71}{{\tiny $m=7$}}
	\hspace{-0.5cm}\includegraphics[width=0.45\linewidth]{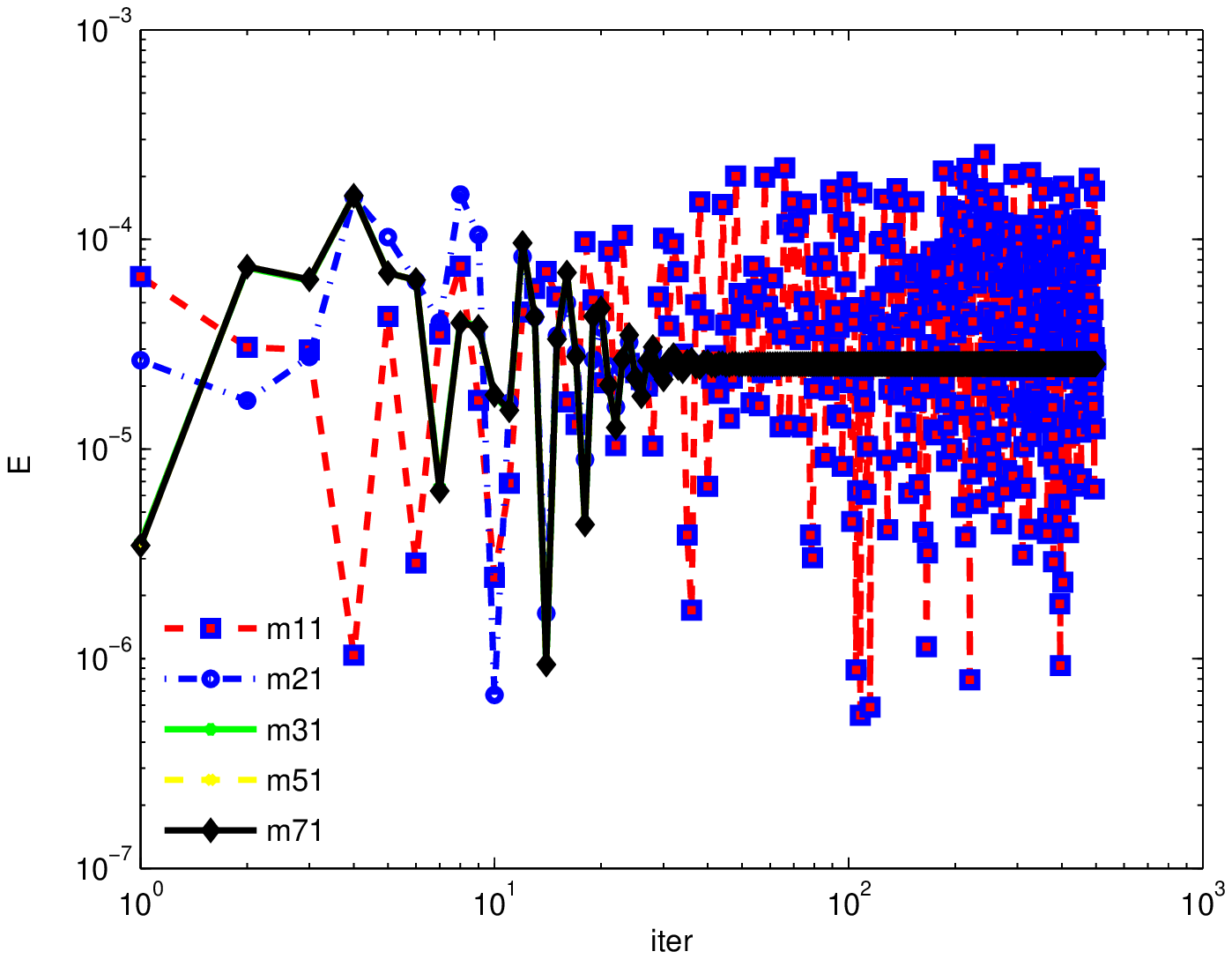}}
\end{center}
\caption{\label{GraduO1} \textcolor{black}{ Plot of   $\min | \nabla u|^2$ and $\min | u|^2$ versus the number of iterations for different values of $k_1=\pi \times 10^{m}$, where $u$ is the numeric solution of the Helmholtz problem.   (a) and (b) Case of 50 points on the boundary.  (c) and (d) Case of 100 points on the boundary. (e) and (f) Case of 200 points on the boundary.}}\label{GraduO1}
\end{figure}

\begin{figure}[h!]
\begin{center}
\subfigure[]{
\psfrag{m11}{{\tiny $m=1$}}
\psfrag{m21}{{\tiny $m=2$}}
\psfrag{m31}{{\tiny $m=3$}}
\psfrag{m51}{{\tiny $m=5$}}
\psfrag{m71}{{\tiny $m=7$}}
	\psfrag{iter}[c]{{\tiny Number of iterations}}
	\psfrag{E}[c]{{\tiny $\min | \nabla u|^2$}}
	\includegraphics[width=0.45\linewidth]{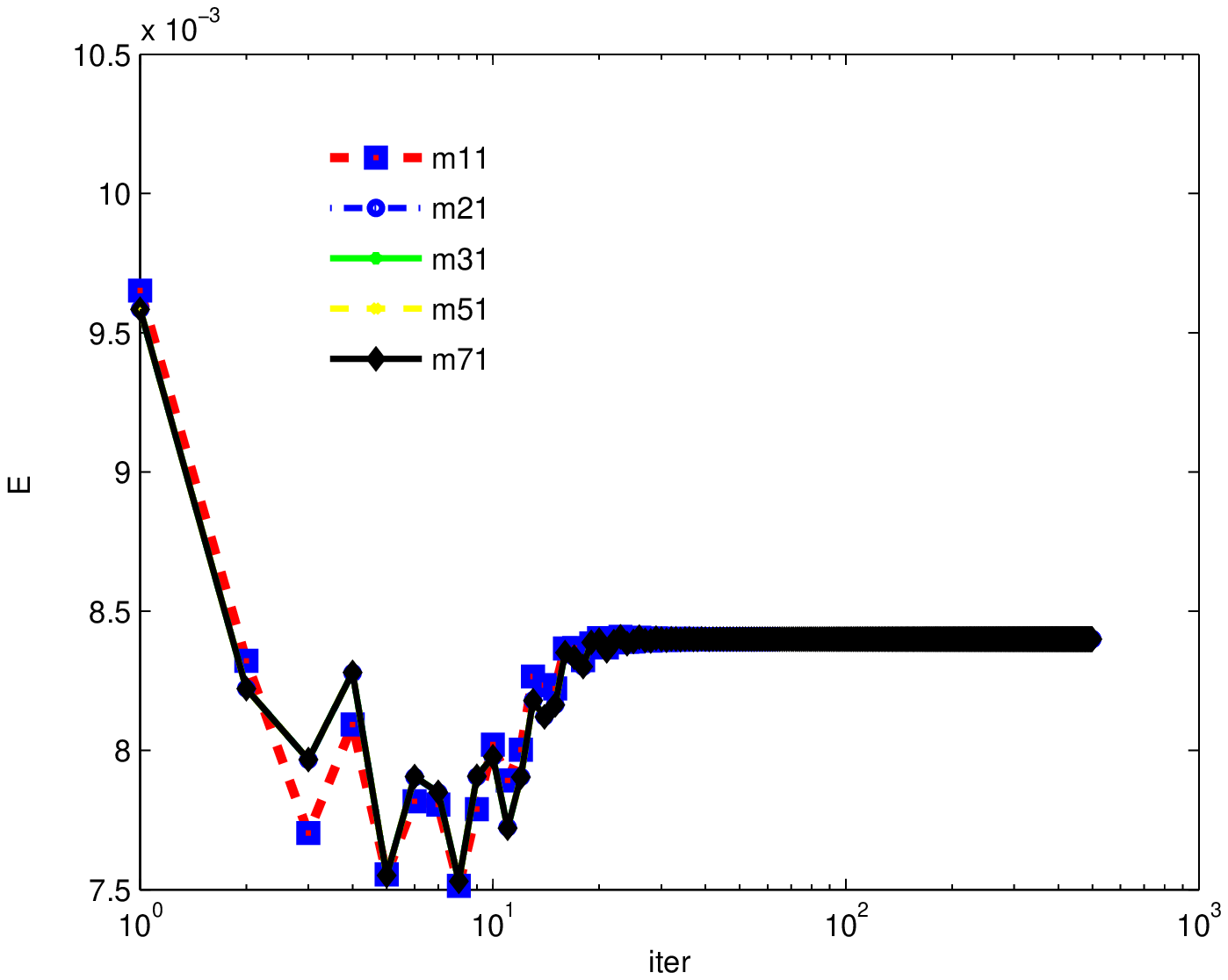}}
\subfigure[]{
\psfrag{m11}{{\tiny $m=1$}}
\psfrag{m21}{{\tiny $m=2$}}
\psfrag{m31}{{\tiny $m=3$}}
\psfrag{m51}{{\tiny $m=5$}}
\psfrag{m71}{{\tiny $m=7$}}
	\psfrag{iter}[c]{{\tiny Number of iterations}}
	\psfrag{E}[c]{{\tiny $\min |u|^2$}}
	\includegraphics[width=0.45\linewidth]{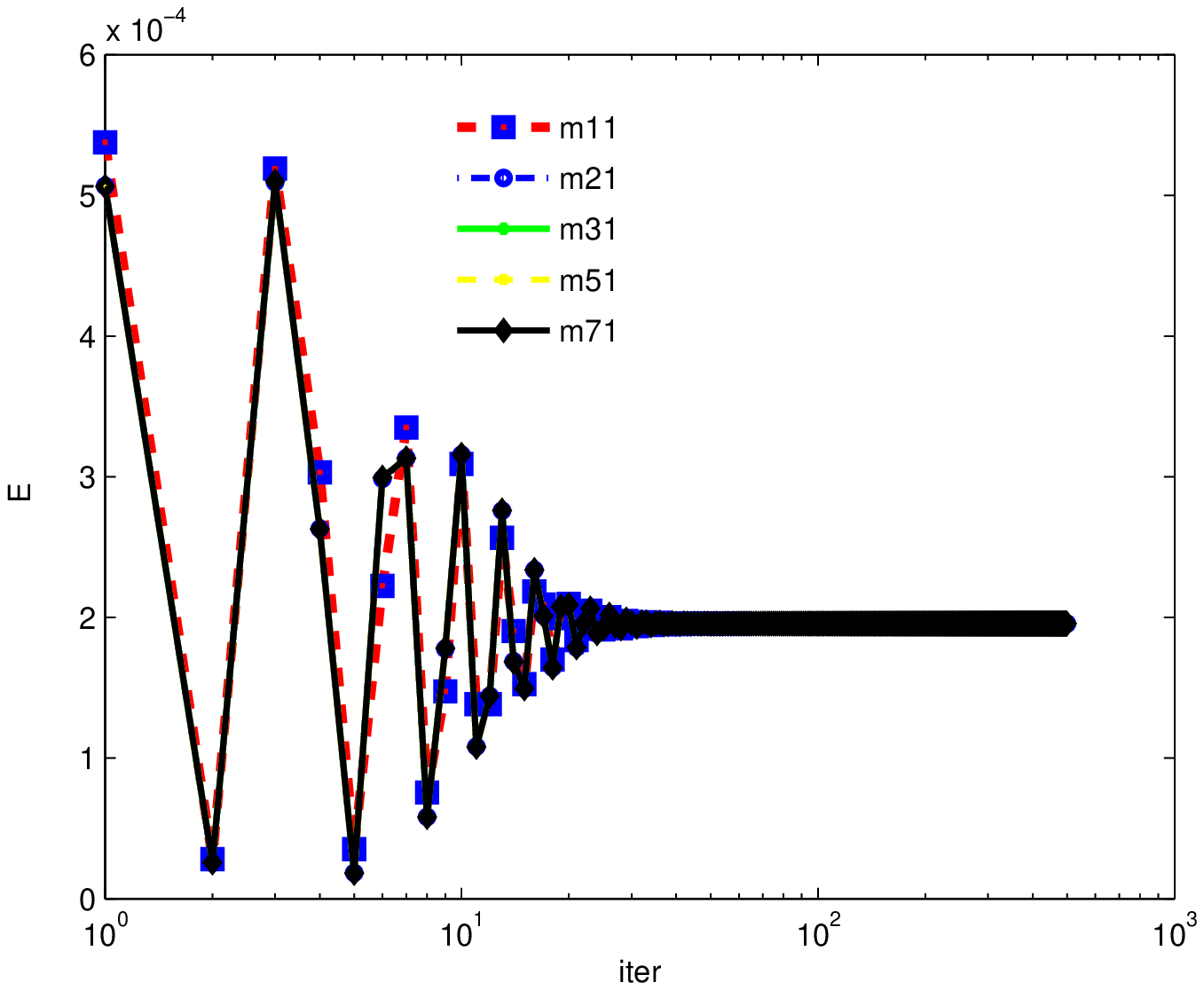}}
\\
\subfigure[]{
\psfrag{iter}[c]{{\tiny Number of iterations}}
	\psfrag{E}[c]{{\tiny $\min | \nabla u|^2$}}
\psfrag{m11}{{\tiny $m=1$}}
\psfrag{m21}{{\tiny $m=2$}}
\psfrag{m31}{{\tiny $m=3$}}
\psfrag{m51}{{\tiny $m=5$}}
\psfrag{m71}{{\tiny $m=7$}}
	\includegraphics[width=0.45\linewidth]{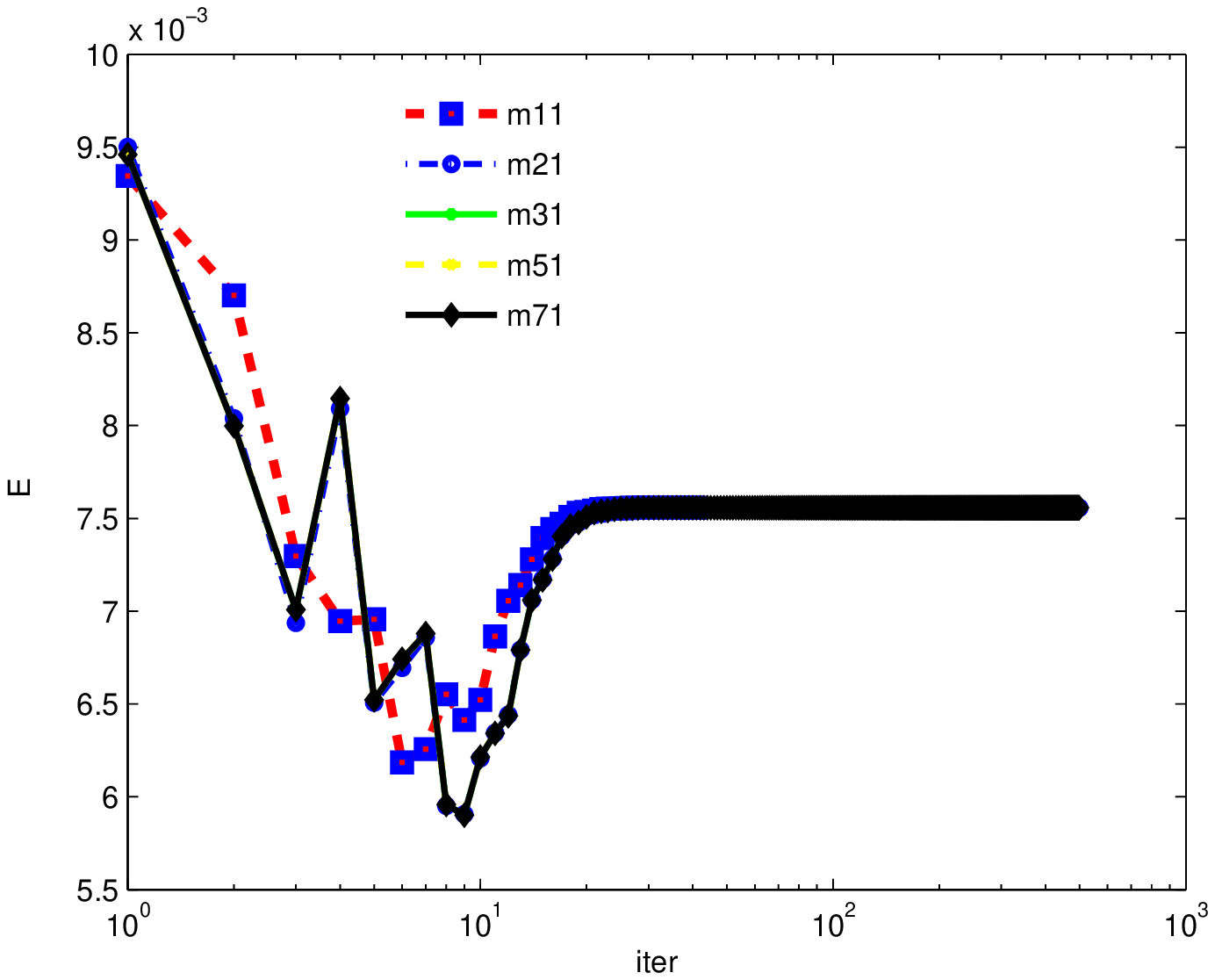}}
\subfigure[]{
\psfrag{iter}[c]{{\tiny Number of iterations}}
	\psfrag{E}[c]{{\tiny $\min |u|^2$}}
\psfrag{m11}{{\tiny $m=1$}}
\psfrag{m21}{{\tiny $m=2$}}
\psfrag{m31}{{\tiny $m=3$}}
\psfrag{m51}{{\tiny $m=5$}}
\psfrag{m71}{{\tiny $m=7$}}
	\includegraphics[width=0.45\linewidth]{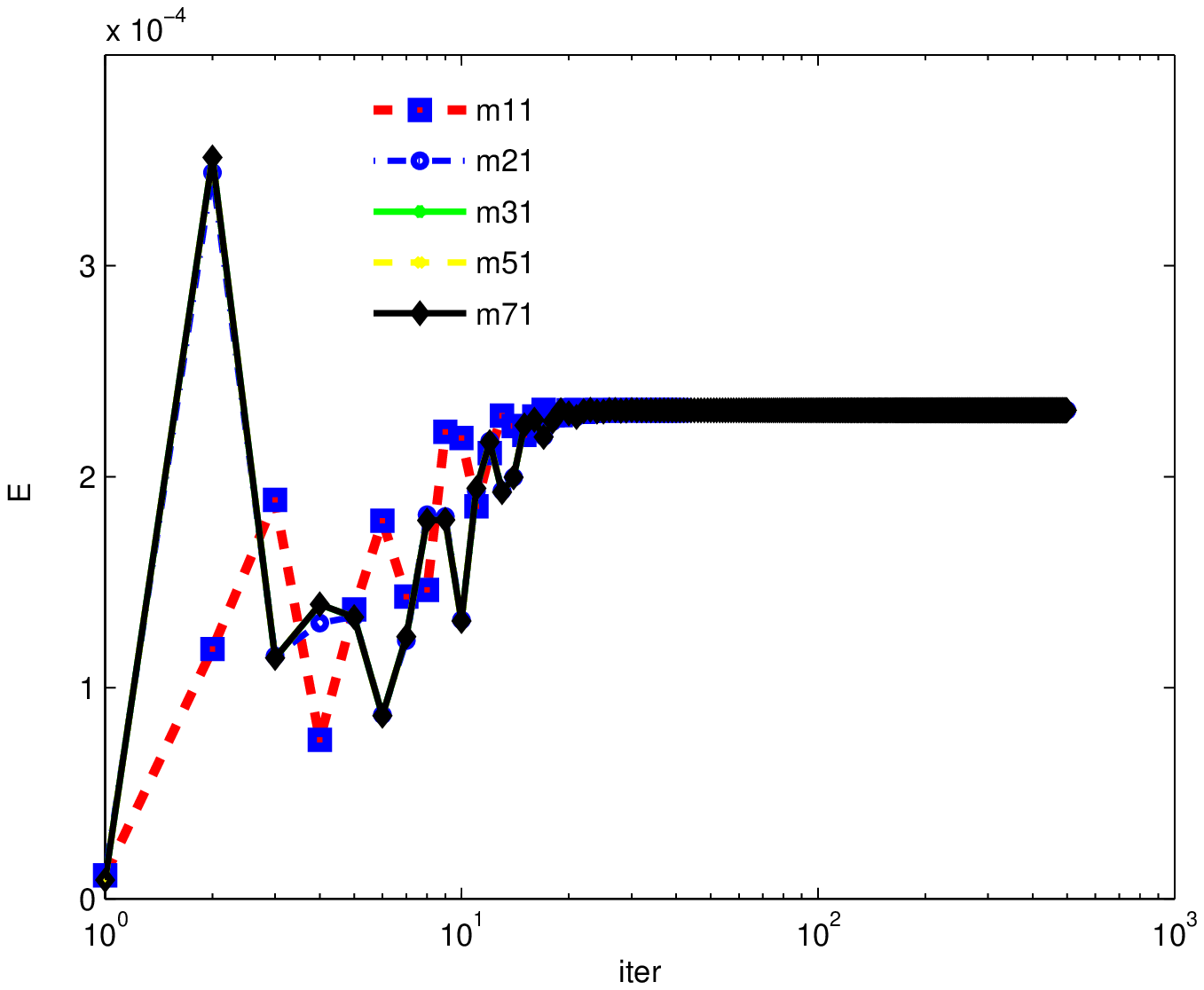}}
\\
\subfigure[]{
\psfrag{iter}[c]{{\tiny Number of iterations}}
	\psfrag{E}[c]{{\tiny $\min | \nabla u|^2$}}
\psfrag{m11}{{\tiny $m=1$}}
\psfrag{m21}{{\tiny $m=2$}}
\psfrag{m31}{{\tiny $m=3$}}
\psfrag{m51}{{\tiny $m=5$}}
\psfrag{m71}{{\tiny $m=7$}}
	\hspace{-0.5cm}\includegraphics[width=0.45\linewidth]{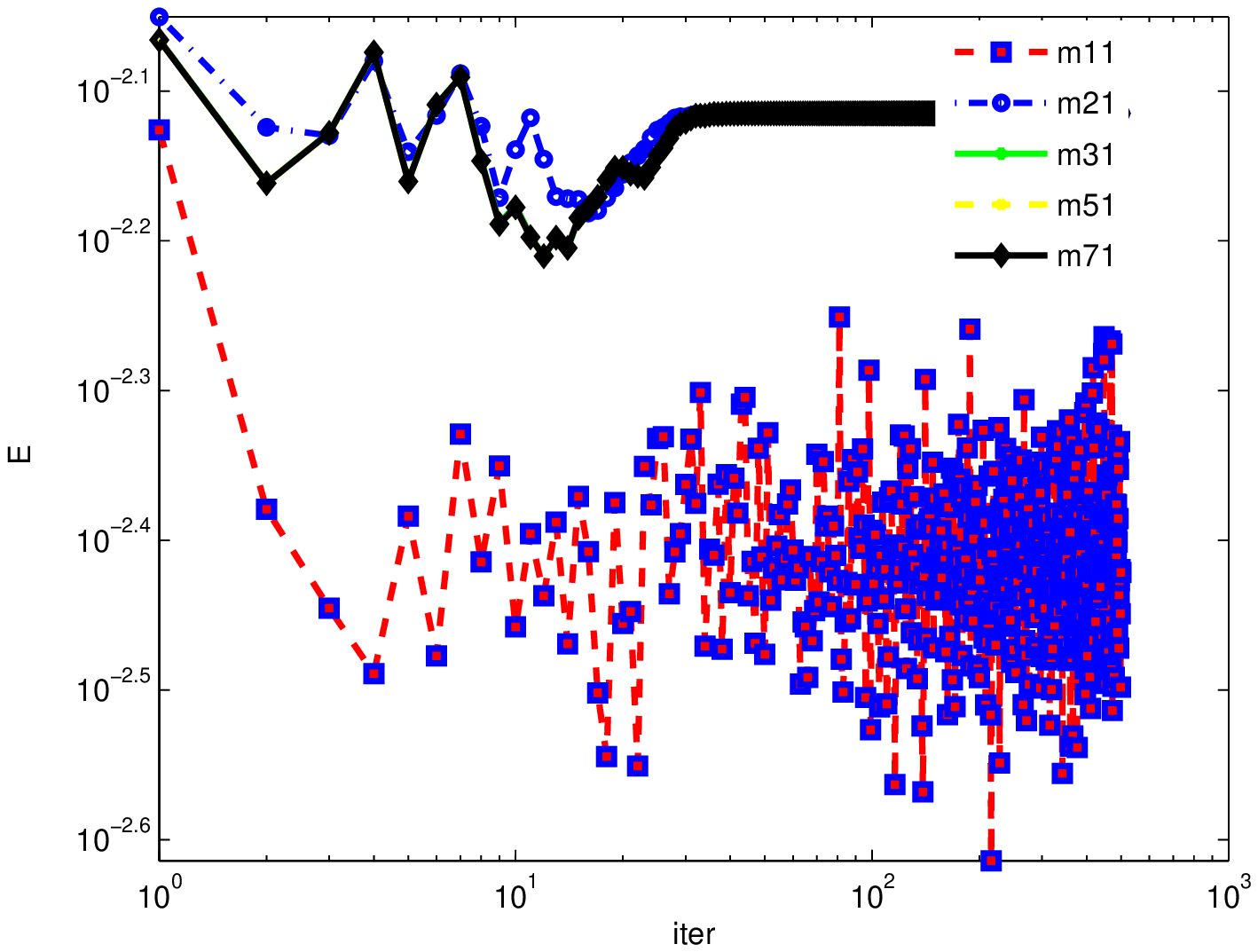}}
\subfigure[]{
\psfrag{iter}[c]{{\tiny Number of iterations}}
	\psfrag{E}[c]{{\tiny $\min | u|^2$}}
\psfrag{m11}{{\tiny $m=1$}}
\psfrag{m21}{{\tiny $m=2$}}
\psfrag{m31}{{\tiny $m=3$}}
\psfrag{m51}{{\tiny $m=5$}}
\psfrag{m71}{{\tiny $m=7$}}
	\hspace{-0.5cm}\includegraphics[width=0.45\linewidth]{./MinUOmega2Mesh200enFoctionDeITERtoutN.eps}}
\end{center}
\caption{\label{uO2} \textcolor{black}{ Plot of $\min | \nabla u|^2$ and  $\min | u|^2$ versus the number of iterations for different values of $k_2=\pi \times 10^{-m}$, where $u$ is the numeric solution of the Helmholtz problem.  (a) and (b) Case of 50 points on the boundary.  (c) and (d) Case of 100 points on the boundary. (e) and (f) Case of 200 points on the boundary.}}\label{uO2}
\end{figure}

\textcolor{black}{Let us analyse the explanation of these results. We notice that there are two necessary conditions to be satisfied to ensure the convergence of the algorithm by perturbations:}
\begin{enumerate}
 \item \textcolor{black}{Using the approximation $\gamma(x)=\frac{J(x)}{|\nabla u(x)|^2}$ and $q(x)=\frac{j(x)}{| u(x)|^2}$, we need to ensure that there exist $\delta_1>0$ and $\delta_2>0$ such that for each iteration step for $l=1, 2,\ldots$,  $|\nabla  u^l_{k_i}(x)|^2>\delta_1$ and $| u^l_{k_i}(x)|^2>\delta_2$ (where $k_i$ are the chosen frequencies). In other words, we need that the sequences $\{|\nabla u^l_{k_i}(x)|^2: l=1,2\ldots \}$ and $\{|u^l_{k_i}(x)|^2 : l=1,2\ldots\}$ have some uniform positive lower bound.}
\item \textcolor{black}{The corrector functions to update the initial guess for $\gamma$ and $q$ should be small enough ($|\tilde{u}_1^l|=\delta|u_1|\ll 1$) and for $l\to \infty$,  $|\tilde{u}_1^l|$ should tends to $0$.}
\end{enumerate}
{\color{black}Indeed, if the first condition does not hold, we have a division by zero and the algorithm has no any sense.
In the second condition, the smallness of the correctors functions $\tilde{u}_1^l$ is the basic assumption for deriving the approximate systems~(\ref{GU1cora})-(\ref{GU1corb}) and~(\ref{QU1cor}) which avoid all the terms of the second order on $\delta$. If $|\tilde{u}_1^l|$ is not small enough, we cannot do it any more  and  the solutions of  systems~(\ref{GU1cora})-(\ref{GU1corb}) and~(\ref{QU1cor}) have   no any sense.
Moreover, the algorithm converges if and only if $|\tilde{u}_1^l|\to 0$ for $l\to \infty$.

Figure~\ref{u1O1} shows the decay behaviour of the upper bound of $|\tilde{u}_1^l|^2$ for the corrector $\tilde{u}_1^l$ from the conductivity update algorithm (see system~(\ref{GU1cora})-(\ref{GU1corb})) for different frequencies and meshes. We observe that we have a good convergence corresponding to the logarithmic decay of $|\tilde{u}_1^l|^2$ for all frequencies and meshes with $50$ and $100$ boundary points, but we have a divergence result corresponding to the non-decay of  $|\tilde{u}_1^l|^2$ for the frequency $k_1=10 \pi$ and for the mesh with $200$ boundary points. The corrector function $\tilde{u}_1^l$ for the reconstruction of $q$ in our numerical tests for $k_2=\pi \times 10^{-m}$, $m=1, 2,3,5,7$, is equal to zero. This means that at each iteration step we update $q_0$ by  $\frac{j(x)}{|u^l_{k_i}(x)|^2}$.
 
To understand why for $k_1=10 \pi$ and $k_2=0.1 \pi$  the algorithm diverges for a mesh of $200$ boundary points and converges for a mesh of $50$ or $100$ boundary points, let us verify the first condition of the convergence. Figure~\ref{GraduO1} (respectively Figure~\ref{uO2}) shows the lower bounds of $|\nabla  u^l_{k_i}(x)|^2$ and $| u^l_{k_i}(x)|^2$ for different $k_1$ (respectively $k_2$) and for different meshes. We notice that we have for all cases, except the case for $k_1=10 \pi$, $k_2=0.1 \pi$ and for the mesh with $200$ boundary points, that the sequences $\{\min\limits_x|\nabla u^l_{k_i}(x)|^2: l=1,2\ldots \}$ and $\{\min\limits_x|u^l_{k_i}(x)|^2 : l=1,2\ldots\}$ converge for $l\to \infty$ to a positive constant.
Therefore, we see that for $k_1=10 \pi$ and $k_2=0.1 \pi$ we obtain a divergence of the quantities $\frac{J(x)}{|\nabla u^l_{k_i}(x)|^2}$ and $\frac{j(x)}{| u^l_{k_i}(x)|^2}$.
The divergence does not take place for a small number of boundary points because of a lower order of the precision. For example, we notice that with the growth of the number of boundary points (\textit{i.e.} with the growth of the precision) the limits of the sequences $\{\min\limits_x|\nabla u^l_{k_i}(x)|^2: l=1,2\ldots \}$ and $\{\min\limits_x|u^l_{k_i}(x)|^2 : l=1,2\ldots\}$ becomes more and more smaller, as  illustrated on Figure~\ref{mmin} for $\min\limits_x|u^l_{k_1}|^2$. In the case the mesh of $400$ boundary points, the divergence stops the numeric test by the error of the division by zero.
}


\begin{figure}[h!]
\begin{center}
 \psfrag{iter}[c]{{\tiny Number of iterations}}
	\psfrag{E}[c]{{\tiny $\min | u|^2$}}
\psfrag{m11}{{\tiny $m=1$, mesh $50$}}
\psfrag{m21}{{\tiny $m=2$,mesh $50$}}
\psfrag{m31}{{\tiny $m=3$, mesh $50$}}
\psfrag{m51}{{\tiny $m=5$,mesh $50$}}
\psfrag{m71}{{\tiny $m=7$, mesh $50$}}
\psfrag{m12}{{\tiny $m=1$, mesh $100$}}
\psfrag{m22}{{\tiny $m=2$,mesh $100$}}
\psfrag{m32}{{\tiny $m=3$, mesh $100$}}
\psfrag{m52}{{\tiny $m=5$,mesh $100$}}
\psfrag{m72}{{\tiny $m=7$, mesh $100$}}
\psfrag{m13}{{\tiny $m=1$, mesh $200$}}
\psfrag{m23}{{\tiny $m=2$,mesh $200$}}
\psfrag{m33}{{\tiny $m=3$, mesh $200$}}
\psfrag{m53}{{\tiny $m=5$,mesh $200$}}
\psfrag{m73}{{\tiny $m=7$, mesh $200$}}
	\hspace{-0.5cm}\includegraphics[width=0.45\linewidth]{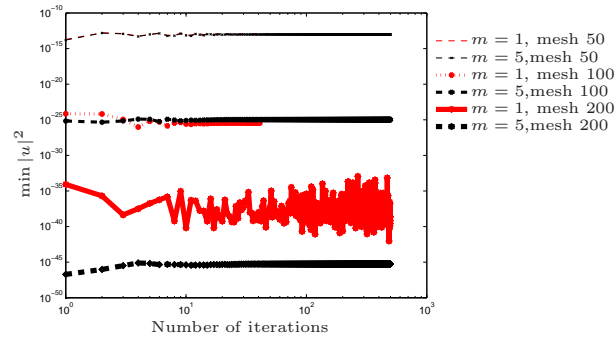}
\end{center}
 \caption{\label{mmin}  \textcolor{black}{Plot of $\min | u|^2$ versus the number of iterations for different values of $k_1=\pi \times 10^{m}$, $m=1,5$, where $u$ is the numeric solution of the Helmholtz problem.}}\label{mmin}
\end{figure}
\appendices
\section{\textcolor{black}{Proof of Lemma~\ref{P1}}}\label{App}
{\color{black}

Let $D(x)=F(u)f(xa)+G(u)(bx-1)$, where $a=\frac{\tilde{\gamma}}{\gamma}$ with $b=\frac{\tilde{q}}{q}$ are unknown and $f(x)=\frac{(x-1)^2}{x+1}$. Using the linearity
of the second term, and by  introducing $N(x)=~F(u)f(xa)-D(x)$, we see that
\begin{eqnarray*}
N(x) & = & \frac{N(x_{2})-N(x_{1})}{x_{2}-x_{1}}x+\frac{x_{2}N(x_{1})-x_{1}N(x_{2})}{x_{2}-x_{1}}\\
 & = & \frac{N(x_{1})(x_{2}-x)+N(x_{2})(x-x_{1})}{x_{2}-x_{1}}.
\end{eqnarray*}
 By returning to $D$, and by  introducing the function \[
d(x_{1},x_{2},x)=D(x)-\frac{D(x_{1})(x_{2}-x)+D(x_{2})(x-x_{1})}{x_{2}-x_{1}}\]
 we have \[
d(_{1},x_{2},x)=F(u)\left[f(xa)-\frac{f(ax_{1})(x_{2}-x_{1}) +f(ax_{2})(x-x_{1})}{x_2-x_1}\right],\]
 which is also \[
d(x_{1},x_{2},x)=F(u)\left[4a^2  \frac{x(x-x_{1}-x_{2})+x_1x_2}{a^{3}xx_{1}x_{2}+a^{2}(x(x_1+x_2)+x_1x_2)+a(x_{1}+x_{2}+x)+1}\right].\]
 Let us define\[
Q(x_{1},x_{2},x_{3},a)=4a^2  \frac{x_3(x_3-x_{1}-x_{2})+x_1x_2}{a^{3}x_3x_{1}x_{2}+a^{2}(x_3(x_1+x_2)+x_1x_2)+a(x_{1}+x_{2}+x_3)+1}.\]

We have obtained \[
d(x_{i},x_{j},x_{k})=F(u)Q(x_{i},x_{j},x_{k},a).\]
 Note that $Q(x_{i},x_{j},x_{k},a)=Q(x_{j},x_{i},x_{k},a)$, but other
permutation do not in general yield the same values.

As a consequence, from $n$ distinct measurements, we obtain $3C_{n}^{3}$
identities, that is, $3C_{n}^{3}$ formulas of the form 
\begin{equation}
\frac{1}{F(u)}=Q(x_{i},x_{j},x_{k},a)\frac{1}{d(x_{i},x_{j},x_{k})}.\label{10}
\end{equation}
 The value of $a$ can thus be deducted by intersection.

Note that $Q$, as a function of $a$, has only two roots equal to zero (for $x_3(x_3-x_{1}-x_{2})+x_1x_2 \ne 0$). By
an appropriate choice of $x_{i}$, $x_{j}$, $x_{k}$, we can set
 $a\in(0,\infty).$

We see that the equation becomes \[
Q(x_{i},x_{j},x_{k},a)=\frac{d(x_{i},x_{j},x_{k})}{d(x'_{1},x'_{j},x'_{k})}Q(x'_{1},x'_{j},x'_{k},a).\]
 Provided that the function $a\mapsto\frac{Q(x_{i},x_{j},x_{k},a)}{Q(x'_{i},x'_{j},x'_{k},a)}$
is bijective, $a$ is  determined uniquely. By  using relation~(\ref{10}),  this
defines $F$, and therefore $N$ and finally $G$.

Consequently, to determine  $a$ and $b$, it is sufficient to choose four different points $x_{1}$, $x_2$, $x_3$ and $x_4$ to obtain a bijective 
function on $(0,\infty)$ of the form
\[
a\mapsto\frac{Q(x_{1},x_{2},x_{3},a)}{Q(x_{1},x_{2},x_{4},a)}.\]
}

\label{lastpage}


\begin{thebibliography}{00}
\bibitem{AmFink}\textit{H. Ammari, E. Bonnetier, Y. Capdeboscq, M. Tanter, and M. Fink,}  
Electrical Impedance Tomography by Elastic Deformation
SIAM J. Appl. Math. , 68(6), (2008), 1557--1573.
\bibitem{AmKg}\textit{H. Ammari, H.Kang} Layer Potential Techniques in Imaging.Mathematical Surveys and Monographs, 153, Am. Math. Soc., Providence, 2009. 
\bibitem{yveV}\textit{Y. Capdeboscq and M. Vogelius,} A general
representation formula for boundary voltage perturbations caused by
internal conductivity inhomogeneities of low volume fraction. Math.
Modeling Num. Anal., 37 (2003), 159-173.
\bibitem{yveP}\textit{Y. Capdeboscq,} Private Communication. 2008.
\bibitem{AY}\textit{H. Ammari,Y. Capdeboscq, A. Rozanova-Pierrat,} Microwave imaging by elastic perturbations. in preparation.
\bibitem{Free}\textit{F. Hecht, O.Pironneau, K. Ohtsuka, A. Le Hyaric,} FreeFem++, http://
www.freefem.org/ (2007).
\bibitem{CGFK}\textit{H. Ammari, Y. Capdeboscq, H. Kang, A. Kozhemyak} Mathematical models and reconstruction methods in magneto-acoustic imaging. Preprint.
\end{thebibliography}
\end{document}